\documentclass[12pt,oneside,reqno]{amsart}

\usepackage[usenames,dvipsnames,svgnames,table]{xcolor}
\usepackage{mathtools}
\usepackage{amsfonts,amsthm,amsmath,amssymb}
\usepackage{enumitem}
\usepackage{hyperref}

\usepackage[paper=letterpaper, margin=2.5cm]{geometry}
\usepackage{valuediagrams}
\setkeys{valuediagram}{
    classic color=LimeGreen,
    quantum color=FireBrick,
    fontsize=\tiny,
    xscale=0.4,
    yscale=0.2,
    critical windows=None,
    critical window color=Yellow!50,
    critical window padding=0.1,
    node size=3pt,
    padding=0.5,
    xmin=1,
    xmax=None,
    height=None,
    thickness={very thick},
    arrow={-},
    baseline=0pt,
    tikzpicture options={},
}%

\usepackage{xargs}
\usepackage[colorinlistoftodos,prependcaption,textsize=tiny,linecolor=red,backgroundcolor=red!25,bordercolor=red]{todonotes}
\newcommandx{\franco}[2][1=]{\todo[linecolor=blue,backgroundcolor=blue!25,bordercolor=blue,#1]{#2 ---Franco}}
\newcommandx{\laura}[2][1=]{\todo[linecolor=orange,backgroundcolor=orange!25,bordercolor=orange,#1]{#2 ---Laura}}
\newcommandx{\caro}[2][1=]{\todo[linecolor=purple,backgroundcolor=purple!25,bordercolor=purple,#1]{#2 ---Caro}}
\newcommandx{\nantel}[2][1=]{\todo[linecolor=red,backgroundcolor=red!25,bordercolor=red,#1]{#2 ---Nantel}}
\newcommandx{\frank}[2][1=]{\todo[linecolor=green,backgroundcolor=green!35,bordercolor=green,#1]{#2 ---Frank}}

\newtheorem{theorem}{Theorem}[section]
\newtheorem{proposition}[theorem]{Proposition}
\newtheorem{corollary}[theorem]{Corollary}
\newtheorem{lemma}[theorem]{Lemma}

\theoremstyle{definition}
\newtheorem{example}[theorem]{Example}
\newtheorem{remark}[theorem]{Remark}

\numberwithin{equation}{section}

\newcommand{\Fl}{{\mathbb F}\ell}
\newcommand{\Fln}{{\mathbb F}\ell_n}
\newcommand{\CC}{{\mathbb C}}
\newcommand{\NN}{{\mathbb N}}
\newcommand{\ZZ}{{\mathbb Z}}

\newcommand{\frakS}{{\mathfrak S}}
\newcommand{\calF}{{\mathcal F}}
\newcommand{\calL}{{\mathcal L}}
\newcommand{\calM}{{\mathcal M}}
\newcommand{\qu}{{\mathbf u}}
\newcommand{\qv}{{\mathbf v}}
\newcommand{\bzero}{{\mathbf 0}}
\newcommand{\bx}{{\mathbf x}}
\newcommand{\precdot}{{\prec\!\!\!\cdot\,}}

\newcommand{\sg}{{\mbox{\rm sg}}}

\newcommand{\ptsp}{\hspace{0.55pt}}
\newcommand{\pb}{\raisebox{-1pt}{\ptsp\color{magenta}\rule{0.75pt}{11pt}}\ptsp}  
\newcommand{\p}{\raisebox{-0.5pt}{\ptsp\color{magenta}\rule{0.75pt}{9pt}}\ptsp}  
\newcommand{\pp}{\raisebox{-0.3pt}{\ptsp\color{magenta}\rule{0.75pt}{6pt}}\ptsp}  

\DeclareMathOperator{\het}{ht}
\DeclareMathOperator{\supp}{supp}
\DeclareMathOperator{\nd}{end}

\DeclareMathOperator{\kactionop}{{\scriptstyle \bullet}}
\newcommand{\kaction}[1][k]{\kactionop_{#1}}

\newcommand{\rc}[1]{{\color{red}#1}}
\newcommand{\defcolor}[1]{{\color{RoyalBlue}#1}}
\newcommand{\demph}[1]{\defcolor{{\sl #1}}}
\newcommand{\ncycle}{\mathfrak o}

\def\R{\mathbf{R}}
\def\C{\mathbf{C}}
\def\F{\mathbf{F}}
\def\T{\mathbf{T}}

\copyrightinfo{}{}

\title[Quantum Murnaghan--Nakayama rule]{A quantum Murnaghan--Nakayama rule\\ for the flag manifold}

\author[Benedetti]{C.~Benedetti}
\address{Carolina Benedetti, Universidad de Los Andes, Bogot\'a, Colombia}
\email{c.benedetti@uniandes.edu.co}

\author[Bergeron]{N.~Bergeron}
\address{Nantel Bergeron, York University, Toronto, Ontario, Canada}
\email{bergeron@yorku.ca}

\author[Colmenarejo]{L.~Colmenarejo}
\address{Laura Colmenarejo, North Carolina State University, Raleigh, North Carolina, USA}
\email{lcolmen@ncsu.edu}

\author[Saliola]{F.~Saliola}
\address{Franco~Saliola, Universit\'e du Qu\'ebec \`a Montr\'eal, Montr\'eal, Qu\'ebec, Canada}
\email{saliola.franco@uqam.ca}

\author[Sottile]{F.~Sottile}
\address{Frank Sottile, Texas A\&M University, College Station, Texas, USA}
\email{sottile@tamu.edu}

\begin{document}

\begin{abstract}
    In this paper, we give a rule for the multiplication of a Schubert class by a tautological
    class in the (small) quantum cohomology ring of the flag manifold.
    As an intermediate step, we establish a formula for the multiplication of
    a Schubert class by a quantum Schur polynomial indexed by a hook partition.
    This entails a detailed analysis of chains and intervals in the
    quantum Bruhat order.
    This analysis allows us to use results of Leung--Li  and of Postnikov  to reduce quantum
    products by hook Schur polynomials to the (known) classical product.
\end{abstract}

\maketitle

\section{Introduction}

The \demph{Murnaghan--Nakayama rule} is a method for computing values of the
irreducible characters of symmetric groups \cite{Mu37,Na41a,Na41b}.
Under the identification of the character ring with the ring of symmetric
functions, it translates into a rule for the multiplication of a Schur
symmetric function by a Newton power sum.
As explained in~\cite[Sect.~2]{MorrisonSottile}, in the cohomology ring of a Grassmannian, the
Murnaghan--Nakayama formula is a formula for the intersection product of a Schubert variety
with a component of the Chern character of the tautological bundle (a tautological class), as these tautological classes are
essentially Newton power sums.
Consequently, any formula describing multiplication by a tautological class will also be called a Murnaghan--Nakayama rule.

Murnaghan--Nakayama rules have been described for many different algebras
that carry a notion of a Newton power sum.
        Fomin and Green gave a version for noncommutative symmetric functions, which
        led to formulas for characters of representations associated to stable Schubert
        and Grothendieck polynomials~\cite{FoGr98}.
        Assaf and McNamara gave a skew version~\cite{A-McN}, which Konvalinka generalized to a skew
        rule for multiplication by a `quantum' (perturbed by a parameter $q$) power sum
        function~\cite{Ko12}.
        Bandlow, \emph{et al.} gave a version in the cohomology of an affine Grassmannian~\cite{BSZ11},
        and Lee gave a version in the cohomology of the affine flag manifold~\cite{Lee}.
       %
       %
        Tewari gave a version for noncommutative Schur functions~\cite{Tewari},
        and Berg \emph{et al.} for the immaculate noncommutative symmetric functions \cite{BBSSZ}.
        Wildon gave a plethystic version~\cite{Wildon}.
        Ross gave a version for loop Schur functions~\cite{Ro14}, providing
        a fundamental step in the orbifold Gromov--Witten/Donaldson--Thomas
        correspondence~\cite{RoZo13}.
        Morrison and Sottile gave Murnaghan--Nakayama rules in
        the cohomology ring of a flag manifold and in the quantum cohomology of the
        Grassmannian~\cite{MorrisonSottile}.
        Khanh, \emph{et al.} gave a version in the Grothendieck ring of the Grassmannian~\cite{KHST}.
        Khanh also gave a version for $K$-$k$-Schur and $k$-Schur functions~\cite{Khanh}, generalizing~\cite{BSZ11}.
        Fan, \emph{et al.}  gave a Murnaghan--Nakayama formula for Chern--Schwartz--MacPherson classes of
        Schubert cells in the flag variety~\cite{FGX}.
        %

We first derive a formula for the multiplication of a Schubert polynomial by a hook Schur polynomial and use this to give a
new proof of the 
Murnaghan--Nakayama formula in the cohomology of a flag variety that was originally derived in~\cite{MorrisonSottile}.
This hook multiplication formula is related to the formula of M{\'e}sz{\'a}ros et al.~\cite[Thm.~2]{MPP}, but it does not use
the noncommutative Fomin--Kirillov algebra~\cite{FominKirillov}.

Our main result (Corollary~\ref{C:quantumMN}) is a Murnaghan--Nakayama rule in the small quantum cohomology
ring of the flag manifold, proving a formula conjectured by Morrison at FPSAC
2014~\cite{Morrison}.
This follows from a  formula (Theorem~\ref{thm:q_hook_product}) for multiplication by a quantum Schur polynomial indexed by
a hook partition. 
M{\'e}sz{\'a}ros et al.~\cite[Thm.~15]{MPP} also gave a formula for multiplication by a hook quantum Schur
polynomial, which would imply our formula.
Unfortunately, their proof is not valid, as we explain in Section~\ref{SS:R_MPP}---this necessitates the arguments we give.
Our arguments involve a detailed analysis of chains and intervals in the \demph{quantum Bruhat order}, a partial order relevant to
multiplication in the quantum cohomology of a flag manifold.
We also use the ``quantum-equals-classical'' result of Leung and Li \cite{LeungLi2012} as well as
Postnikov's cyclic symmetry \cite{Postnikov_symmetry} to show that structure constants for quantum products by hook quantum
Schur 
polynomials are equal to structure constants for the classical product from~\cite{So96}, and thereby deduce our results.

In Section~\ref{sec:classic-mn-rule} we develop some background on the flag variety and Bruhat order, giving a new proof of the
Murnaghan--Nakayama formula in the cohomology of the flag variety.
Section~\ref{S:qCoh} develops background on the quantum cohomology ring of the flag variety, and gives a proof of the 
Murnaghan--Nakayama formula in this setting, which relies on results established in two subsequent
sections.
We extract some results of Leung and Li \cite{LeungLi2012} in Section~\ref{quantum-equals-classical}, including a version
of their  ``quantum-equals-classical'' result (Lemma~\ref{L:Quantum_is_Classical}).
The final Section~\ref{left-operators} develops a detailed combinatorics of chains in the quantum Bruhat order to prove the second
result used in our proof of the quantum Murnaghan--Nakayama rule (Theorem~\ref{Th:nonZero}).
It is the technical heart of this paper.

A theme in the paper~\cite{BS-Bruhat} is the relation between intervals in the $k$-Bruhat order $\leq_k$ (a suborder of the Bruhat order on
$S_n$) and Littlewood--Richardson coefficients $c^w_{u,v}$ when $v$ is a Grassmannian permutation.
For example, the isomorphism type of an interval $[u,w]_k$ in the $k$-Bruhat order only depends upon $wu^{-1}$, and $c^w_{u,v}$ also only
depends upon $wu^{-1}$ when $v$ is Grassmannian.
Only substantially weaker versions of these and related results hold for quantum cohomology.
While an interval $[u,q^\alpha w]_k^q$ in the quantum $k$-Bruhat order (see Section~\ref{quantum-k-bruhat-order}) is not determined by
$wu^{-1}$, its isomorphism type is preserved by some operations including a cyclic shift (Corollary~\ref{C:equivalences}) that is related to
Postnikov's cyclic shift and is distinct from the cyclic shift in~\cite{BS-Bruhat} (which does not preserve isomorphism type).
When $v$ is Grassmannian, we show that the numerical part $N^{w,\alpha}_{u,v}$ of a quantum Littlewood--Richardson coefficient (see
Section~\ref{quantum-equals-classical}) only depends upon $wu^{-1}$ (but $q^\alpha$ may change); this is
Theorem~\ref{Th:quantumIndependence}.

\begin{remark}\label{rem:fpsac}
The results of this paper were announced in a published FPSAC conference abstract~\cite{FPSAC23}. Here, we provide all the proofs, clarify some definitions and statements, and correct typos.
\end{remark}

\section{A Murnaghan--Nakayama rule for (classic) cohomology}
\label{sec:classic-mn-rule}

We review combinatorial models for the cohomology ring of the flag manifold.
We omit the precise definitions of the flag manifold and
its cohomology as the combinatorial models will suffice for our purposes; the
interested reader can consult~\cite{Fu97}.

Then we deduce a formula for the product of a Schubert polynomial and a hook
Schur function from which we deduce the Murnaghan--Nakayama rule in the
cohomology of the flag manifold, resulting in a different proof than that given by Morrison and Sottile \cite{MorrisonSottile}.
Our proof of the Murnaghan--Nakayama rule in quantum cohomology is analogous.

\subsection{Cohomology of the flag manifold}
\label{cohomology-ring-of-the-flag-manifold}

Let $n\in\NN$ be  a positive integer and write \defcolor{$\Fln$} for the manifold of complete flags in $\CC^n$.
We let \defcolor{$H^*\Fln$} denote its cohomology ring with coefficients in $\ZZ$.
This is a free $\ZZ$-module with a distinguished basis
\defcolor{$\{\frakS_w\}_{w \in S_n}$} of \demph{Schubert classes} which are indexed by elements of the symmetric group
\defcolor{$S_n$}. 
Borel~\cite{Borel} described the ring structure of $H^*\Fln$ as a quotient of a polynomial ring,
 \begin{equation}
    \label{borel-presentation}
    H^* \Fln\ \simeq\ \ZZ[x_1,\dotsc,x_n] / I_n\,,
    \qquad I_n\ =\ \big\langle e_a(x_1,\dotsc,x_n)\mid a \in [n] \big\rangle\,,
 \end{equation}
where $e_a$ is the $a$th elementary symmetric polynomial and
$\defcolor{[n]}\vcentcolon=\{1,\dots,n\}$.
Lascoux and Sch\"utzen\-ber\-ger~\cite{LaSch82} defined a set of polynomials in
$\ZZ[x_1, \ldots, x_n]$ called \demph{Schubert polynomials} that correspond to
the Schubert classes $\frakS_w$ under Borel's isomorphism.
We write $\frakS_w(x)$ for the Schubert polynomial corresponding to $\frakS_w$
on the set of variables $x: = \{x_1, \ldots, x_n\}$.

Every Schur polynomial $s_\lambda(x_1, \ldots, x_k)$ is a Schubert polynomial.
Let $w\in S_n$ be a permutation with a unique descent at position $k$.
Then the sequence $(w(k){-}k,\dotsc, w(2){-}2,w(1){-}1)$ is a partition contained in the \demph{rectangular
partition $R_{k,n-k}$} of $k(n{-}k)$ into $k$ parts with each part of
size $n{-}k$.
This is a bijection and
for $\lambda \subseteq R_{k, n-k}$ we write $\defcolor{v(\lambda,k)}\in S_n$ for
the corresponding permutation.
Then $\frakS_{v(\lambda, k)}(x) = s_\lambda(x_1, \ldots, x_k)$.
For example, when $n=7$,  
\begin{equation}
  v((3,1,0),3)\ =\ 136\p2457
  \qquad\mbox{and}\qquad
  v((1,0,0,0),4)\ =\ 1235\p467\ =\ (4,5)\,.
\end{equation}
Observe that we may write permutations using both one-line and cycle notation.
The ``\,$\p$\,'' is written after the $k$-th position in a permutation in one-line notation.
  
It remains an important open problem in Schubert calculus to
find a combinatorial rule expressing $\frakS_u \cdot \frakS_v$
in terms of the basis $\{ \frakS_w \}_{w \in S_n}$ of $H^*\Fln$.
That is, one seeks a combinatorial rule for the coefficients
$\defcolor{c_{u,v}^{w}}\in\ZZ$ in the product
\begin{equation}
    \label{eq:classical-coeffs}
    \frakS_u \cdot \frakS_v\ =\ \sum_{w} c_{u, v}^{w} \frakS_w\,.
\end{equation}
A remarkable property of Schubert polynomials is that this computation can
be performed directly in the polynomial ring as the resulting coefficients
$c_{u, v}^w$ are the same.

Although this problem remains open in general, several special cases are known.
\demph{Monk's Formula}~\cite{Monk} treats the case where $v$ is the simple transposition $(k, k{+}1)$.
For any $u \in S_n$,
  \begin{equation}
    \label{classical-Monk-formula}
    \frakS_u(x) \cdot \frakS_{(k,k+1)}(x)
       \ =\ 
    \sum_{\substack{1 \leq i \leq k < j \leq n \\ \ell(u) + 1 = \ell(u(i,j))}} \frakS_{u (i,j)}(x)\,,
 \end{equation}
where $\defcolor{\ell(u)}\vcentcolon=\#\{i<j\mid u(i)>u(j)\}$ denotes the length of the permutation $u$.
In this case, $\frakS_{(k, k+1)}(x) = s_{(1)}(x_1, \ldots, x_k)$.
Proposition~\ref{pieri-formulas-version-2} below describes the special case of
multiplying by $s_{(b, 1^{a-1})}(x_1, \ldots, x_k)$
in terms of the \emph{$k$-Bruhat order}, which we define next.
Here, $(b,1^{a-1})$ is the \demph{hook partition} of $a{+}b{-}1$ into $a$ parts, the first of size $b$, and the remaining all of size $1$.

\subsection{$k$-Bruhat order}

The terms appearing in Monk's Formula \eqref{classical-Monk-formula}
define a partial order on $S_n$ graded by $\ell(u)$.
Define the \demph{$k$-Bruhat order $\leq_k$} by the covering relation
\[
    {\defcolor{u\lessdot_k u (i,j)}}
    \qquad\text{if}\quad
    1 \leq i\leq k<j \leq n
    \quad\text{and}\quad
    \ell(u)+1=\ell(u(i,j))\,.
\]
That is, $u\lessdot_k u (i,j)$ if and only if $\frakS_{u(i,j)}(x)$ appears in the product
$\frakS_{u}(x) \cdot \frakS_{(k,k+1)}(x)$.
Since $\ell(u)+1=\ell(u (i,j))$ if and only if $u(i)<u(j)$ and there does not
exist $l$ such that $i < l < j$ and $u(i) < u(l) < u(j)$, we obtain the
following equivalence.

\begin{proposition}
    \label{nothing-in-between-is-in-between}
    Let $u$ and $w$ be two permutations.
    Then $u \lessdot_k w$ if and only if
    \begin{enumerate} 
        \item
            $w = u(i, j)$ for some $i \leq k < j$ such that $u(i) < u(j)$;
            and

        \item
            there does not exist $l$ such that
            $i < l < j$ and $u(i) < u(l) < u(j)$.
    \end{enumerate}
\end{proposition}

Write \defcolor{$[u,w]_k$} for the interval between $u$ and $w$ in the
$k$-Bruhat order.
We may label each cover relation $u\lessdot_k u(i,j)$ with $i\leq k<j$ by $u(i)$ and write
\begin{equation*}
    \defcolor{u \xrightarrow{u(i)} u(i,j)}
    \qquad \text{if} \qquad
    u \lessdot_k u(i,j)\ \mbox{with}\ i\leq k<j\,.
\end{equation*}
A saturated chain 
$\gamma\colon u\xrightarrow{l_1}u_1\xrightarrow{l_2}\dotsb\xrightarrow{l_r}u_r$
in this poset has \demph{length~$r$}
and we denote its maximal element by $\defcolor{\nd(\gamma)} = u_r$.
If in addition there exists a number $a$ such that $l_1>\dotsb >l_a<l_{a+1}<\dotsb<l_{r}$, then we say that $\gamma$ is
a \demph{peakless chain} of \demph{height $a$} and \demph{length $r$}.

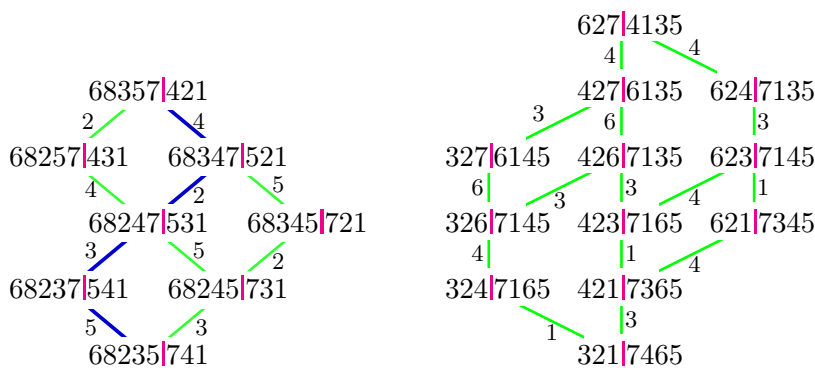
\begin{figure}[htb]
  \centering
   \begin{picture}(135,115)(-25,-2)
    \put(-4.5,0){\includegraphics{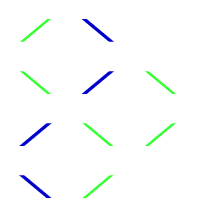}}
                 \put(5.5,100){\small$68357\p421$}
    \put(-24.5, 75){\small$68257\p431$}     \put(35.5, 75){\small$68347\p521$} 
                 \put(5.5, 50){\small$68247\p531$}     \put(65.5,50){\small$68345\p721$} 
    \put(-24.5, 25){\small$68237\p541$}     \put(35.5,25){\small$68245\p731$} 
                 \put(5.5,  0){\small$68235\p741$}
    \put( 3, 89){\scriptsize$2$} \put( 45, 89){\scriptsize$4$} 
    \put( 4, 63){\scriptsize$4$} \put( 45, 61){\scriptsize$2$} \put( 75, 65){\scriptsize$5$}
    \put( 4, 40){\scriptsize$3$} \put( 45, 40){\scriptsize$5$} \put( 75, 36){\scriptsize$2$}
    \put( 4, 11){\scriptsize$5$} \put( 46, 11){\scriptsize$3$}
   \end{picture}  
    \qquad
   \begin{picture}(144,135)(-53,-2)
    \put(-33.5,0){\includegraphics{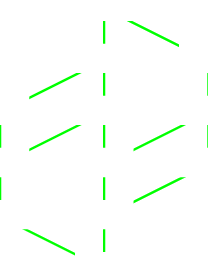}}
                                      \put(  0,125){\small$627\p4135$}
                                      \put(  0,100){\small$427\p6135$}   \put(50,100){\small$624\p7135$} 
    \put(-50, 75){\small$327\p6145$}  \put(  0, 75){\small$426\p7135$}   \put(50, 75){\small$623\p7145$}
    \put(-50, 50){\small$326\p7145$}  \put(  0, 50){\small$423\p7165$}   \put(50, 50){\small$621\p7345$}
    \put(-50, 25){\small$324\p7165$}  \put(  0, 25){\small$421\p7365$}
                                      \put(  0,  0){\small$321\p7465$}

            \put( 10, 114){\scriptsize$4$} \put( 42, 117){\scriptsize$4$}
    \put( -17, 92){\scriptsize$3$} \put( 10, 89){\scriptsize$6$} \put( 68, 89){\scriptsize$3$}
                                   \put( 18, 64){\scriptsize$3$} \put( 42, 60){\scriptsize$4$} \put( 68, 64){\scriptsize$1$}
    \put( -40, 64){\scriptsize$6$} \put( -9, 59){\scriptsize$3$}
    \put( -40, 39){\scriptsize$4$} \put( 18, 39){\scriptsize$1$} \put( 42, 35){\scriptsize$4$}
    \put( -12,  9){\scriptsize$1$} \put( 18, 14){\scriptsize$3$}
   \end{picture} 

\caption{Two intervals in $k$-Bruhat orders on $S_n$.}
   \label{F:k-Bruhat} 
\end{figure}

Figure~\ref{F:k-Bruhat} shows two labeled intervals,
one in the $5$-Bruhat order on $S_8$ and the other in the $3$-Bruhat order on $S_7$.
The interval on the left has a peakless chain (colored in blue); its `middle' chain has labels $5,3,2,4$, so it is peakless of height 3,
but the other interval has no peakless chains.
\subsection{Grassmannian--Bruhat order}\label{SS:GBO}

We will make use of a second partial order on $S_n$ that is defined using the
$k$-Bruhat order. Define the \demph{Grassmannian--Bruhat order} on $S_n$ by
\begin{equation*}
    \defcolor{\eta\preceq\zeta}
    \qquad\text{if}\quad
    \begin{array}{c}
        \text{there exists $u \in S_n$ and $1\leq k < n$ with $u\leq_k \eta u\leq_k\zeta u$.}
    \end{array}
\end{equation*}
This is a graded partial order on $S_n$ with minimal element the identity permutation $e$.
Its rank function is $\defcolor{\calL(\zeta)}\vcentcolon=\ell(\zeta u)-\ell(u)$ for any $u \in S_n$ and $k \in \NN$ with
$u\leq_k\zeta u$. 
The cover relations inherit a labeling from the $k$-Bruhat orders.
As shown in~\cite[Sect.~3.2]{BS-Bruhat}, neither $\preceq$, $\calL$, nor the labeling depend on the choice of $k$ or $u$.
For a permutation $\zeta\in S_n$, its \demph{height} is 
$\defcolor{\het(\zeta)} \vcentcolon= \#\{i\in[n]\mid i<\zeta(i)\}$ and its  
\demph{support} is $\defcolor{\supp(\zeta)} = \{i\in[n]\mid\zeta(i)\neq i\}$.
Then for any $\zeta \in S_n$ and $k$ satisfying $\het(\zeta)\leq k\leq n-\#\supp(\zeta)+\het(\zeta)$,
there is a $u\in S_n$ with $u\leq_k \zeta u$~\cite[Thm.~3.1.5]{BS-Bruhat},
which among other things implies that $e \preceq \zeta$.

The order $\preceq$ enjoys another symmetry.
Permutations $\zeta\in S_n$ and $\eta\in S_m$ are \demph{shape-equivalent} if we have 
$\supp(\zeta)\subseteq\{i_1<\dotsb<i_r\}$ and $\supp(\eta)\subseteq\{j_1<\dotsb<j_r\}$ and for all $a,b$, $\zeta(i_a)=i_b$
if and only if $\eta(j_a)=j_b$.
For example, the cycles  $(2,3,5,7,4)$ and $(1,2,4,5,3)$ are shape-equivalent.
Shape-equivalence 
induces an isomorphism on intervals~\cite[Thm.~3.2.3({\it v})]{BS-Bruhat}.
It also preserves the relative order of labels on chains,
so that a chain in one interval is peakless if and only if the corresponding chain in a shape-equivalent interval is peakless.

\begin{figure}[htb]
  \centering
   \begin{picture}(96,115)(-6,-2)
    \put(  0,-1){\includegraphics{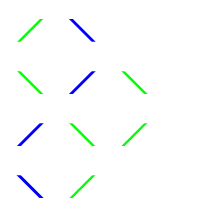}}
    \put(  9,100){\small$(23574)$}
    \put(-13, 75){\small$(3574)$}      \put(30, 75){\small$(234)(57)$} 
    \put(  7.5, 50){\small$(34)(57)$}  \put(65,50){\small$(234)$} 
    \put( -8, 25){\small$(57)$}        \put(42,25){\small$(34)$} 
    \put( 25,  0){\small$e$} 
    \put( 4, 90){\scriptsize$2$} \put( 44, 90){\scriptsize$4$} 
    \put( 4, 62){\scriptsize$4$} \put( 44, 62){\scriptsize$2$} \put( 70, 63){\scriptsize$5$}
    \put( 5, 39){\scriptsize$3$} \put( 44, 39){\scriptsize$5$} \put( 70, 36){\scriptsize$2$}
    \put( 5, 10){\scriptsize$5$} \put( 44, 10){\scriptsize$3$}

   \end{picture}
    \qquad
   \begin{picture}(130,135)(-46,-2)
    \put(-33.5,0){\includegraphics{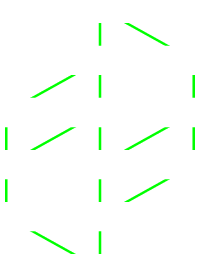}}
                                   \put(-7,125){\small$(174)(36)$}
                                   \put(-5,100){\small$(17634)$}  \put(42,100){\small$(14)(36)$} 
    \put(-45.5, 75){\small$(1764)$}  \put(-0.5, 75){\small$(1634)$}    \put(44.5, 75){\small$(1364)$}
    \put(-43, 50){\small$(164)$}   \put( 2, 50){\small$(134)$}     \put(47, 50){\small$(364)$}
    \put(-40, 25){\small$(14)$}    \put( 5, 25){\small$(34)$}
                            \put(14,  0){\small$e$}
    
            \put( 9, 114){\scriptsize$4$} \put( 40, 117){\scriptsize$4$}
    \put( -17, 92){\scriptsize$3$} \put( 18, 89){\scriptsize$6$} \put( 64, 89){\scriptsize$3$}
                                   \put( 18, 64){\scriptsize$3$} \put( 41, 60){\scriptsize$4$} \put( 64, 64){\scriptsize$1$}
    \put( -39, 64){\scriptsize$6$} \put( -8, 59){\scriptsize$3$}
    \put( -39, 39){\scriptsize$4$} \put( 18, 39){\scriptsize$1$} \put( 41, 35){\scriptsize$4$}
    \put( -13,  9){\scriptsize$1$} \put( 18, 13){\scriptsize$3$}
   \end{picture} 
    \qquad
   \begin{picture}(96,115)(-6,-2)
    \put(  0, -1){\includegraphics{minimal4_narrow.eps}}
    \put(  9,100){\small$(12453)$}
    \put(-13, 75){\small$(2453)$}      \put(30,75){\small$(123)(45)$} 
    \put(  7.5, 50){\small$(23)(45)$}  \put(65,50){\small$(123)$} 
    \put( -8, 25){\small$(45)$}        \put(42,25){\small$(23)$} 
    \put( 25,  0){\small$e$} 

    \put( 4, 90){\scriptsize$1$} \put( 44, 90){\scriptsize$3$} 
    \put( 4, 62){\scriptsize$3$} \put( 44, 61){\scriptsize$1$} \put( 70, 63){\scriptsize$4$}
    \put( 5, 39){\scriptsize$2$} \put( 44, 39){\scriptsize$4$} \put( 70, 36){\scriptsize$1$}
    \put( 5, 10){\scriptsize$4$} \put( 44, 10){\scriptsize$2$}
  \end{picture}
\caption{Intervals in the Grassmannian--Bruhat order.}
   \label{F:GBorder} 
\end{figure}

We illustrate this.
The intervals in Figure~\ref{F:k-Bruhat} correspond to the permutations $\zeta_1=(2,3,5,7,4)$ and $\zeta_2=(1,7,4)(3,6)$
under the map that associates the permutation $v u^{-1}$ with the interval $[u, v]$.
Figure~\ref{F:GBorder} shows labeled intervals $[e,\zeta_i]_\preceq$ for these two permutations, as well as for $\eta=(1,2,4,5,3)$, which is
shape-equivalent to $\zeta_1$.
The commas are omitted in the cycles for display purposes.
The peakless chains in the intervals $[e,\zeta_1]_\preceq$ and  $[e,\eta]_\preceq$ are highlighted.

\subsection{Minimal permutations}
\label{minimal-permutations} 

A permutation $\zeta$ is
\demph{minimal} if $\calL(\zeta) = \#\supp(\zeta) - s(\zeta)$,
where $\defcolor{s(\zeta)}$ is the number of (nontrivial) cycles in the factorization of $\zeta$ into disjoint cycles.
Minimality is a central concept in this paper.
This terminology (introduced in~\cite{BS-Bruhat}) arises from the inequality $\calL(\zeta) \geq \#\supp(\zeta) - s(\zeta)$, so that a
minimal permutation has the minimal possible rank given its cycle structure.
For example, $\zeta_1 = (2,3,5,7,4)$ is minimal, because
\begin{equation*}
    \calL(\zeta_1) = \ell(68357421) - \ell(68235741) = 22 - 18 = 4 = \#\supp(\zeta_1) - s(\zeta_1);
\end{equation*}
and $\zeta_2 = (1,7,4)(3,6)$ is not minimal, because
\begin{equation*}
    \calL(\zeta_2) = \ell(6274135) - \ell(3217465) = 12 - 7 = 5 > 3 = \#\supp(\zeta_2) - s(\zeta_2).
\end{equation*}

Related to minimality is the notion of crossing.
Let $l_1<l_2$ and $m_1<m_2$ be two disjoint pairs of numbers, and suppose that $l_1<m_1$.
There are three possibilities for $l_2$.
Either $l_2<m_1$ or $m_1<l_2<m_2$ or $m_2<l_2$.
When the second occurs, we say that the pairs are \demph{crossing}, otherwise they are \demph{noncrossing}.
This refers to arcs connecting elements of each pair, when the numbers are marked along a number line.
\[
   \begin{picture}(110,55)(0,-24)
    \put(0,0){\includegraphics{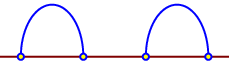}}
    \put(6,-9){$l_1$}  \put(36,-9){$l_2$}  \put(66,-9){$m_1$}  \put(96,-9){$m_2$}
    \put(24,-24){noncrossing}
   \end{picture}
   \qquad
   \begin{picture}(110,55)(0,-24)
    \put(0,0){\includegraphics{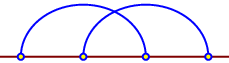}}
    \put(6,-9){$l_1$}  \put(36,-9){$m_1$}  \put(66,-9){$l_2$}  \put(96,-9){$m_2$}
    \put(34,-24){crossing}
   \end{picture}
   \qquad
   \begin{picture}(110,55)(0,-24)
    \put(0,0){\includegraphics{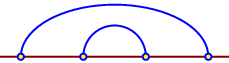}}
    \put(6,-9){$l_2$}  \put(36,-9){$m_1$}  \put(66,-9){$m_2$}  \put(96,-9){$l_2$}
    \put(24,-24){noncrossing}
   \end{picture}
\]
Let $A,B\subset[n]$ be disjoint and suppose that each has at least two elements.
If there exist a pair of elements of $A$ and a pair of elements of $B$ that are crossing, then $A$ and $B$ are
\demph{crossing}; otherwise they are  \demph{noncrossing}.

A product $\zeta\cdot\eta$ of two permutations is \demph{noncrossing}
when $\supp(\zeta)$ and $\supp(\eta)$ are disjoint and noncrossing.
A permutation $\zeta$ is \demph{irreducible} if it cannot be written
nontrivially as a noncrossing product.
Any permutation factors uniquely as a noncrossing product of irreducible
permutations.
In~\cite{BS-Bruhat}, $s(\zeta)$ counted the number of irreducible factors in the
noncrossing factorization of $\zeta$ for minimal permutations.
By \cite[Lem.~3.3.1]{BS-Bruhat} we have that
$s(\zeta)$ is the number of cycles (this is also explained in the proof below).
We give a characterization of minimal permutations in terms of
peakless chains in the Grassmannian--Bruhat order.

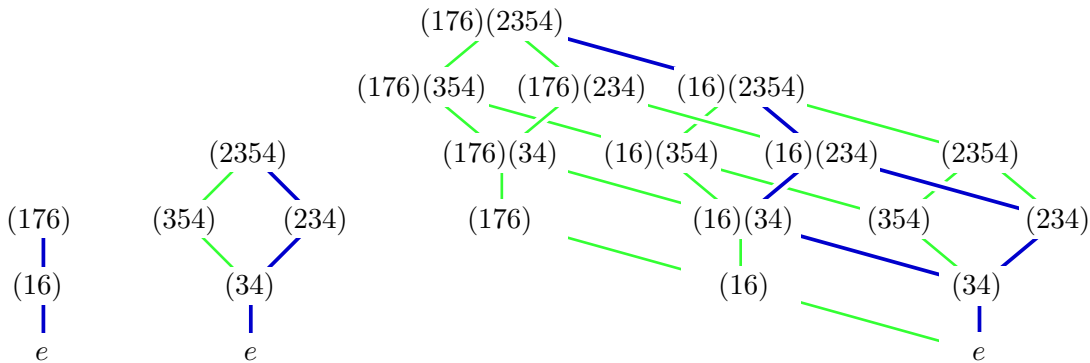
\begin{figure}[htb]
 \centering
 \begin{picture}(26,80)(-11,0)
  \put(0,0){\includegraphics{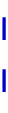}}  
  \put(-12,50){\small$(176)$}  \put(-10,25){\small$(16)$}   \put(-1,0){\small$e$}
 \end{picture}
 \qquad
 \begin{picture}(74,85)(-1,0)     
  \put(10,0){\includegraphics{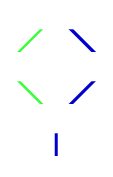}}
             \put(21,75){\small$(2354)$}
  \put(-0.5,50){\small$(354)$} \put(49,50){\small$(234)$}
             \put(27,25){\small$(34)$}
             \put(34.5,0){\small$e$}
 \end{picture}
 \begin{picture}(275,135)(-13,0)
  \put(-1,0){\includegraphics{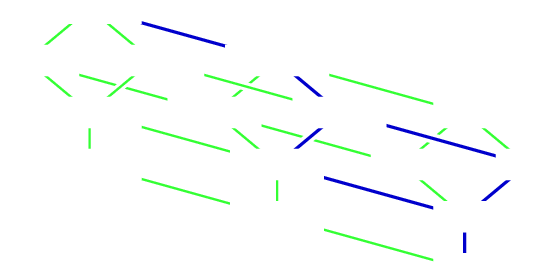}}   
             \put(11,125){\small$(176)(2354)$}
  \put(-13,100){\small$(176)(354)$} \put(47.5,100){\small$(176)(234)$}
             \put(19.5,75){\small$(176)(34)$}
             \put(29,50){\small$(176)$}

             \put(108,100){\small$(16)(2354)$}
  \put(80.5,75){\small$(16)(354)$} \put(141,75){\small$(16)(234)$}
             \put(114,50){\small$(16)(34)$}
             \put(124,25){\small$(16)$}

             \put(208,75){\small$(2354)$}
  \put(180,50){\small$(354)$} \put(240,50){\small$(234)$}
             \put(212,25){\small$(34)$}
             \put(220,0){\small$e$}
 \end{picture}

 \caption{If $\zeta\cdot\eta$ is a noncrossing product, then the interval
     $[e,\zeta\cdot\eta]_\preceq$ is isomorphic to the product of the intervals
     $[e,\zeta]_\preceq$ and $[e,\eta]_\preceq$, and the chains in
     $[e,\zeta\cdot\eta]_\preceq$ are shuffles of chains in each of the
 factors.}
 \label{F:product}
\end{figure}

\begin{lemma}\label{L:properties-noncrossing-products}
    Let $\zeta$ be a permutation and $\zeta=\zeta_1\dotsb \zeta_s$ its unique
    factorization into a noncrossing product of irreducible permutations.
    \begin{enumerate}
        \item
            $\calL(\zeta)= \calL(\zeta_1)+\dotsb+\calL(\zeta_s).$

        \item
            $[e,\zeta]_\preceq\ \simeq\  [e,\zeta_1]_\preceq\ \times\ \dotsb\ \times\  [e,\zeta_s]_\preceq.$

        \item
            $[e,\zeta]_\preceq$ has a peakless chain if and only if each interval $[e,\zeta_i]_\preceq$ has a peakless chain.

        \item
            $\zeta$ is a minimal permutation if and only if each $\zeta_i$ is a minimal permutation.

        \item
            If $\zeta$ is minimal and irreducible, then $\zeta$ is a cycle.
    \end{enumerate}
\end{lemma}

\begin{proof}
    (1) By \cite[Lem.~3.3.1]{BS-Bruhat}, $\calL(\zeta)= \calL(\zeta_1)+\dotsb+\calL(\zeta_s)$ as this is a noncrossing product.

    (2) \cite[Lem.~3.3.3]{BS-Bruhat} implies that $[e,\zeta]_\preceq\ \simeq\  [e,\zeta_1]_\preceq\ \times\ \dotsb\ \times\  [e,\zeta_s]_\preceq$.

    (3) If each interval $[e, \zeta_i]_{\preceq}$ has a peakless chain
        $\gamma_i$, then shuffling these chains appropriately gives a peakless
        chain $\gamma$ in $[e,\zeta]_\preceq$ by (2). Conversely, if $\gamma$ is
        a peakless chain in  $[e,\zeta]_\preceq$, restricting it to each
        $[e,\zeta_i]_\preceq$ gives a peakless chain $\gamma_i$.
        Figure~\ref{F:product} illustrates this
        for the noncrossing product
        of the two minimal permutations $\zeta=(1,7,6)$ and $\eta=(2,3,5,4)$. The
        peakless chains are indicated, and commas in cycles are omitted for display
        purposes.

    (4) follows from (1) and (2).

    (5) Suppose that $\zeta$ is minimal and irreducible.
        If $\zeta$ is not a single cycle, then $\zeta=\eta\cdot \eta'$,
        where $\eta$ and $\eta'$ have disjoint supports that are crossing.
        By \cite[Lem.~3.3.1]{BS-Bruhat}, $\calL(\zeta)> \calL(\eta) + \calL(\eta')$, which contradicts that $\zeta$ is minimal.
        This inequality is illustrated by $\zeta=(1,7,4)(3,6)$.
        We showed that $\calL(\zeta)=5$, but $\calL(1,7,4) = 2$ and $\calL(3,6)=1$.
\end{proof}

\begin{proposition}\label{P:peaklessChains}
 A permutation $\zeta$ is minimal if and only if the interval $[e,\zeta]_\preceq$ has a peakless chain.
 The number of peakless chains of height $a$ in $[e,\zeta]_\preceq$ is the binomial coefficient
 $\left(\begin{smallmatrix} s(\zeta)-1 \\ \het(\zeta)-a \end{smallmatrix}\right)$.
\end{proposition}
\begin{proof}
    We begin by showing that the first statement follows from the corresponding
    result for irreducible permutations.
    Let $\zeta=\zeta_1\dotsb \zeta_s$ be the unique
    factorization of $\zeta$ into a noncrossing product of irreducible permutations.
    By Lemma~\ref{L:properties-noncrossing-products},
    $\zeta$ is minimal if and only if each $\zeta_i$ is minimal.
    Since $\zeta_i$ is irreducible,
    $\zeta$ is minimal if and only if each interval $[e, \zeta_i]_{\preceq}$ has
    a peakless chain.
    By Lemma~\ref{L:properties-noncrossing-products},
    this holds if and only if $[e, \zeta]_{\preceq}$ has a peakless chain.

We now prove the result for irreducible permutations.
Suppose that $\zeta$ is minimal and irreducible.
By Lemma~\ref{L:properties-noncrossing-products}, $\zeta$ is a single cycle.
In~\cite[Lem.~6.7]{BS-iso} it is proved that if $\zeta$ is a minimal cycle,
then there is a unique peakless chain in the interval $[e, \zeta]_\preceq$,
which proves the forward direction.
Conversely, suppose that $[e, \zeta]_\preceq$ has a peakless chain.
Let \defcolor{$a$} be the height of that chain and set
$\defcolor{b}\vcentcolon=\calL(\zeta){-}a{+}1$.
Replacing $\zeta$ by a shape-equivalent permutation if necessary, we may assume
that $\supp(\zeta)=[n]$ (so that $\zeta\in S_n$) and thus $\calL(\zeta) \geq
\#\supp(\zeta) - s(\zeta) = n-1$.

Since $e \preceq \zeta$, there is a $k < n$ and $u \in S_n$ such that $u \leq_k \zeta u$ and $[e, \zeta]_\preceq \simeq [u, \zeta u]_k$.
By~\cite[Thm.~3.1.5]{BS-Bruhat}, we have
\begin{equation*}
    \{i \mid i < \zeta(i)\}\ \subseteq\ \{u^{-1}(1), \ldots, u^{-1}(k)\}
    \quad\text{and}\quad
    \{i \mid i > \zeta(i)\}\ \subseteq\ \{u^{-1}(k + 1), \ldots, u^{-1}(n)\}
\end{equation*}
so that
\begin{equation*}
    \het(\zeta)\ =\ \#\{i \mid i < \zeta(i)\} \leq k
    \quad\text{and}\quad
    n - \het(\zeta)\ =\ \#\{i \mid i > \zeta(i)\} \leq n - k.
\end{equation*}
Together, these imply that $\het(\zeta) = k$.
Since there is a peakless chain in $[e,\zeta]_\preceq$
and $[e,\zeta]_\preceq\simeq[u,\zeta u]_k$,
it follows from Formula~\eqref{eq:schubert-times-hook-via-peakless-chains}
below that the coefficient $c^{\zeta u}_{u, v((b,1^{a-1}),k)}$ is nonzero.
This implies that $v((b,1^{a-1}),k)\leq_k \zeta u$.
As $\zeta u\in S_n$, we have $v((b,1^{a-1}),k)\in S_n$.
Thus $a\leq k<n$ and $b\leq n{-}k$, so that $\calL(\zeta)=a{+}b{-}1\leq n{-}1$.
This implies that $\calL(\zeta)=n{-}1$, $a=k$, and $b=n{-}k$.
As $n=\#\supp(\zeta)$, we conclude that $\zeta$ is a minimal cycle
and that $[e, \zeta]_\preceq$ has a unique peakless chain.
Note that since $s(\zeta) = 1$ and $\het(\zeta) = a$, we have
$\left(\begin{smallmatrix} s(\zeta)-1 \\ \het(\zeta)-a \end{smallmatrix}\right)
= \left(\begin{smallmatrix} 0 \\ 0 \end{smallmatrix}\right) = 1$.

Finally, suppose that $\zeta=\zeta_1\dotsb\zeta_s$ is the noncrossing factorization of a minimal permutation $\zeta$ into disjoint minimal
cycles $\zeta_1,\dotsc,\zeta_s$.
Recall from Lemma~\ref{L:properties-noncrossing-products} that
$[e,\zeta]_\preceq \simeq  [e,\zeta_1]_\preceq \times \dotsb \times  [e,\zeta_s]_\preceq.$
We count the shuffles of peakless chains $\gamma_i$ in $[e,\zeta_i]_\preceq$ with a given height.
We need only consider the sequence of labels, and note that as the cycles $\zeta_i$ have disjoint supports, the set of labels on
the chains $\gamma_i$ are disjoint.
Partition a peakless chain (represented by its labels)
\begin{equation}\label{Eq:Peakless_Partition}
    l_1\ >\ \dotsb\ >\ l_{a-1} \ >\  l_a\ <\ l_{a+1}\ <\ \dotsb\ <\ l_{a+b-1}
\end{equation}
into three subsets.
Call $l_1> \dotsb>l_{a-1}$ the decreasing segment, $l_a$ the minimum, and
$l_{a+1} <\dotsb<l_{a+b-1}$ the increasing segment.

In a peakless chain~\eqref{Eq:Peakless_Partition} of height $a$ in $[e,\zeta]_\preceq$, its decreasing segment contains
the decreasing segments of each $\gamma_i$, and the same for increasing segments.
Its minimum $l_a$ is the minimum of all the minima of the $\gamma_i$.
The only choices possible are distributing the remaining $s{-}1$ minima between the decreasing and increasing segments
of~\eqref{Eq:Peakless_Partition}, which may be done independently.
Since $\het(\zeta)=\het(\zeta_1)+\dotsb+\het(\zeta_s)$, we see that placing $\het(\zeta)-a$ of these remaining minima in
the increasing segment gives a peakless chain of height $a$.
This completes the proof.
\end{proof}

\subsection{A Murnaghan--Nakayama rule for $H^*\Fln$}

We recall the following special case of
Formula~\eqref{eq:classical-coeffs}, which is~\cite[Cor. 9]{So96}.

\begin{proposition}
    \label{pieri-formulas-version-2}
    Let $u \in S_n$,  $a\leq k$, and $b\leq n{-}k$. Then
    \begin{equation}
        \label{eq:schubert-times-hook-via-peakless-chains}
        \frakS_u(x)\cdot s_{(b,1^{a-1})}(x_1,\dotsc,x_k)\ =\ \sum_\gamma \frakS_{\nd(\gamma)}(x)\,,
    \end{equation}
    the sum over all peakless chains starting at $u$ of height $a$ and length $a{+}b{-}1$ in the $k$-Bruhat order.
\end{proposition}

Using Proposition~\ref{P:peaklessChains}, we can explicitly describe the coefficients
appearing in \eqref{eq:schubert-times-hook-via-peakless-chains}.

\begin{proposition}\label{P:classical_hook}
    Let $u \in S_n$, $a\leq k$ and $b\leq n{-}k$.
    Then
    \begin{equation}
        \label{Eq:Schubert-times-hook}
        \frakS_u(x) \cdot s_{(b, 1^{a-1})}(x_1,\dotsc,x_k)\ =\
        \sum
        \left(\begin{smallmatrix} s(\zeta)-1 \\ \het(\zeta)-a \end{smallmatrix}\right)
        \frakS_{\zeta u}(x) \,,
    \end{equation}
    the sum over all minimal permutations $\zeta \in S_n$ such that
    $u\leq_k\zeta u$ and $\calL(\zeta)=a+b-1$.
\end{proposition}

We use Proposition~\ref{P:classical_hook} to derive a formula for the product of a Schubert polynomial and
a power sum symmetric polynomial.
Let $\defcolor{p_r(x_1,\dotsc,x_k)}\vcentcolon=x_1^r+\dotsb+x_k^r$ denote the $r$th power sum polynomial.
Expanded in terms of Schur polynomials~\cite{Macdonald}, we have
 \begin{equation}\label{Eq:PowerSumHook}
    p_r(x_1, \ldots, x_k) =
    s_{(r)}(x_1, \ldots, x_k) - s_{(r-1,1)}(x_1, \ldots, x_k) + \cdots + (-1)^{r+1} s_{(1^r)}(x_1, \ldots, x_k).
\end{equation}

\begin{corollary}[Murnaghan--Nakayama rule for $H^*\Fln$]
    \label{C:classical-schubert-times-powersum}
    Let $u \in S_n$.
    Then
     \begin{equation*}
      \frakS_u(x) \cdot p_r(x_1,\dotsc,x_k)\ =\ \sum (-1)^{\het(\zeta)+1} \frakS_{\zeta u}(x)\,,
     \end{equation*}
    the sum over all minimal cycles $\zeta \in S_n$ such that
    $u\leq_k\zeta u$ and $\calL(\zeta)=r$.
\end{corollary}
\begin{proof}
    Expanding the  product $\frakS_u(x) \cdot p_r(x_1,\dotsc,x_k)$ using Equations
    \eqref{Eq:Schubert-times-hook}
    and
    \eqref{Eq:PowerSumHook}
    shows
    that
    the coefficient of $\frakS_{\zeta u}$ for $u\leq_k\zeta u$ is zero unless
    $\calL(\zeta)=r$ and $\zeta$ is minimal.
    When $\zeta$ is the noncrossing product of $s(\zeta)$ minimal cycles, this coefficient is
    \begin{equation*}
        (-1)^{\het(\zeta)+1} \sum_{a=1}^{r} (-1)^{a-\het(\zeta)}\binom{s(\zeta)-1}{\het(\zeta)-a}\ =\
        \begin{cases}
            0,  & \text{ if } s(\zeta)\neq 1,\\
            (-1)^{\het(\zeta)+1},& \text{ if } s(\zeta)=1.
        \end{cases}
        \qedhere
    \end{equation*}
\end{proof}



This proof differs from that given in~\cite{Morrison, MorrisonSottile}.

\section{A Murnaghan--Nakayama rule for quantum cohomology}
\label{S:qCoh}

The (small) quantum cohomology ring \defcolor{$qH^*\Fln$} of the flag manifold
is a deformation of the ordinary cohomology ring $H^*\Fln$ whose product encodes the three-point
Gromov--Witten invariants of $\Fln$.
We describe this ring and extend the main
objects from Section~\ref{sec:classic-mn-rule} to the quantum setting.

\subsection{Quantum cohomology ring of the flag manifold}
As a $\ZZ$-module, we have
 \begin{equation*}
    qH^*\Fln\ =\  \ZZ[q_1,\dotsc,q_{n-1}] \otimes_\ZZ H^*\Fln\,,
 \end{equation*}
where $q_1,\dotsc,q_{n-1}$ are indeterminates.
To simplify notation, write \defcolor{$\ZZ[q]$} for $\ZZ[q_1,\dotsc,q_{n-1}]$, and use multiindex notation, for
$\alpha\in\NN^{n-1}$, we have $\defcolor{q^\alpha}\vcentcolon= q_1^{\alpha_1}\dotsb q_{n-1}^{\alpha_{n-1}}$.
The Schubert classes $\{\frakS_w\}_{w \in S_n}$ in $H^*\Fln$
form a $\ZZ[q]$-basis of $qH^*\Fln$.
We write $\defcolor{\frakS_u \ast \frakS_v}$ for the product of Schubert
classes $\frakS_u$ and $\frakS_v$ in $qH^*\Fln$ to distinguish it from
the product $\defcolor{\frakS_u \cdot \frakS_v}$  in $H^*\Fln$.

The ring structure of $qH^*\Fln$ is determined by the \demph{quantum Monk formula}~\cite{FGP}.
For every $u\in S_n$ and $1\leq k<n$ we have,
\begin{equation}
    \label{quantum-Monk-formula}
    \frakS_u * \frakS_{(k,k+1)}
    \ =\ 
    \sum_{\substack{1 \leq i \leq k < j \leq n \\ \ell(u) + 1 = \ell(u(i,j))}} \frakS_{u (i,j)}\ 
    \ +\ 
    \sum_{\substack{1 \leq i \leq k < j \leq n \\ \ell(u) + 1 = \ell(u(i,j)) + 2(j - i)}} q_{i,j} \frakS_{u (i,j)}
\end{equation}
where $\defcolor{q_{i,j}}\vcentcolon=q_iq_{i+1}\dotsb q_{j-1}$ for $i < j$.
The summand on the left above comes from the classical Monk
formula~\eqref{classical-Monk-formula}.
For example, in $qH^*\Fl_4$ we have,
 \begin{equation}
   \label{Eq:q-Monk}
    \frakS_{1432} * \frakS_{(2,3)}
    \ =\  \frakS_{3412} \ +\ \frakS_{2431}
    \ +\ q_{2}\frakS_{1342} \ +\  q_{2,3} \frakS_{1234}\,.
 \end{equation}

Givental and Kim~\cite{GiventalKim}
proved that $qH^*\Fln \simeq \ZZ[q][x] / I^q_n$ for a certain ideal $I^q_n$,
establishing an analog to Borel's presentation \eqref{borel-presentation}.
Fomin, Gelfand, and Postnikov~\cite{FGP} defined \demph{quantum Schubert polynomials}
$\defcolor{\frakS^q_w(x)}$ in $\ZZ[q][x]$ so that the quantum product of
Schubert classes satisfies $\frakS_u \ast \frakS_v \mapsto (\frakS^q_u(x) + I^q_n) (\frakS^q_v(x) + I^q_n)$, for $u,v\in S_n$,
under this isomorphism.
To simplify notation, for $f(x) \in \ZZ[q][x]$, we will write
$\frakS_u \ast f(x) $ for $\frakS^q_u(x) f(x) + I^q_n$.

As explained in Section~\ref{cohomology-ring-of-the-flag-manifold},
every Schur polynomial is a Schubert polynomial:
more precisely, $s_\lambda(x_1, \ldots, x_k) = \frakS_{v(\lambda, k)}(x)$,
where $\lambda \subseteq R_{k, n-k}$.
Analogously, we define the
\demph{quantum Schur polynomial $s_\lambda^q(x_1,\dotsc,x_k)$},
for $\lambda \subseteq R_{k, n-k}$, to be 
${s_\lambda^q(x_1,\dotsc,x_k)} \vcentcolon= \frakS^q_{v(\lambda,k)}(x)$.
Following Equation~\eqref{Eq:PowerSumHook}, the \demph{quantum power sum} is defined to be
\begin{multline}
    \label{Eq:Q-powerSum}
    \qquad  \defcolor{p^q_r(x_1,\dotsc,x_k)}\ \vcentcolon= \\
        s^q_{(r)}(x_1,\dotsc,x_k) - s^q_{(r-1,1)}(x_1,\dotsc,x_k) +
                   \dotsb + (-1)^{r+1} s^q_{(1^r)}(x_1,\dotsc,x_k)\,.
    \qquad
\end{multline}

\subsection{Quantum $k$-Bruhat order}
\label{quantum-k-bruhat-order}

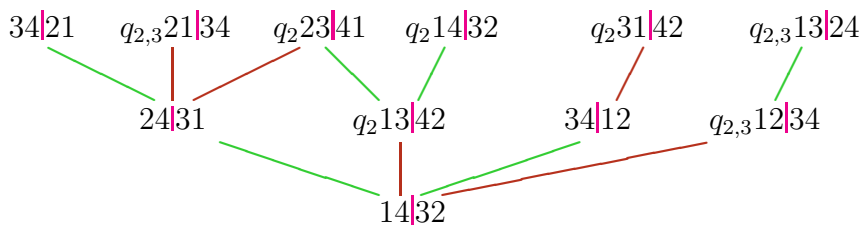
\begin{figure}[htpb]
    \centering
    \begin{picture}(330,80)(10,0)
        \thicklines
        \put( 10, 70){$34\pb21$}          \put( 52, 70){$q_{2,3} 21\pb34$}
        \put(110, 70){$q_{2} 23\pb41$}       \put(160, 70){$q_{2} 14\pb32$}
        \put(230, 70){$q_{2} 31\pb42$}       \put(290, 70){$q_{2,3} 13\pb24$}

        \put( 65,46){\color{LimeGreen}\line(-2,1){40}}
        \put( 72,46){\color{BrickRed}\line(0,1){20}}
        \put( 80,46){\color{BrickRed}\line(2,1){40}}

        \put(150,46){\color{LimeGreen}\line(-1,1){20}}
        \put(165,46){\color{LimeGreen}\line( 1,2){10}}

        \put(240,46){\color{BrickRed}\line( 1,2){10}}

        \put(300,46){\color{LimeGreen}\line( 1,2){10}}

        \put(59.4, 35){$24\p31$}          \put(140, 35){$q_{2} 13\pb42$}
        \put(220, 35){$34\pb12$}         \put(275, 35){$q_{2,3} 12\pb34$}

        \put(150,10){\color{LimeGreen}\line(-3,1){60}}
        \put(158,10){\color{BrickRed}\line( 0,1){20}}
        \put(166,10){\color{LimeGreen}\line(3,1){60}}
        \put(174,10){\color{BrickRed}\line(5,1){100}}
        \put(150,0){$14\pb32$}
    \end{picture}
    \caption{Two levels of the quantum $2$-Bruhat order on $S_4[q]$ above $14\pb32$.}
    \label{fig:quantum-2-Bruhat-poset}
\end{figure}

The terms appearing in the quantum Monk formula~\eqref{quantum-Monk-formula} define a ranked partial order on
the (infinite) set
 \begin{equation*}
    \defcolor{S_n[q]}\  \vcentcolon=\  \big\{ q^\alpha u \mid u\in S_n \text{ and } \alpha\in\NN^{n-1}\big\}\,.
 \end{equation*}
The \demph{quantum $k$-Bruhat order $\leq_k^q$} on $S_n[q]$ is defined by the following cover relations:
\begin{enumerate}
    \item $\defcolor{u\lessdot_k^q u (i,j)}$ if $i\leq k<j$ and $\ell(u)+1=\ell(u (i,j))$,
         so that $u\lessdot_k u (i,j)$;
    \item $\defcolor{u\lessdot_k^q q_{i,j} u (i,j)}$ if $i\leq k<j$ and \rule{0pt}{13pt}
           $\ell(u)+1=\ell(u (i,j))+2(j{-}i)$\,; and
    \item extend $q$-multiplicatively:  $u\leq_k^q q^\beta w$  if and
        only if $q^{\alpha} u\leq_k^q q^{\alpha+\beta} w$ for $u,w\in S_n[q]$ and any $\alpha,\beta\in \NN^{n-1}$. \rule{0pt}{13pt} \vspace{2pt}
\end{enumerate}
We will say that covers of type (1) and their extensions by (3) are \demph{classical}, and covers of type (2) and their
extensions by (3) are \demph{quantum}. 
The \demph{rank function $\ell$} on $S_n[q]$ is given by
$\defcolor{\ell(q^\alpha u)} \vcentcolon= 2 \deg q^\alpha + \ell(u)$, where $\deg q^\alpha$ is the usual degree of the 
monomial $q^\alpha$ in $\ZZ[q]$.
Figure~\ref{fig:quantum-2-Bruhat-poset} shows two levels of the quantum $2$-Bruhat order on $S_4[q]$ starting at the
permutation $14\p32$.
We draw the classical covers in green
and the quantum covers in red.
Compare the first level to the example in~\eqref{Eq:q-Monk}.

As with the $k$-Bruhat order (see Proposition~\ref{nothing-in-between-is-in-between}), there is an intrinsic
formulation of the cover relations for the quantum $k$-Bruhat order.

\begin{proposition}
    \label{prop:quantum-k-bruhat-covers}
    Let $u \in S_n$ and $q^\alpha w \in S_n[q]$.
    Then $u \lessdot_k^q q^\alpha w$ if and only if
    \begin{enumerate}
        \item\label{prop:quantum-k-bruhat-classic-covers}
            $\alpha=0$ and $q^\alpha w =w= u(i, j)$ for some $i \leq k < j$ such that $u(i) < u(j)$, and
            \\
            there does not exist $l$ such that $i < l < j$ and $u(i) < u(l) < u(j)$; or

        \item\label{prop:quantum-k-bruhat-quantum-covers}
            $q^\alpha w = q_{i,j} u(i, j)$ for some $i \leq k < j$ such that $u(i) > u(j)$, and
            \\
            for every $l$ such that $i < l < j$, we have $u(j) < u(l) < u(i)$.
    \end{enumerate}
\end{proposition}

\subsection{Minimal intervals}
\label{sec:Qminim}

As far as we know, there is no analog of the Grassmannian--Bruhat order
related to the quantum $k$-Bruhat order, and so rather than speak of minimal
permutations (as we did in the classical case; see
Section~\ref{minimal-permutations}), we work instead with certain
intervals in the quantum $k$-Bruhat order.
This lack of an analog of the Grassmannian--Bruhat order necessitates the significant work we do in Section~\ref{left-operators}.

Let $u, w \in S_n$ and suppose that $u \leq_k^q q^{\alpha}w$ in the quantum
$k$-Bruhat order on $S_n[q]$.
We say that the interval $[u,q^\alpha w]_k^q$ is \demph{minimal} if
its rank $\ell(q^{\alpha}w)-\ell(u)$ equals $\#\supp(wu^{-1}) - s(wu^{-1})$.
When $q^\alpha=1$, the interval $[u, w]_k^q=[u,w]_k$ is minimal if and only if $wu^{-1}$ is a minimal permutation.
In the proof of Theorem~\ref{thm:q_hook_product} we will show that if $[u,q^\alpha w]_k^q$ is minimal, then $wu^{-1}$ is a
minimal permutation in the sense of Section~\ref{minimal-permutations}.

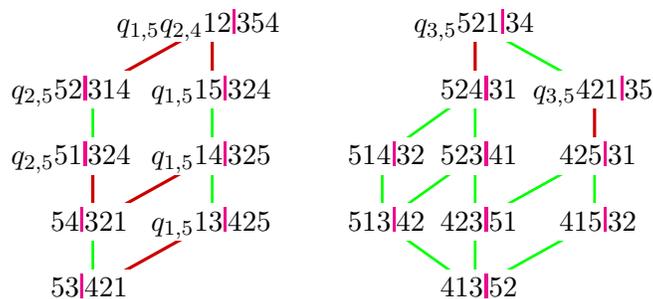
\begin{figure}[htb]
\centering
   \begin{picture}(100,110)(-15,0)
    \put(13,0){\includegraphics{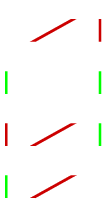}}
                                         \put(25,100){\small$q_{1,5}q_{2,4}12\p354$}
    \put(-15, 75){\small$q_{2,5}52\p314$} \put(38, 75){\small$q_{1,5}15\p324$}  
    \put(-15, 50){\small$q_{2,5}51\p324$} \put(38, 50){\small$q_{1,5}14\p325$}  
    \put(  0, 25){\small$54\p321$}       \put(38, 25){\small$q_{1,5}13\p425$}  
    \put(  0,  0){\small$53\p421$}  
   \end{picture}
   \qquad
   \begin{picture}(120,110)(-35,0)
    \put(-25,0){\includegraphics{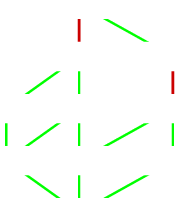}}
                                        \put(-10,100){\small$q_{3,5}521\p34$}  
                                        \put(0, 75){\small$524\p31$}   \put(35,75){\small$q_{3,5}421\p35$}  
    \put(-35,50){\small$514\p32$}       \put(0, 50){\small$523\p41$}   \put(45,50){\small$425\p31$}  
    \put(-35,25){\small$513\p42$}       \put(0, 25){\small$423\p51$}   \put(45,25){\small$415\p32$}  
                                        \put(0,  0){\small$413\p52$}  
  \end{picture}

   \caption{Two minimal quantum intervals.  Note that $q_{1,5}q_{2,4}=q_{1,4}q_{2,5}=q_1q_2^2q_3^2q_4$.}
   \label{F:more_minimal}
\end{figure}

The interval 
on the left in Figure~\ref{F:k-Bruhat} is a minimal interval in the quantum
$5$-Bruhat order in $S_8$.
Figure~\ref{F:more_minimal} displays two minimal intervals in quantum $k$-Bruhat order on $S_5[q]$ for $k=2,3$.

\subsection{The Murnaghan--Nakayama rule for $qH^*\Fln$}

As for Corollary~\ref{C:classical-schubert-times-powersum},
we derive the quantum Murnaghan--Nakayama formula from an intermediate result
describing the quantum product of a Schubert class $\frakS_u$ by a quantum hook Schur
polynomial $s_{(b, 1^{a-1})}^q(x_1,\dotsc,x_k)$.

\begin{theorem}
    \label{thm:q_hook_product}
    Let $u \in S_n$\,, $a \leq k$ and $b \leq n-k$. Then
    \begin{equation*}
        \frakS_u
        \ast
        s^q_{(b,1^{a-1})} (x_1, \dots, x_k)
        = \sum
        \left(\begin{smallmatrix} s(wu^{-1})-1 \\ \het(wu^{-1})-a \end{smallmatrix}\right)
        q^\alpha
        \frakS_{w}\,,
    \end{equation*}
    the sum over all minimal intervals $[u,  q^{\alpha}w]_k^q$ such that $\ell(q^\alpha w) - \ell(u)= a{+}b{-}1$.
\end{theorem}

The proof depends on  Theorem~\ref{Th:nonZero} and Lemma~\ref{L:Quantum_is_Classical}, which will be proved in
subsequent sections.

\begin{proof}
    For  $u,w \in S_n$, a partition $\lambda$, $k \in \NN$ with $1\leq k$ and $\alpha\in \NN^{n-1}$, define
    quantum Littlewood--Richardson coefficients $C^{q^\alpha w}_{u,v(\lambda,k)}\in q^\alpha \NN$ by 
    \begin{equation}
        \label{eq-definition-Cq}
        \frakS_u \ast s^q_\lambda(x_1, \ldots, x_k)
        \ =\ \sum_{q^\alpha w \in S_n[q]} C^{q^\alpha w}_{u,v(\lambda,k)} \frakS_w\,.
    \end{equation}
    We must prove that
    \begin{equation*}
        C^{q^\alpha w}_{u,v((b, 1^{a-1}),k)}\ =\ 
        \begin{cases}
            q^\alpha \left(\begin{smallmatrix} s(wu^{-1})-1 \\ \het(wu^{-1})-a \end{smallmatrix}\right),
            &
            \begin{array}[t]{lr}
                \text{if $[u, q^{\alpha}w]_k^q$ is minimal and} \\[-0.5ex]
                \qquad\qquad\qquad \text{$\ell(q^\alpha w) - \ell(u) = a{+}b{-}1$,}
            \end{array}
            \\
            0, &
            \begin{array}[t]{lr}
                \text{otherwise.}
            \end{array}
        \end{cases}
    \end{equation*}

    Suppose that $[u, q^{\alpha}w]_k^q$ is a minimal interval
    with $\ell(q^\alpha w) - \ell(u) = a{+}b{-}1$.
    By Theorem~\ref{Th:nonZero},
    there exists a hook partition $\gamma$
    with $|\gamma| = \ell(q^\alpha w) - \ell(u)$
    and $C^{q^\alpha w}_{u,v(\gamma,k)} \neq 0$.
    By Lemma~\ref{L:Quantum_is_Classical},
    there exist permutations $y, z \in S_n$ satisfying
    \begin{equation}
        \label{properties-of-y-and-z}
        y\ \leq_k\ z,
        \quad
        \ell(z) - \ell(y)\ =\ \ell(q^\alpha w) - \ell(u) = a{+}b{-}1,
        \quad
        zy^{-1}\ =\ wu^{-1},
    \end{equation}
    and such that
    \begin{equation}
        \label{eq-leung-li-reduction}
        C^{q^\alpha w}_{u,v(\mu,k)}\ =\  q^\alpha c^z_{y,v(\mu,k)}
        \quad \text{for all partitions $\mu$ with $|\mu| = a{+}b{-}1$\,.}
    \end{equation}
    This implies that $z y^{-1}$ is a minimal permutation, since
    \begin{align*}
        \ell(z) - \ell(y)
        = \ell(q^\alpha w) - \ell(u)
        = \#\supp(w u^{-1}) - s(w u^{-1})
        = \#\supp(z y^{-1}) - s(z y^{-1}),
    \end{align*}
    where the middle equality follows from the assumption that
    $[u, q^{\alpha}w]_k^q$ is minimal.
    Then
    \begin{equation}
        C^{q^\alpha w}_{u,v((b, 1^{a-1}),k)}
        \ =\ q^\alpha c^z_{y,v((b, 1^{a-1}),k)}
        \ =\ q^\alpha \left(\begin{smallmatrix} s(zy^{-1})-1 \\ \het(zy^{-1})-a \end{smallmatrix}\right)
        \ =\ q^\alpha \left(\begin{smallmatrix} s(wu^{-1})-1 \\ \het(wu^{-1})-a \end{smallmatrix}\right),
    \end{equation}
    where the first equality follows from \eqref{eq-leung-li-reduction}
    and the second equality follows from Proposition~\ref{P:classical_hook}
    and that $z y^{-1}$ is a minimal permutation.

    Conversely,
    suppose $C^{q^\alpha w}_{u,v((b, 1^{a-1}), k)} \neq 0$.
    By Lemma~\ref{L:Quantum_is_Classical},
    there exist $y, z \in S_n$ satisfying
    \eqref{properties-of-y-and-z} and \eqref{eq-leung-li-reduction}.
    Setting $\mu = (b, 1^{a-1})$ in \eqref{eq-leung-li-reduction} gives
    \begin{equation}
        C^{q^\alpha w}_{u,v((b, 1^{a-1}),k)}
        \ =\ q^\alpha c^z_{y,v((b, 1^{a-1}),k)}.
    \end{equation}
    Since the left hand side is assumed to be nonzero, we have that
    $c^z_{y,v((b, 1^{a-1}),k)}$ is nonzero.
    It now follows from Proposition~\ref{P:classical_hook} that
    $z y^{-1}$ is a minimal permutation so that
    \begin{equation*}
        \ell(q^\alpha w) - \ell(u)
        = \ell(z) - \ell(y)
        = \#\supp(z y^{-1}) - s(z y^{-1})
        = \#\supp(w u^{-1}) - s(w u^{-1}),
    \end{equation*}
    where the first and last equalities follow from \eqref{properties-of-y-and-z}.
    Thus, $[u, q^\alpha w]_k^q$ is a minimal interval such that
    $\ell(q^\alpha w) - \ell(u) = a{+}b{-}1$.
\end{proof}

The following was conjectured in~\cite{Morrison}.
It follows from Theorem~\ref{thm:q_hook_product} in the same way that
Corollary~\ref{C:classical-schubert-times-powersum} follows from Proposition~\ref{P:classical_hook},
and is the signature result in this paper.

\begin{corollary}[Quantum Murnaghan--Nakayama Formula]
    \label{C:quantumMN}
    Let $u \in S_n$. Then
    \begin{equation*}
        \frakS_u \ast p_r^q(x_1,\dotsc,x_k) = \sum (-1)^{\het(wu^{-1})+1} q^\alpha \frakS_{w}
    \end{equation*}
    the sum over all minimal intervals $[u, q^{\alpha}w]_k^q$ of rank $r$ such that $wu^{-1}$ is a
    single cycle.
\end{corollary}

\begin{example}\label{Ex:MN_example}
In $qH^*\Fl_8$, if $u=68235741$, then the product $\frakS_{u}\ast p^q_4(x_1,\dotsc,x_5)$ is
\begin{gather*}
      \frakS_{(2,3,5,7,4)u}
    + \frakS_{(2,4,3,5,7)u}
    + \frakS_{(3,5,6,7,4)u}
    - \frakS_{(2,3,5,6,7)u}
    \\
    + q_{5,7}\frakS_{(2,3,4,7,5)u}
    + q_{5,7}\frakS_{(3,4,6,7,5)u}
    + \defcolor{q_{5,8}\frakS_{(1,6,7,3,5)u}}
    + q_{5,8}\frakS_{(1,7,2,3,5)u}
    \\
    - q_{5,8}\frakS_{(1,6,7,5,4)u}
    - q_{5,8}\frakS_{(1,7,5,3,4)u}
    - q_{5,8}\frakS_{(1,7,4,3,5)u}
    - q_{2,6}\frakS_{(2,3,5,7,8)u}
    \\
    + q_{2,8}\frakS_{(1,7,3,5,8)u}
    + q_{2,8}\frakS_{(1,8,2,3,4)u}
    + q_{2,8}\frakS_{(1,6,7,5,8)u}
    + q_{2,8}\frakS_{(1,7,5,6,8)u}
    -q_{2,8}q_{5,7}\frakS_{(1,7,5,4,8)u} \,.
\end{gather*}
Note that unlike in ordinary cohomology, this is not an identity among quantum Schubert polynomials $\frakS^q_w(x)$,
but rather among their images in the quotient ring $\ZZ[q][x]/I^q_n$.
\end{example}

\subsection{Remarks regarding~\cite[Thm.~15]{MPP}}\label{SS:R_MPP}

Proposition~\ref{P:classical_hook} may be deduced from a stronger result of M\'esz\'aros, Panova, and
Postnikov~\cite[Thm.~8]{MPP}.
There, they establish a formula for a `lift' of the Schur polynomial
$s_{(b,1^{a-1})}(x_1,\dotsc,x_k)$ to the Fomin--Kirillov~\cite{FominKirillov}
algebra ${\mathcal E}_n$, which is a subalgebra (generated by certain Dunkl elements $\theta_i$) of the quotient of a free
associative algebra with generators $x_{i,j}$ for $1\leq i<j\leq n$ by certain relations, which include $x_{i,j}^2=0$.
Their paper contains another result~\cite[Thm.~15]{MPP} which is a quantum analog of~\cite[Thm.~8]{MPP}, giving a formula
for a `quantum lift' of $s^q_{(b,1^{a-1})}(x_1,\dotsc,x_k)$ to the quantum Fomin--Kirillov algebra
$\mathcal{E}^q_n$~\cite{Postnikov}.
This formula would imply Theorem~\ref{thm:q_hook_product}.
Unfortunately, as we now sketch, there is a gap in the proof of~\cite[Thm.~15]{MPP} (we have communicated this to the
authors).

The result~\cite[Thm.~15]{MPP} is deduced from the proof of~\cite[Thm.~8]{MPP} using a technical
result~\cite[Lem.~13]{MPP} which roughly asserts the following:
If an identity involving positive terms holds in ${\mathcal E}_n$, and if no relations  $x_{i,j}^2=0$ were used to deduce
the identity, then the same identity holds in ${\mathcal E}^q_n$.
We are convinced this lemma is correct, but it unfortunately does not apply in the case
invoked to prove~\cite[Thm.~15]{MPP}.
The proof of~\cite[Thm.~8]{MPP} uses in a fundamental way~\cite[Cor.~10]{MPP}, that
$e_a(\theta_1,\dotsc,\theta_a)h_b(\theta_1,\dotsc,\theta_a)=0$ in $\mathcal{E}_{a+b}$ (in their notation).
However, this identity can only be proven directly in  $\mathcal{E}_{a+b}$ using the relations $x_{i,j}^2=0$.
For example, when $a=b=1$, as $\theta_1=x_{1,2}$, it becomes $0=e_1(x_{1,2})h_1(x_{1,2})=x_{1,2}^2$.

%
%
%
\section{Quantum equals classical}
\label{quantum-equals-classical}

A recurring theme in quantum Schubert calculus is that oftentimes the numerical part
$\defcolor{N^{w,\alpha}_{u,v}}\vcentcolon=q^{-\alpha}C^{q^\alpha w}_{u,v}(\in\NN)$ of a quantum
Littlewood--Richardson number (a Gromov--Witten invariant)
is naturally identified with a related classical Littlewood--Richardson number $c^{w'}_{u',v'}$.
This is seen in the quantum Monk formula~\eqref{quantum-Monk-formula} in which these numerical parts are all 1 as in the
classical Monk formula~\eqref{classical-Monk-formula}.
Bertram's quantum Giambelli formula~\cite{Bertram} and the results of Buch, Kresch, and Tamvakis~\cite{BKT} are other
examples. 
Leung and Li~\cite{LeungLi2012} discuss this thoroughly and then establish very general results, including some that
naturally identify different Gromov--Witten invariants on the same flag variety. 
Their result is type-free and expressed in the language of algebraic groups, and they use it to deduce a 
weaker version of the result we need, Lemma~\ref{L:Quantum_is_Classical}.
We first state some results from~\cite{LeungLi2012} for the type A flag variety using the notation from Section~\ref{S:qCoh} and extract
a useful statement from one of their proofs, from which we deduce  Lemma~\ref{L:Quantum_is_Classical}.

Fix $i\in[n{-}1]$.
Let $\defcolor{e_i}\vcentcolon=(0,\dotsc,0,1,0,\dotsc,0)$ be the $i$th standard basis vector, so that $q^{e_i}=q_i$.
Write $\defcolor{s_i}\vcentcolon=(i,i{+}1)$, the $i$th simple transposition.
We will also need the linear map $\varpi_i\colon\ZZ^{n-1}\to\ZZ$, defined by
$\defcolor{\varpi_i}(\alpha)\vcentcolon=-\alpha_{i-1}+2\alpha_i-\alpha_{i+1}$, where $\alpha_0=\alpha_n=0$.   
Thus
 \begin{multline*}
   \qquad  q_{2,4}=q^{(01100)}=q^{e_2+e_3}\,,   \qquad
   6{\color{red}\underline{35}}214 s_2\ =\ 6{\color{red}\underline{53}}214\,,\quad
   63{\color{red}\underline{52}}14 s_3\ =\ 63{\color{red}\underline{25}}14\,, \\
   \varpi_2(\underline{122}32)\ =\ -1+4-2=1\,,\quad \mbox{ and }\quad  \varpi_4(12\underline{232})\ =\ -2+6-2=2\,.\qquad 
 \end{multline*}
 For $i\in[n{-}1]$, we also define a map $\defcolor{\sg_i}\colon S_n\to\{0,1\}$ which records whether or not $u\lessdot_i
 us_i$, 
 \[
    \sg_i(u)\ \vcentcolon=\ \left\{ \begin{array}{rcl}
       1& &\mbox{if\ } u s_i\lessdot_i  u     \quad\mbox{equivalently, if\ } u(i+1)<u(i)
      \\
       0& &\mbox{if\ } u    \lessdot_i  u s_i \quad\mbox{equivalently, if\ } u(i)<u(i+1)
    \end{array}\right. .
 \]                 
 That is, $\sg_i(u)=1$ when $u$ has a descent at $i$ and it equals $0$ when $u$ has an ascent at $i$. 
  
\begin{proposition}~\cite[Theorem 1.1]{LeungLi2012}\label{P:LL_Th1.1}
  Let $x,y,z\in S_n$ and $\alpha\in\NN^{n-1}$. 
  \begin{enumerate}
    \item We have that $N^{z, \alpha}_{x,y}=0$ unless
           $\sg_i(z)+\varpi_i(\alpha) \leq \sg_i(x)+\sg_i(y)$ for all $i\in[n{-}1]$.\vspace{2pt}
           
    \item Suppose that   $\sg_i(z)+\varpi_i(\alpha) = \sg_i(x)+\sg_i(y)=2$ for some $i\in[n{-}1]$, then
  \[
           N^{z,\alpha}_{x,y}\ =\ N^{z,\alpha-e_i}_{x s_i, y s_i}\ =\ \left\{ \begin{array}{lcl}
                 N^{z s_i, \alpha-e_i}_{x, ys_i}& &\mbox{if\/\ } \sg_i(z)=0\vspace{5pt}\\
                 N^{z s_i, \alpha}_{x, ys_i}&  &\mbox{if\/\ } \sg_i(z)=1\end{array}\right. .
 \]
 \end{enumerate}
\end{proposition}

Note that if $N^{w,\alpha}_{u,v}\neq 0$, then by part (1), we have $\varpi_i(\alpha)\leq 2$ for each $i\in[n{-}1]$.
This also constrains $u,v,w$.
For example, if $N^{w,\alpha}_{u,v(\lambda,k)}\neq 0$ and $i\neq k$,  then $\sg_i(v(\lambda,k))=0$ and
Proposition~\ref{P:LL_Th1.1}(1) becomes 
\begin{equation}\label{Eq:Grass_bound}
  \sg_i(w)\ +\ \varpi_i(\alpha)\ \leq\ \sg_i(u)\ \leq\ 1\,.
\end{equation}

We first recall two lemmas from~\cite{LeungLi2012}.

\begin{proposition}\label{P:LL_lemmas}
   Suppose that $\alpha\in\NN^{n-1}$ is a nonnegative nonzero  exponent vector.
   \begin{enumerate}
   \item  \emph{(\cite[Lemma 2.8]{LeungLi2012})}
        There is an $i\in[n{-}1]$ with $\varpi_i(\alpha)>0$.
    \item  \emph{(\cite[Lemma 2.10]{LeungLi2012})}
        If $\varpi_i(\alpha)>0$ for a unique $i\in[n{-}1]$, then $\varpi_i(\alpha)\geq 2$.
   \end{enumerate}
\end{proposition}

If $N^{w,\alpha}_{u,v(\lambda,k)}\neq 0$ and $i\neq k$ in Proposition~\ref{P:LL_lemmas}(1), then $\varpi_i(\alpha)=1$ as
$\varpi_i(\alpha)>1$ 
contradicts~\eqref{Eq:Grass_bound}, and this implies that  $\sg_i(u)=1$ and $\sg_i(w)=0$.

Theorem 1.2 in~\cite{LeungLi2012} is their ``quantum equals classical'' result.
It states that if $N^{w, \alpha}_{u,v(\lambda,k)}\neq 0$, then there exist $w',u'\in S_n$ such
that  $N^{w, \alpha}_{u,v(\lambda,k)}=N^{w', {\bzero}}_{u',v(\lambda,k)}= c^{w'}_{u',v(\lambda,k)}$.
A consequence of their proof is that $w u^{-1}=w'(u')^{-1}$, but we require more than that.
We state and prove the result we will need; this is extracted from the proof of Theorem 1.2 in~\cite{LeungLi2012}.

\begin{proposition}\label{P:LLThm1.2}
  Suppose that $u,w\in S_n$, $\lambda$ is a partition, $\alpha\in\NN^{n-1}$ is nonzero, and $N^{w,\alpha}_{u,v(\lambda,k)}\neq 0$.
  \begin{enumerate}
  \item There exists $i\in[n{-}1]$ such that $\sg_i(u)=1$, $\sg_i(w)=0$, and
  \[
    \mbox{either }\ \varpi_i(\alpha) = 1\ \mbox{if }i\neq k
    \quad\mbox{or}\quad
    \varpi_i(\alpha) = 2\ \mbox{if }i= k    \,.
  \]

   \item For any partition $\mu$ with $|\mu|=|\lambda|$, with $i$ as in (1), we have
        $ N^{w,\alpha}_{u,v(\mu,k)}=N^{ws_i, \alpha-e_i}_{us_i,v(\mu,k)}$.
   \end{enumerate}
 \end{proposition}
 \begin{proof}
     We first prove statements (1) and (2) in the case when $\mu=\lambda$.
   By Proposition~\ref{P:LL_lemmas}(1), there is an $i\in[n{-}1]$ such that $\varpi_i(\alpha)>0$.
   Suppose that $i$ is not unique with this property.
   Then we may assume that $i\neq k$, from which it follows that $\sg_i(v(\lambda, k)) = 0$.
   As $N^{w,\alpha}_{u,v(\lambda,k)}\neq 0$, Equation~\eqref{Eq:Grass_bound} implies that $\sg_i(u)=1=\varpi_i(\alpha)$ and $\sg_i(w)=0$.
   Thus $\sg_i(ws_i)=1=\sg_i(v(\lambda,k)s_i)$, so that
   \[
     \sg_i(ws_i)\ +\ \varpi_i(\alpha)\ =\ \sg_i(u)\ +\ \sg_i(v(\lambda,k) s_i)\ =\ 2\,.
   \]
   By Proposition~\ref{P:LL_Th1.1}(2) with $x=u$, $y=v(\lambda,k) s_i$, and $z=w s_i$, we have that 
   \[
     N^{ ws_i,\alpha}_{u, v(\lambda,k) s_i}\ =\ N^{ws_i, \alpha-e_i}_{u s_i, v(\lambda,k)}\ =\ N^{w,\alpha}_{u,v(\lambda,k)}\,.
   \]
   (As $s_i^2=1$ and $\sg_i(z)=\sg_i(ws_i)=1$).
   
   Now suppose that there is a unique $i$ with $\varpi_i(\alpha)>0$.
   Then by Proposition~\ref{P:LL_lemmas}(2), we have $\varpi_i(\alpha)\geq 2$.
   As $N^{w,\alpha}_{u,v(\lambda,k)}\neq  0$,  Proposition~\ref{P:LL_Th1.1}(1) implies that 
   \[
     \sg_i(w)\ =\ 0\,,\ \ \ 
     \varpi_i(\alpha)\ =\ 2\,,\ \ \mbox{and}\ \
     \sg_i(u)\ =\ \sg_i(v(\lambda,k))\ = \ 1\,.
   \]
   Since $\sg_i(v(\lambda,k))=1$, we must have that $i=k$, as $v(\lambda,k)$ has a unique descent at $k$.
   As $\sg_i(w)=0$,  Proposition~\ref{P:LL_Th1.1}(2) with $x=v(\lambda,k)$, $y=u$, and $z=w$ implies that
   $N^{w,\alpha}_{v(\lambda,k), u}=N^{w s_i, \alpha-e_i}_{v(\lambda,k), u s_i}$.
   As commutativity implies that $N^{z,\alpha}_{x,y}=N^{z,\alpha}_{y,x}$, we have
   $N^{ws_i, \alpha-e_i}_{u s_i,v(\lambda,k)}=N^{w,\alpha}_{u,v(\lambda,k)}$.

   This proves statement (1) and statement (2), in the case when $\mu=\lambda$.
   The assumption that $|\mu|=|\lambda|$ implies that for all $i\in[n{-}1]$, $\sg_i(v(\lambda,k))=\sg_i(v(\mu,k))$, as
   $v(\mu,k)$ is the identity permutation when $|\mu|=0$, and otherwise it has a unique descent at $k$.
   The arguments we gave for $N^{w,\alpha}_{u,v(\lambda,k)}$ will then also hold for $N^{w,\alpha}_{u,v(\mu,k)}$, which completes the proof.
 \end{proof}

\begin{lemma} \label{L:Quantum_is_Classical}
  Suppose that $u,w\in S_n$, $\lambda$ is a partition, $\alpha\in\NN^{n-1}$, and $C^{q^\alpha w}_{u,v(\lambda,k)}\neq 0$.
  Then there exist $y,z\in S_n$ with $y\leq_k z$, $\ell(z)-\ell(y)=|\lambda|$, and $zy^{-1}=wu^{-1}$
  such that for every partition $\mu$ with $|\mu|=|\lambda|$, we have
  $q^{-\alpha}C^{q^\alpha w}_{u,v(\mu,k)}=N^{w,\alpha}_{u,v(\mu,k)}=c^{z}_{y,v(\mu,k)}$. 
\end{lemma}
\begin{proof}
  We use induction on $|\alpha|$.
  Let $\alpha\in\NN^{n-1}$.
  When $|\alpha|=0$, $\alpha={\bzero}$ and $C^{q^{\bf 0}w}_{u,v(\mu,k)}=c^w_{u,v(\mu,k)}$.
  Suppose that $\alpha\neq\bzero$.
  Let $i\in[n{-}1]$ be the index guaranteed by Proposition~\ref{P:LLThm1.2}(1).
  By Proposition~\ref{P:LLThm1.2}(2), if $\mu$ is a partition with $|\mu|=|\lambda|$, then 
  $N^{w,\alpha}_{u,v(\mu,k)}=N^{ws_i, \alpha-e_i}_{us_i,v(\mu,k)}$.
  In particular,   $wu^{-1}=ws_i(u s_i)^{-1}$  and $N^{ws_i, \alpha-e_i}_{us_i,v(\lambda,k)}\neq 0$, which implies that
  $us_i\leq_k^q q^{\alpha-e_i}w s_i$.
  This completes the proof.
\end{proof}

\begin{example}\label{Ex:q-minimal}
  We illustrate these results.
  In Example~\ref{Ex:MN_example}, the term $q_{5,8}\frakS_{(1,6,7,3,5)u}$ appears in the product
  $\frakS_{u}\ast p^q_4(x_1,\dotsc,x_5)$.
  As $78251\p346=(1,6,7,3,5)u$, by Equation~\eqref{Eq:Q-powerSum}, there is a hook partition $\lambda\subseteq R_{5,3}$ with
  $|\lambda|=4$ such that $N^{78251\pp346, q_{5,8}}_{68235\pp741, v(\lambda,5)}\neq 0$.
  Since $q_{5,8}=q^\alpha$ for $\alpha=(0,0,0,0,1,1,1)$, we have $\varpi_7(\alpha)=1$, so that by Proposition~\ref{P:LLThm1.2}(2), 
  $N^{78251\pp3\rc{46}, q_{5,\rc{8}}}_{68235\pp7\rc{41}, v(\lambda,5)}= N^{78251\pp3\rc{64}, q_{5,\rc{7}}}_{68235\pp7\rc{14}, v(\lambda,5)}$,
  as $q^{\alpha-e_7}=q_{5,7}$.
  We may apply Proposition~\ref{P:LLThm1.2}(2) again as $\varpi_6(\alpha-e_7)=1$ to obtain that 
  $N^{78251\pp\rc{36}4, q_{5,\rc{7}}}_{68235\pp\rc{71}4, v(\lambda,5)}= N^{78251\pp\rc{63}4, q_{5,6}}_{68235\pp\rc{17}4, v(\lambda,5)}$.

  If $\beta=\alpha-e_7-e_6$, so that $q^\beta=q_{5,6}=q_5$, then $i=5=k$ is the only index with $\varpi_i(\beta)>0$ and we have 
  $\varpi_5(\beta)=2$.
  By Proposition~\ref{P:LLThm1.2}(2) again, 
  $N^{7825\rc{1\pp6}34, q_5}_{6823\rc{5\pp1}74, v(\lambda,5)}= c^{7825\rc{6\pp1}34}_{6823\rc{1\pp5}74, v(\lambda,5)}$.
  Figure~\ref{F:qminimal} shows the interval $[68231\p574, 78256\p134]_5$, which has one peakless chain (the vertical one in blue, with
  labels $6,3,1,3$). 
\begin{figure}[htb]
 \centering
\begin{picture}(145,122)(-50,-2)
    \put(-31,0){\includegraphics{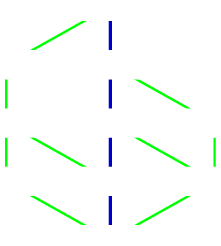}}
                                      \put(0,112){\small$78256\p134$}
    \put(-50, 84){\small$78236\p154$} \put(0, 84){\small$7825{\color{red}\underline{3}}\p164$} 
    \put(-50, 56){\small$78235\p164$} \put(0, 56){\small$7825{\color{red}\underline{1}}\p364$} \put(50,56){\small$68253\p174$}      
    \put(-50, 28){\small$68235\p174$} \put(0, 28){\small$782{\color{red}\underline{3}}1\p564$} \put(50,28){\small$68251\p374$} 
                                      \put(0,  0){\small${\color{red}\underline{6}}8231\p574$} 
\end{picture}
\qquad\qquad
\begin{picture}(175,122)(-60,-2)
    \put(-41.5,0){\includegraphics{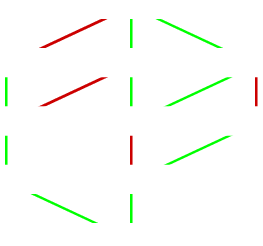}}
                                      \put(-6,112){\small$q^\alpha78251\p346$}
    \put(-60, 84){\small$78256\p341$} \put(-6, 85){\small$q^\alpha78231\p546$} \put(54,84){\small$q^\alpha68251\p347$} 
    \put(-60, 56){\small$78236\p541$} \put(-6, 56){\small$q^\alpha68231\p547$} \put(60,56){\small$68257\p341$}      
    \put(-60, 28){\small$78235\p641$} \put( 0, 28){\small$68237\p541$}
                                      \put( 0,  0){\small$68235\p741$} 
\end{picture}

 \caption{Two minimal intervals in the quantum 5-Bruhat order on $S_8[q]$.}
  \label{F:qminimal}
\end{figure}
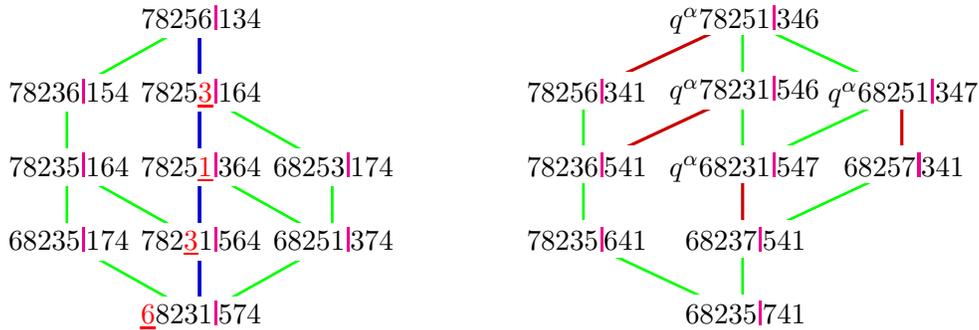


As this has height 3 (as does the minimal permutation $(1,6,7,3,5)$), we see that \newline
$c^{7825\rc{6\pp1}34}_{6823\rc{1\pp5}74, v((2,1,1,0,0),5)}=1$, and thus all of the other coefficients considered above are 1.

Results in~\cite{BS-Bruhat} also show that $c^{78256\pp134}_{68231\pp574, v((2,2,0,0,0),5)}=1$, and thus
 \begin{eqnarray*}
    1\ =\ c^{78256\pp134}_{68231\pp574, v((2,2,0,0,0),5)}
    &=&  N^{78251\pp634, q_{5,6}}_{68235\pp174, v((2,2,0,0,0),5)}\\
    &=&  N^{78251\pp364, q_{5,7}}_{68235\pp714, v((2,2,0,0,0),5)}
    \ =\ N^{78251\pp346, q_{5,8}}_{68235\pp741, v((2,2,0,0,0),5)}\,,
 \end{eqnarray*}
by Lemma~\ref{L:Quantum_is_Classical}.
For $\lambda\neq(2,1,1,0,0), (2,2,0,0,0)$ with $|\lambda|=4$, these coefficients vanish.\hfill$\diamond$ 
\end{example}

Observe that the intervals may change when applying a reduction as in Proposition~\ref{P:LLThm1.2}(2), even though these
reductions do not change $wu^{-1}$.
Indeed, the two intervals in Figure~\ref{F:qminimal} demonstrate this.
(In~\cite{BS-Bruhat}, it was shown that the classical interval $[u,w]_k$ depends only on $wu^{-1}$.)
We show that the coefficient $N^{w,\alpha}_{u,v(\lambda,k)}$ depends only upon $wu^{-1}$.

\begin{theorem}\label{Th:quantumIndependence}
  For a partition $\lambda$, permutations $w,u,\overline{w},\overline{u}$ in $S_n$, exponent vectors
  $\alpha,\overline{\alpha}\in\NN^{n-1}$ and indices $k,l<n$
  if both coefficients $N^{w,\alpha}_{u,v(\lambda,k)}$ and
  $N^{\overline{w},\overline{\alpha}}_{\overline{u},v(\lambda,l)}$ are nonzero and $wu^{-1}=\overline{w}\,\overline{u}^{-1}$,
  then the coefficients are equal.
 \end{theorem}
\begin{proof}
  Suppose that $N^{w,\alpha}_{u,v(\lambda,k)}\neq 0$.
  By Lemma~\ref{L:Quantum_is_Classical}, there exist $y,z\in S_n$ with $y\leq_k z$ such that $\ell(z)-\ell(y)=|\lambda|$,
  $zy^{-1}=wu^{-1}$  and    $N^{w,\alpha}_{u,v(\lambda,k)}=c^{z}_{y,v(\lambda,k)}$.

  Suppose also that $N^{\overline{w},\overline{\alpha}}_{\overline{u},v(\lambda,l)}\neq 0$.
  By Lemma~\ref{L:Quantum_is_Classical} again, there exist  $\overline{y},\overline{z}\in S_n$ with
  $\overline{y}\leq_l \overline{z}$ such that $\ell(\overline{z})-\ell(\overline{y})=|\lambda|$,
  $\overline{z}\,\overline{y}^{-1}=\overline{w}\,\overline{u}^{-1}$ and   
  $N^{\overline{w},\alpha}_{\overline{u},v(\lambda,l)}=c^{\overline{z}}_{\overline{y},v(\lambda,l)}$.

  Since $w u^{-1}= \overline{w}\,\overline{u}^{-1}$, we have that
  $zy^{-1}=\overline{z}\,\overline{y}^{-1}$.
  Then~\cite[Thm.\ E]{BS-Bruhat} implies that $c^{z}_{y,v(\lambda,k)}=c^{\overline{z}}_{\overline{y},v(\lambda,l)}$,
  which completes the proof.
\end{proof}

%
\section{Left operators}
\label{left-operators}

We will prove the following quantum analog of the first part of Proposition~\ref{P:peaklessChains}.

\begin{theorem}
    \label{Th:nonZero}
    Let $u,w\in S_n$ and $q^\alpha w\in S_n[q]$ be such that $[u,q^\alpha w]^q_k$ is a minimal interval.
    Then there is a hook partition $\lambda$ with $|\lambda|=\ell(q^\alpha w)-\ell(u)$ such that
    $C^{q^\alpha w}_{u,v(\lambda,k)}$ is nonzero.
\end{theorem}

Our proof of Theorem~\ref{Th:nonZero} is given at the end of this section.
We sketch the organization of this section and that proof.
We first extend the approach in~\cite{BS-Monoid} (see also~\cite{Benedetti_Bergeron}) of studying saturated chains in the
$k$-Bruhat order using operators.
We introduce operators $\qv$ with an action $\kaction$ on $S_{n}[q]$ and an equivalence relation in Section~\ref{SS:monoid} that encodes
saturated chains in the quantum $k$-Bruhat order as follows:
Let $\qv$ be an operator, $u,w\in S_n$, and $\alpha\in\NN^{n-1}$ be such that $u\leq^q_k \qv\kaction u = q^\alpha w$.
Elements of the equivalence class of $\qv$ are in bijection with chains in $[u,q^\alpha w]^q_k$.
The operator $\qv$ is \demph{classical} if $q^\alpha=1$, in which case $\qv$ is an operator from~\cite{BS-Monoid}.

We show in Section~\ref{SS:diagrammatic-notation} that these operators and the equivalence relation exhibit symmetries that
are enjoyed by the quantum $k$-Bruhat order, the most important being Postnikov's cyclic shift~\cite{Postnikov_symmetry},
which induces a cyclic symmetry of the coefficients $C^{q^\alpha w}_{u,v}$.
Relations of degree two are studied in Section~\ref{S:degreeTwo}.
In Section~\ref{crossing-products} we associate a graph to an operator $\qv$, and show that if
$\qv\kaction u=q^\alpha w$ and if the interval $[u,q^\alpha w]_k^q$ is minimal, then the graph of $\qv$ is a (noncrossing) forest.
In Section~\ref{ss:Paths} we study operators whose graph is a path.
In Section~\ref{SS:trees}, 
we use induction to prove that a noncrossing forest $\qv$ is equivalent to a composition $\R\C$ of what we term a row ($\R$) and
a column ($\C$), and that $\qv$ (and hence $\R\C$) may be cyclically shifted to a classical operator $\qv'$, which is also
a composition $\R'\C'$ of a row $\R'$ and a column $\C'$ in the classical sense of~\cite{Benedetti_Bergeron}.
The product $\R'\C'$ corresponds to a peakless chain in the classical interval corresponding to $\qv'$.
Given this, Postnikov's cyclic symmetry~\cite[Thm.~4]{Postnikov_symmetry} and
Proposition~\ref{pieri-formulas-version-2} enable us to deduce Theorem~\ref{Th:nonZero}.

\subsection{Definition of the left operators}
\label{SS:monoid}

We first recall the results of~\cite{BS-Monoid} concerning saturated chains in the classical $k$-Bruhat order
and then extend part of its approach to the quantum $k$-Bruhat order.

\subsubsection{Classical left operators}
\label{SS:classical-left-operators}

Let \defcolor{$\calF_n$} be the free monoid generated by the $\binom{n}{2}$ symbols
$\{\defcolor{\qv_{a,b}} \mid a<b \mbox{ with }a,b\in[n]\}$.
Let \defcolor{$S_n\cup\{\bzero\}$} be the symmetric group $S_n$ with an
absorbing element \defcolor{$\bzero$} adjoined.
We often omit the comma between indices, thus
$\qv_{ab}=\qv_{a,b}$.

Both the Grassmannian--Bruhat order and the $k$-Bruhat order lead to actions of
$\calF_n$ on $S_n\cup\{\bzero\}$.
The \defcolor{action of $\calF_n$ on $S_n\cup\{\bzero\}$ induced by the
Grassmannian--Bruhat order} is defined on generators $\qv_{ab}$ for $a<b$ in $[n]$ by $\qv_{ab}\cdot\bzero=\bzero$, and for
$\zeta\in S_n$, 
\begin{equation*}
   \qv_{ab} \cdot \zeta\ =\ \begin{cases}
            (a,b) \zeta, & \text{ if }\zeta\precdot (a,b)\zeta,\\
            \bzero,& \text{ otherwise,}
       \end{cases}
\end{equation*}
and extended via composition to $\calF_n$.
For $k< n$, the \defcolor{action of $\calF_n$ on
$S_n\cup\{\bzero\}$ induced by the $k$-Bruhat order} is similarly defined  on generators $\qv_{ab}$ for
$a<b$ in $[n]$ by $\qv_{ab}\kaction\bzero=\bzero$, and for $u\in S_n$,
\begin{equation*}
   \qv_{ab}\kaction u\ =\ \begin{cases}
            (a,b) u, & \text{ if }u \lessdot_k (a,b)u,\\
            \bzero,& \text{ otherwise.}
       \end{cases}
\end{equation*}

Let \defcolor{$\calM_n$} be the monoid with $\bzero$ that is the quotient of $\calF_n\cup\{\bzero\}$ by the
following relations.
\begin{equation}
    \label{Eq:monoid-relations}
    \begin{array}{rrclll}
        (i)   & \qv_{ab}\qv_{cd}           & \equiv & \qv_{cd}\qv_{ab}                     &  & \text{($b<c$ or $a<c<d<b$)} \\ [0.5ex]
        (ii)  & \qv_{ab}\qv_{cd}           & \equiv & \qv_{cd}\qv_{ab}\ \equiv\ \bzero     &  & \text{($a\leq c<b\leq d$)}  \\ [0.5ex]
        (iii) & \qv_{bc}\qv_{cd}\qv_{ac} & \equiv & \qv_{bd}\qv_{ab}\qv_{bc}                 &  & \text{($a<b<c<d$)}          \\ [0.5ex]
        (iv)  & \qv_{ac}\qv_{cd}\qv_{bc} & \equiv & \qv_{bc}\qv_{ab}\qv_{bd}                 &  & \text{($a<b<c<d$)}          \\ [0.5ex]
        (v)   & \qv_{bc}\qv_{ab}\qv_{bc} & \equiv & \qv_{ab}\qv_{bc}\qv_{ab}\ \equiv\ \bzero &  & \text{($a<b<c$)}
    \end{array}
\end{equation}
If a composition  $\qv$ is equivalent to $\bzero$ ($\qv\equiv\bzero$), then we say that \demph{$\qv$ is zero} and
if $\qv$ is not equivalent to $\bzero$, then we say that \demph{$\qv$ is nonzero}.
We summarize the main results of~\cite{BS-Monoid}.

\begin{proposition}\label{Prop:monoid}
    Let $e$ denote the identity element of $S_n$.
    \begin{enumerate}
        \item
            The map $\calF_n\cup\{\bzero\}\to S_n\cup\{\bzero\}$ defined by
            $\bx\mapsto \bx \cdot e$ is surjective, factors through $\calM_n$,
            and induces a bijection between $\calM_n$ and $S_n \cup \{\bzero\}$.
        \item
            For $\zeta\in S_n$, the fiber $\{\bx\in\calF_n\mid \bx \cdot
            e = \zeta\}$ is in bijection with the set of saturated chains in
            $[e,\zeta]_\preceq$.
        \item\label{Prop_3}
            For $u,w\in S_n$ and $k<n$, the set $\{ \bx\in\calF_n\mid \bx\kaction u=w\}$ is
            in bijection with the set of saturated chains in $[u,w]_k$.
    \end{enumerate}
\end{proposition}

\subsubsection{Quantum left operators}
\label{SS:left-quantum-operators}

Unlike what we just described, we do not know of an analog of the Grassmannian--Bruhat
order or of the monoid $\calM_n$ for the quantum Bruhat order on $S_n$.
Some structure does extend, most notably the two- and three-term relations in~\eqref{Eq:monoid-relations} and certain symmetries.
With the goal of an analog of Proposition~\ref{Prop:monoid}\eqref{Prop_3}, we define 
\demph{(left) operators on $S_n[q]\cup\{\bzero\}$} as follows.
For $1\leq k<n$, $\{a,b\}\subseteq [n]$ with  $a\ne b$, and $u\in S_n$ define
\begin{align} \label{Eq:Def-k-action}
    \defcolor{{\qv}_{ab}} \kaction u
    &=
    \begin{cases}
        (a,b) u, \hphantom{q_{i,j}} & \text{if $a<b$ and $u \lessdot_k^q (a,b)u\,,$}            \\[0.4ex]
          q_{i,j} (a,b) u\,, & \text{if $a>b$ and  $u\lessdot_k^q q_{i,j} (a,b) u$\,,\ \ 
           ($u(i)=a$, $u(j)=b$)} \\[0.4ex]
        \bzero\,, & \text{otherwise\,,}
    \end{cases}
\end{align}
and extend to  $S_n[q] \cup \{\bzero\}$ by setting
$\defcolor{\qv_{ab} \kaction (q^\alpha u) \vcentcolon= q^\alpha \, (\qv_{ab} \kaction u)}$.
Note that if ${\qv}_{ab} \kaction u\neq\bzero$, then $u\lessdot_k^q{\qv}_{ab} \kaction u$ is a cover, and all covers in the
quantum Bruhat order occur in this way.

\begin{remark}\label{Rem:k-action}
We make a few remarks and definitions about these operators.

\begin{enumerate}[label=(\roman*)]
  \item
      The operators $\qv_{ab}$ with $a<b$ are those
      of~\S\ref{SS:classical-left-operators}, and are subject to the
      relations in~\eqref{Eq:monoid-relations}.

  \item
        We say that ${\qv}_{ab}$ is \demph{classical} if $a<b$
        and \demph{quantum} if $a>b$.

  \item  \label{Rem:k-action:zero}
    A composition $\qv$ of operators \demph{is zero} ($\defcolor{\qv\equiv \bzero}$) if for every $u\in S_n$
    and $1\leq k<n$, we have that $\qv\kaction u= \bzero$, using~\eqref{Eq:Def-k-action}.

  \item  \label{Rem:k-action:not-zero}
    If we do not have $\qv\equiv\bzero$ in this sense, then $\qv$ is \demph{nonzero}.
    In this case, there exist $u,w\in S_n$, $1\leq k<n$, and $\alpha\in\NN^{n-1}$
    such that $\qv\kaction u= q^\alpha w$.
    Note that $u\leq_k^q q^\alpha w$.

  \item  \label{Rem:k-action:chains}
    In the situation of~\ref{Rem:k-action:not-zero}, where $\qv\kaction u=q^\alpha w$, if
    $\qv=\qv_{a_rb_r}\dotsb\qv_{a_2b_2}\qv_{a_1b_1}$, then there exist $w_0=u,w_1,\dotsc,w_r=w$ in $S_n$ and
    $\alpha_0=0, \alpha_1,\dotsc,\alpha_r=\alpha$ in $\NN^{n-1}$ such that for all $1\leq i<j\leq r$,
    \begin{equation}\label{Eq:sub-chain}
      \qv_{a_jb_j} \dotsb \qv_{a_ib_i} \kaction  (q^{\alpha_{i-1}}w_{i-1})\ =\ q^{\alpha_j} w_j\,.
    \end{equation}
    In particular, the composition $\qv_{a_jb_j} \dotsb \qv_{a_ib_i}$ is nonzero.
    We also have that 
    \[
      u\ \lessdot_k^q\ q^{\alpha_1}w_1\ \lessdot_k^q\ q^{\alpha_2}w_2
      \ \lessdot_k^q\ \dotsb\  \lessdot_k^q\ q^{\alpha_r}w_r=q^\alpha w
    \]
    is a saturated chain in the interval $[u,q^\alpha w]_k^q$.

  \item\label{Rem:k-action:zero_compose}
    A consequence of the definition~\ref{Rem:k-action:zero} and the observation~\ref{Rem:k-action:chains} is that if
    $\qv\equiv\bzero$, then $\qv''\qv\qv'\equiv\bzero$ for any operators $\qv',\qv''$.
    
  \item\label{Rem:k-action:uk_equiv}
    In the situation of~\ref{Rem:k-action:not-zero}, with $\qv\kaction u=q^\alpha w$, suppose that 
    we have an operator $\qv'$ such that $\qv'\kaction u=q^\alpha w$.
    Then we say that the operators $\qv$ and $\qv'$ are \demph{$(u,k)$-equivalent}.
    We explore this equivalence relation in Sections~\ref{SS:Shape-equivalence} through~\ref{S:degreeTwo}. 

    In this case, by the same reasoning as~\ref{Rem:k-action:chains}, $\qv'$ represents a saturated chain in $[u,q^\alpha w]_k^q$.
    When $\qv\neq\qv'$, the operators induce distinct chains.
    As $[u,q^\alpha w]_k^q$ is a ranked poset, $\qv$ and $\qv'$ have the same length, so that $(u,k)$-equivalence is ranked.

  \item
    Suppose that $u\leq_k^q q^\alpha w$ and we have     a saturated chain in $[u,q^\alpha w]_k^q$, 
    \[
      u\ =\ q^{\alpha_0}w_0\ \lessdot_k^q\ q^{\alpha_1}w_1\ \lessdot_k^q\ q^{\alpha_2}w_2
      \ \lessdot_k^q\ \dotsb\  \lessdot_k^q\ q^{\alpha_r}w_r=q^\alpha w\,,
    \]
    For each $i=1,\dotsc,r$, consider an operator $\qv_{a_ib_i}$ satisfying $\qv_{a_ib_i}\kaction q^{\alpha_{i-1}}w_{i-1}=q^{\alpha_i}w_i$.
    Then we have that $(w_{i-1})^{-1}w_i=(s,t)$ with $s\leq k<t$,  $a_i=w_{i-1}(s)$, and $b_i=w_{i-1}(t)$.
    Furthermore, if $a_i<b_i$, then $q^{\alpha_i}=q^{\alpha_{i-1}}$, but if $a_i>b_i$, then  $q^{\alpha_i}=q_{s,t}q^{\alpha_{i-1}}$.
    Hence, the operator $\defcolor{\qv}\vcentcolon=\qv_{a_rb_r}\dotsb \qv_{a_1b_1}$ is such that
    $\qv\kaction u=q^\alpha w$.
    
  \item \label{Rem:k-action:subword}
    Finally, suppose that $\qv''\qv\qv'$ is a nonzero operator and $\qv''\qv\qv'\kaction u=q^\alpha w$.
    If we write $\qv'\kaction u=q^\beta x$ and $\qv\qv'\kaction u=q^\gamma y$, then $\qv\kaction x = q^{\gamma-\beta} y$.
    Suppose that $\qu$ is $(x,k)$-equivalent to $\qv$, so that $\qu\kaction x = q^{\gamma-\beta} y$.
    Then $\qv''\qv\qv'$ is $(u,k)$-equivalent to $\qv''\qu\qv'$.
    
\end{enumerate}

Items~\ref{Rem:k-action:zero_compose} and~\ref{Rem:k-action:subword} show that equivalence of operators may be induced by
subwords.\hfill$\diamond$ 
\end{remark}

A consequence of Remark~\ref{Rem:k-action} is the following generalization of Proposition~\ref{Prop:monoid}\eqref{Prop_3}.

\begin{proposition}\label{Prop:chains_and_operators}
  Suppose that $u\leq_k^q q^\alpha w$.
  Then there is a bijection between saturated chains in the interval $[u,q^\alpha w]_k^q$ and nonzero operators $\qv$ 
  such that $\qv\kaction u=q^\alpha w$, and this set of operators forms a single $(u,k)$-equivalence class of nonzero operators.
\end{proposition}

\subsection{Flattening and $(u,k)$-equivalence}
\label{SS:Shape-equivalence}

Equivalence among classical operators,~\eqref{Eq:monoid-relations} and
Proposition~\ref{Prop:monoid}, depends only upon the relative order of the indices, and not their actual values.
We establish a version for quantum operators.
We begin with some definitions.

Let $a=(a_1,\dotsc,a_t)$ be a sequence of integers and set $m\vcentcolon=\#\{a_1,\dotsc,a_t\}$, the number of distinct
integers in this sequence.
The \demph{flattening} of $a$ is the sequence $(c_1,\dotsc,c_t)$ where each $c_i\in[m]$ and $c_i>c_j$ if and only if $a_i>a_j$.
For example, the flattening of $(3,6,1,6,8,3,1)$ is $(2,3,1,3,4,2,1)$.
Sequences are \demph{shape-equivalent} if they have the same flattening.
The flattening of an operator $\qv_{a_r b_r}\dotsb\qv_{a_1 b_1}$ is  $\qv_{c_r d_r}\dotsb\qv_{c_1 d_1}$, where $(c_r,d_r,\dotsc,c_1,d_1)$  is
the flattening of $(a_r,b_r,\dotsc,a_1,b_1)$.
Operators $\qv$ and $\qv'$ are \demph{shape-equivalent} if they have the same flattening.
This extends the notion of shape-equivalence of permutations $\zeta$ from Section~\ref{SS:GBO}.

The \demph{support} of an operator $\qv=\qv_{a_r b_r}\dotsb\qv_{a_1 b_1}$ is $\supp(\qv)\vcentcolon=\{a_r,b_r,\dotsc,a_1,b_1\}$.
The operator $\qv$ is \demph{minimal} if the factorization of 
$\zeta=\zeta(\qv)=(a_r,b_r)\dotsb(a_1,b_1)$ into transpositions has minimal length in that
$r=\#\supp(\zeta)-s(\zeta)$.
As $r$ is the minimal possible number of transpositions whose product is $\zeta$, we have that
$\supp(\zeta(\qv))=\{a_r,b_r,\dotsc,a_1,b_1\}=\supp(\qv)$.
Furthermore, if $\qv\kaction u=q^\alpha w$, then as $\zeta(\qv)=wu^{-1}$, the interval $[u,q^\alpha w]_k^q$ is minimal in the
sense of Section~\ref{sec:Qminim}.
Later, in the proof of Theorem~\ref{Th:nonZero} in Section~\ref{sec:final}, we will show that the permutation $\zeta(\qv)$ is a
minimal permutation in the sense of Section~\ref{minimal-permutations}.

Shape-equivalence preserves equivalence to $\bzero$ and for minimal operators, $(u,k)$-equivalence.

\begin{theorem}\label{Th:shape-equivalence}
  Let $\qv$ be an operator with flattening $\qv'$.
  Then $\qv\equiv\bzero$ if and only if $\qv'\equiv\bzero$.

  Suppose that $\qv$ is nonzero and minimal, and that $\qu$ is an operator that is $(u,k)$-equivalent to $\qv$
  for some $u\in S_n$ and $k<n$.
  Then $\qu$ is minimal and $\supp(\qv)=\supp(\qu)$.
  Furthermore, if $\qv'$ is the flattening of $\qv$ and $\qu'$ the flattening of $\qu$, then there exist
  $x\in S_m$ and $l<m$, where $m=\#\supp(\qv)$, such that $\qv'$ is $(x,l)$-equivalent to $\qu'$.
\end{theorem}


We deduce Theorem~\ref{Th:shape-equivalence} at the end of this subsection after establishing some preliminary lemmas and definitions.
Let $s\in[n]$.
Define the truncation map $\defcolor{\tau_s}\colon[n]\to[n{-}1]$ by
 \[
   \tau_s(j)\ =\ \left\{ \begin{array}{rcl} j&\ &\mbox{if }j<s\\ j{-}1&&\mbox{if } j\geq s\end{array}\right. \ .
 \]
This extends to operators.
For distinct $a,b,s\in[n]$, define $\defcolor{\tau_s(\qv_{a,b})}\vcentcolon=\qv_{\tau_s(a),\tau_s(b)}$.
If $\qv$ is an operator and $s\not\in\supp(\qv)$, then we may apply $\tau_s$ to every operator $\qv_{a,b}$ in $\qv$, and we
define $\tau_s(\qv)$ to be the composition of the $\tau_s(\qv_{a,b})$.
For example, $\tau_3(\qv_{45}\qv_{51}\qv_{25})=\qv_{34}\qv_{41}\qv_{24}$.

For $u\in S_n$ and $r\in[n]$, let $\defcolor{u/r}\in S_{n-1}$ be the permutation obtained by deleting the $r$th column and
$u(r)$th row of the permutation matrix of $u$.
In one-line notation, this is obtained by omitting the $r$th value of $u$ and then flattening, so that
$(41\underline{3}652)/3=31542$.
Omitting the $s$th component of  $\alpha\in\NN^{n-1}$ gives $\alpha/s\in\NN^{n-2}$, which is the sequence
$(\alpha_1,\dotsc,\alpha_{s-1},\alpha_{s+1},\dotsc,\alpha_{n-1})$.
Thus $(3,1,\underline{2},4,0)/3 = (3,1,4,0)$.

We also have an expansion map $\defcolor{\iota_s}\colon[n{-}1]\to[n]$ defined by
 \[
   \iota_s(j)\ =\ \left\{ \begin{array}{rcl} j&\ &\mbox{if }j<s,\\ j{+}1&&\mbox{if } j\geq s.\end{array}\right. 
 \]
This extends to operators in the same way as $\tau_s$ extended to operators.
For $r,s\in[n{+}1]$ and $v\in S_n$, define a permutation $u=\varepsilon_{r,s}(v)\in S_{n+1}$ by $u(r)=s$ and $u/r=v$.
For $\alpha\in\NN^{n-1}$, define $\varepsilon_{r,s}(\alpha)\vcentcolon=(\alpha_1,\dotsc,\alpha_{r{-}1},s,\alpha_r,\dotsc,\alpha_{n-1})\in\NN^{n}$.

\begin{lemma}\label{L:relative-labels}
  Let $\qv$ be a nonzero operator with $\qv\kaction u=q^\alpha w$.
  \begin{enumerate}[label=(\roman*)]
  \item\label{L:r-l:i}
    For $s\not\in\supp(\qv)$, set $r=u^{-1}(s)$.
    Then $r=w^{-1}(s)$ and $\tau_s(\qv) {\kactionop}_{\tau_r(k)} (u/r) = q^{\alpha/r} w/r$.
  \item \label{L:r-l:ii}
    Suppose that $r\leq u^{-1}(a)$ for all $a\in\supp(\qv)$ or $r>u^{-1}(a)$ for all $a\in\supp(\qv)$.
    Then for any $s\in[n{+}1]$, $\iota_s(\qv)\kactionop_{\iota_r(k)} \varepsilon_{r,s}(u) = q^{\varepsilon_{r,0}(\alpha)}\varepsilon_{r,s}(w)$.
  \end{enumerate}
\end{lemma}

\begin{proof}
  Both statements follow (by composition) from the case when $\qv=\qv_{ab}$ is a single operator.
  As $u\lessdot_k q^\alpha w=\qv_{ab}\kaction u$, we have
  $w=(a,b)u=u(i,j)$, where $i=u^{-1}(a)\leq k < u^{-1}(b)=j$ and either $q^\alpha=1$ when $a<b$ or $q^\alpha=q_{i,j}$ when $a>b$.
  The statements follow from the definitions of $\tau_s$, $\iota_s$, $/r$, $\varepsilon_{r,s}$, and the
  characterizations of covers in Propositions~\ref{nothing-in-between-is-in-between} and~\ref{prop:quantum-k-bruhat-covers}, in a
  case-by-case analysis depending upon the relative positions of $s$ with $a,b$ and $r$ with $i,j$.
  While we omit that analysis, we illustrate three cases in Example~\ref{Ex:three_cases} below.

  Note that when $a<b$, $\qv_{ab}$ is classical, which is treated in~\cite{BS-Bruhat,BS-Monoid}.
\end{proof}

\begin{example}\label{Ex:three_cases}
  Let $u=154\pb32$.
  Then $\qv_{53}\kactionop_3 u = q_2q_3 13452$.
  We illustrate Lemma~\ref{L:relative-labels}\ref{L:r-l:i}.
  Consider each case of $s\in\{1,2,4\}=[5]\smallsetminus\supp(\qv_{53})$.
  \begin{itemize}
    \item  If $s=1$ then $r=1$, $\tau_1(3)=2$, and $u/1=4321$.
     Thus, we have $\qv_{42}\kactionop_2 43\pb 21= q_1q_2 2341$.
    \item  If $s=2$  then $r=5$, $\tau_5(3)=3$, and $u/5=1432$.
     Thus, we have $\qv_{42}\kactionop_3 143\pb 2= q_2 q_3 1234$.
    \item  If  $s=4$  then $r=3$, $\tau_3(3)=2$, and $u/3=1432$.
     Thus, we have $\qv_{43}\kactionop_2 14\pb32= q_2 1342$.
  \end{itemize}

   For Lemma~\ref{L:relative-labels}\ref{L:r-l:ii}, $r$ is either at most 2 or greater than 4.
  \begin{itemize}
     \item $r=1$.  Let $s=4$.  Then $\iota_4(\qv_{53})=\qv_{63}$, $\iota_1(3)=4$, and $\varepsilon_{1,4}(u)=416532$.
   We have $\qv_{63}\kactionop_4 4165\pb 32=q_3q_4 413562$.
    \item   $r=5$.  Let $s=3$.  Then $\iota_3(\qv_{53})=\qv_{64}$, $\iota_5(3)=3$, and $\varepsilon_{5,3}(u)=165432$.
   We have $\qv_{64}\kactionop_3 165\pb 432=q_2q_3 145632$. 
   \hfill$\diamond$
  \end{itemize}
\end{example}

\begin{corollary}\label{C:single_equivalences}
  Let $\qv$ be a composition of operators.
  \begin{enumerate}[label=(\arabic*)]
    
  \item \label{C:SE_1}
    For $s\in[n]\smallsetminus\supp(\qv)$, we have that $\qv\equiv\bzero$ if and only if $\tau_s(\qv)\equiv\bzero$.

  \item
    For $s\in[n{+}1]$, we have that $\qv\equiv\bzero$ if and only if $\iota_s(\qv)\equiv\bzero$.

  \item  \label{C:SE_3}
    Suppose that $\qv$ is nonzero and $\qv'$ is a composition that is $(u,k)$-equivalent to $\qv$ for some $u\in S_n$ and $k<n$. 
    
        \begin{enumerate}[label=(\roman*)]
        \item  \label{C:SE_3_i}
          Let $s\in[n]\smallsetminus(\supp(\qv)\cup\supp(\qv'))$ and set $r\vcentcolon= u^{-1}(s)$.
          Then $\tau_s(\qv)$ and $\tau_s(\qv')$ are $(u/r,\tau_r(k))$-equivalent.

        \item
          For any $s\in[n{+}1]$ and any $r$ such that either $r\leq u^{-1}(a)$ for all $a\in\supp(\qv)\cup\supp(\qv')$ or
          $r>u^{-1}(a)$ for all $a\in\supp(\qv)\cup\supp(\qv')$, we have that 
          $\iota_s(\qv)$ and $\iota_s(\qv')$ are $(\varepsilon_{r,s}(u), \iota_r(k))$-equivalent.

        \end{enumerate}   
  \end{enumerate}
\end{corollary}

\begin{proof}
  The third statement holds by Lemma~\ref{L:relative-labels}.
  This also proves the (contrapositive of) the implication $\Leftarrow$ in the first two statements.
  As in (1), $\qv=\iota_s(\tau_s(\qv))$ and in (2), $s\not\in\supp(\iota_s(\qv))$ and $\qv=\tau_s(\iota_s(\qv))$, the other
  directions also hold.  
\end{proof}

\begin{proof}[Proof of Theorem~\ref{Th:shape-equivalence}]
  As flattening is induced by applying truncation operators, the first statement of Theorem~\ref{Th:shape-equivalence} follows
  from Corollary~\ref{C:single_equivalences}~\ref{C:SE_1}.

  Suppose that $\qv$ is nonzero and minimal.
  Let $\qu$ be an operator that is $(u,k)$-equivalent to $\qv$.
  Then there exist $w\in S_n$ and $\alpha\in\NN^{n-1}$ such that $\qv\kaction u=\qu\kaction u=q^\alpha w$.
  As $\zeta(\qv)=wu^{-1}$, $\zeta(\qu)=\zeta(\qv)$.
  As $(u,k)$-equivalence is ranked, both $\qv$ and $\qu$ involve the same number of elementary operators $\qv_{ab}$.
  This implies that $\qu$ is minimal and that $\supp(\qv)=\supp(\qu)$.

  Consequently, any sequence of truncation operators that flattens $\qv$ also flattens $\qu$.
  The existence of $x$ and $l$ in the final part of Theorem~\ref{Th:shape-equivalence} follows by applying 
  Corollary~\ref{C:single_equivalences}~\ref{C:SE_3}~\ref{C:SE_3_i} along  this same sequence of truncation operators.  
\end{proof}

\subsection{Diagrammatic notation}
\label{SS:diagrammatic-notation}
\begin{figure}[htpb]
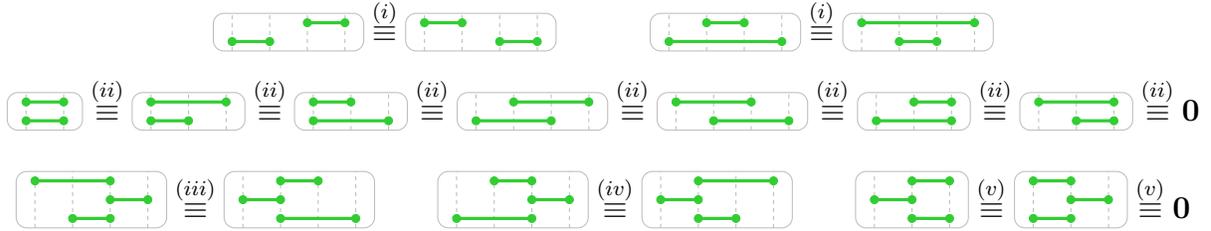

    \begin{gather*}
        \begin{array}{rc@{\hskip3pt}c@{\hskip3pt}c@{\hskip3em}c@{\hskip3pt}c@{\hskip3pt}c}
            & \valuediagram[scale=1.25]{{1,2},{3,4}}
            & \stackrel{(i)}{\equiv}
            & \valuediagram[scale=1.25]{{3,4},{1,2}}
            & \valuediagram[scale=1.25]{{1,4},{2,3}}
            & \stackrel{(i)}{\equiv}
            & \valuediagram[scale=1.25]{{2,3},{1,4}}
        \end{array}
        \\[1ex]
        \begin{array}{rc@{\hskip3pt}c@{\hskip3pt}c@{\hskip3pt}c@{\hskip3pt}c@{\hskip3pt}c@{\hskip3pt}c@{\hskip3pt}c@{\hskip3pt}c@{\hskip3pt}c@{\hskip3pt}c@{\hskip3pt}c@{\hskip3pt}c@{\hskip3pt}c@{\hskip3pt}c}
            & \valuediagram[scale=1.25]{{1,2},{1,2}}
            & \stackrel{(ii)}{\equiv}
            & \valuediagram[scale=1.25]{{1,2},{1,3}}
            & \stackrel{(ii)}{\equiv}
            & \valuediagram[scale=1.25]{{1,3},{1,2}}
            & \stackrel{(ii)}{\equiv}
            & \valuediagram[scale=1.25]{{1,3},{2,4}}
            & \stackrel{(ii)}{\equiv}
            & \valuediagram[scale=1.25]{{2,4},{1,3}}
            & \stackrel{(ii)}{\equiv}
            & \valuediagram[scale=1.25]{{1,3},{2,3}}
            & \stackrel{(ii)}{\equiv}
            & \valuediagram[scale=1.25]{{2,3},{1,3}}
            & \stackrel{(ii)}{\equiv}
            & \bzero
        \end{array}
        \\[2ex]
        \begin{array}{rc@{\hskip3pt}c@{\hskip3pt}c@{\hskip2em}c@{\hskip3pt}c@{\hskip3pt}c@{\hskip2em}c@{\hskip3pt}c@{\hskip3pt}c@{\hskip3pt}c@{\hskip3pt}c}
            & \valuediagram[scale=1.25]{{2,3},{3,4},{1,3}}
            & \stackrel{(iii)}{\equiv}
            & \valuediagram[scale=1.25]{{2,4},{1,2},{2,3}}
            & \valuediagram[scale=1.25]{{1,3},{3,4},{2,3}}
            & \stackrel{(iv)}{\equiv}
            & \valuediagram[scale=1.25]{{2,3},{1,2},{2,4}}
            & \valuediagram[scale=1.25]{{2,3},{1,2},{2,3}}
            & \stackrel{(v)}{\equiv}
            & \valuediagram[scale=1.25]{{1,2},{2,3},{1,2}}
            & \stackrel{(v)}{\equiv}
            & \bzero
        \end{array}
    \end{gather*}
    \caption{Diagrammatic notation for the relations~\eqref{Eq:monoid-relations} among classical operators.}
    \label{fig:classical-relations}
\end{figure}

We introduce a \defcolor{diagrammatic notation} for compositions of operators $\qv_{ab}$, which  will facilitate our study of relations
among them.
Figure~\ref{fig:classical-relations} depicts the relations among classical operators in~\eqref{Eq:monoid-relations} using this notation.
Depict a classical operator $\qv_{ab}$ with $a<b$ by a green line segment $[a,b]$
and a quantum operator $\qv_{ab}$ with $a>b$ by a dotted red line segment $[b,a]$,
drawn in a window with the positions of $1,\dotsc,n$ indicated by vertical
dashed lines as follows,
\newcommand{\labelhack}{\begin{picture}(62,2)
   \put(0 ,-4){{\small 1}}\put(11.5,-4){{\small 2}}
   \put(23,-4){{\small 3}}\put(34.5,-4){{\small 4}}
   \put(46,-4){{\small 5}}\put(57.5,-4){{\small 6}}\end{picture}}
\[
\begin{array}{c@{\qquad}c@{\qquad}c}
    \labelhack& \labelhack& \labelhack   \\
    \valuediagram[padding=1, xmin=1, xmax=6]{{2,4}}
    &
    \valuediagram[padding=1, xmin=1, xmax=6]{{1,5}}
    &
    \valuediagram[padding=1, xmin=1, xmax=6]{{6,3}}
    \\[1ex]
    \qv_{24} & \qv_{15} & \qv_{63}
\end{array}\ .
\]
A composition of operators is depicted by stacking intervals.
For example,
\newcommand{\llabelhack}{\begin{picture}(85,2)
   \put(0 ,-4){{\small 1}}\put(11.5,-4){{\small 2}}
   \put(23,-4){{\small 3}}\put(34.5,-4){{\small 4}}
   \put(46,-4){{\small 5}}\put(57.5,-4){{\small 6}}
   \put(69,-4){{\small 7}}\put(80.5,-4){{\small 8}}\end{picture}}
\[
    \begin{array}{c@{\qquad}c}
    \raisebox{-1pt}{\llabelhack} & \raisebox{2pt}{\llabelhack} \\
     \valuediagram[height=2, padding=1]{{1,2},{4,2},{6,8}}
    &
     \valuediagram[padding=1]{{1,4},{8,4},{5,7},{4,2}}
     \\[1ex]
    \qv_{68}\qv_{42}\qv_{12} & \qv_{42}\qv_{57}\qv_{84}\qv_{14}
  \end{array}\ .
\]
Since the operators act on the left, a composition acts by the rightmost
operator first, which corresponds to reading the associated diagram from the
bottom to the top.

\begin{remark}\label{Rem:relative-labels}
  We may omit labeling the vertical dashed lines and dashed lines not in the support.
  For instance, the diagram $\valuediagram[height=.9, padding=.5]{{2,3},{3,1}}$
  represents any composition of the form ${\qv}_{ca} {\qv}_{bc} $ with $a<b<c$.
%
%
%
  One reason is that equivalence among classical operators,~\eqref{Eq:monoid-relations} and
  Proposition~\ref{Prop:monoid}, only depends upon the relative order of the indices.
  The other reason is that we are primarily concerned with equivalence to $\bzero$ and $(u,k)$-equivalence for minimal operators.
  By Theorem~\ref{Th:shape-equivalence}, these equivalences only depend upon the relative order of the indices.
 \hfill$\diamond$
\end{remark}

\subsection{Symmetries of the quantum $k$-Bruhat order}
\label{sec:actions-on-left-operators}

Operations that preserve the $k$-Bruhat order (its symmetries) played a role in many results in~\cite{BS-Bruhat}.
We identify three operations preserving the quantum $k$-Bruhat order, which are symmetries of that order.
These induce three operations on the operators $\qv_{ab}$ that preserve equivalence to $\bzero$ and $(u,k)$-equivalence, and which are
fundamental to our proof of Theorem~\ref{Th:nonZero}. 
We first define these operations and state the main result of this section before studying the  quantum $k$-Bruhat order.

Let \defcolor{$\ncycle$} be the $n$-cycle $(1, 2, \dotsc, n)$.
For a composition $\qv=\qv_{a_r b_r}\cdots \qv_{a_1 b_1}$, define
\[
    \defcolor{\ncycle(\qv)}\ =\ \ncycle(\qv_{a_r b_r}\cdots \qv_{a_1 b_1})
    \   \vcentcolon=\
    \qv_{\ncycle(a_r),\ncycle(b_r)} \cdots \qv_{\ncycle(a_1),\ncycle(b_1)}\,.
\]
This cyclically shifts indices.
Note that  $\ncycle(\qv_{n,b}) = \qv_{1,b+1}$ and $\ncycle(\qv_{a,n}) = \qv_{a+1,1}$, interchanging
classical and quantum operators, while other applications of $\ncycle$ preserve the type of the operator.
When $n = 3$, we have $\ncycle = (1,2,3)$ and
\[
  \arraycolsep=0.1cm
  \begin{array}{ccccccccccc}
    \valuediagram[padding=1, xmin=1, xmax=3]{{2,3}} &
     \raisebox{-3pt}{$\xrightarrow{\,\ncycle\,}$}&
     \valuediagram[padding=1]{{3,1}} &
     \raisebox{-3pt}{$\xrightarrow{\,\ncycle\,}$}&
     \valuediagram[padding=1, xmax=3]{{1,2}}
    &\raisebox{-2pt}{\quad and\quad}&
     \valuediagram[padding=1, xmin=1, xmax=3]{{3,2}} &
     \raisebox{-3pt}{$\xrightarrow{\,\ncycle\,}$}&
     \valuediagram[padding=1]{{1,3}} &
     \raisebox{-3pt}{$\xrightarrow{\,\ncycle\,}$}&
     \valuediagram[padding=1, xmax=3]{{2,1}}
   \\[1ex]
     {\qv}_{23}&& {\qv}_{31}&& {\qv}_{12}&&
     {\qv}_{32}&& {\qv}_{13}&& {\qv}_{21}
 \end{array}
\]
Each triple cycles back to the original operator with another application of $\ncycle$.
Let \defcolor{$w_0$} be the longest element in $S_n$, so that $w_0(i)=n{+}1{-}i$ for $i\in[n]$.
Define
\[
  \defcolor{w_0(\qv)}\ =\ w_0(\qv_{a_r b_r}\cdots \qv_{a_1 b_1})
  \ \vcentcolon=\
  \qv_{w_0(b_r), w_0(a_r)} \cdots \qv_{w_0(b_1), w_0(a_1)}\,.
\]
This replaces an index $c$ by $w_0(c)=n{+}1{-}c$ and reverses the order of indices on each $\qv_{ab}$.
It preserves quantum and classical operators as $a<b\Leftrightarrow n{+}1{-}b<n{+}1{-}a$.
Finally, define  $\defcolor{\rho(\qv_{a_r b_r}\cdots \qv_{a_1b_1})}\vcentcolon= \qv_{a_1b_1}\cdots\qv_{a_r b_r}$, 
which reverses the order of the operators.

We introduce notation to quantify how these operations act on monomials in the $q_i$.
If $q^\alpha=q_1^{\alpha_1}\dotsb q_{n-1}^{\alpha_{n-1}}$, then set
$\defcolor{w_0(q^\alpha)}\vcentcolon=\defcolor{q^{w_0(\alpha)}}=q_1^{\alpha_{n-1}}\dotsb q_{n-1}^{\alpha_{1}}$.
For $\ncycle$, this requires Laurent monomials in the $q_i$.
When $i<j$, we defined $q_{i,j}$ to be $q_iq_{i+1}\dotsb q_{j-1}$.
When $i>j$, set $q_{i,j}=q_{j,i}^{-1}$.
For $u,w\in S_n$, define $\defcolor{q^{\ncycle(u,w)}}\vcentcolon=q_{w^{-1}(n),u^{-1}(n)}$, which keeps track of how the monomial $q^\alpha$
in elements of $S_n[q]$ is affected by the cyclic shift.
Note that $q^{\ncycle(u,w)}=q_{1,u^{-1}(n)}\cdot q_{w^{-1}(n),1}$, which shows that it is multiplicative in that
$q^{\ncycle(u,w)}=q^{\ncycle(u,v)}\cdot q^{\ncycle(v,w)}$, for any $u,v,w\in S_n$.
Observe that for $u\in S_n$ and $i<j$, if we set $w\vcentcolon= u(i,j)$, then 
\begin{equation}\label{Eq:qOcycle}
   q^{\ncycle(u,w)}\ =\ \left\{ \begin{array}{rcl} 
      q_{i,j}^{-1} &\quad&\mbox{if } u(i)=n\,,\\
      q_{i,j} &\quad&\mbox{if } u(j)=n\,,\\
      1 &\quad&\mbox{otherwise}\,. \end{array}\right.
\end{equation}
For $u\in S_n$, we will often write $\ncycle u$ for $\ncycle(u)$.

\begin{lemma}\label{L:symmetry_qkBO}
  Let $u\in S_n$, $1\leq i\leq k<j\leq n$, and $\alpha\in\NN^{n-1}$.
  Set $\defcolor{w}\vcentcolon= u(i,j)$.
  Then the following are equivalent:
  \begin{align*}
    & (1)\ u\ \lessdot_k^q\ q^\alpha w, 
    & &(3)\ w_0uw_0\ \lessdot_{n-k}^q\ q^{w_0(\alpha)} w_0 w w_0, \\
    & (2)\ \ncycle u\ \lessdot_k^q\ q^{\ncycle(u,w)} q^\alpha \ncycle w, 
    & &(4)\ w w_0\ \lessdot_{n-k}^q\ q^{w_0(\alpha)} u w_0.
 \end{align*}
\end{lemma}

That the operator $\ncycle$ preserves the quantum $k$-Bruhat order in this way follows from results of Postnikov~\cite{Postnikov_symmetry};
for the order it is his Corollary 12 and the exponents  $q^{\ncycle(u,w)}$ are related to his operator $O$ referred to as ``twisted cyclic
shift'' in {\it loc.\ cit.}
We deduce the following corollary concerning the action of the operators $\qv_{ab}$.

\begin{corollary}\label{C:symmetries-q-ops}
  Let $u,w\in S_n$,  $q^\alpha w\in S_n[q]$, and $a\neq b$ in $[n]$.
  The following are equivalent.
  \begin{align*}
    & (1)\ u\ \lessdot_k^q\ \qv_{ab}\kaction u = q^\alpha w,
    & &(3)\ w_0uw_0\ \lessdot_{n-k}^q\ \ w_0(\qv_{ab})\kactionop_{n-k} w_0uw_0 =  q^{w_0(\alpha)} w_0 w w_0,\\
    & (2)\ \ncycle u\ \lessdot_k^q\ \ncycle(\qv_{ab})\kaction \ncycle u = q^{\ncycle(u,w)} q^\alpha \ncycle w,
    & &(4)\ w w_0\ \lessdot_{n-k}^q\ \qv_{ab}\kactionop_{n-k} ww_0 = q^{w_0(\alpha)} u w_0.
 \end{align*}
\end{corollary}

\begin{proof}[Proof of Lemma~\ref{L:symmetry_qkBO}]
 We will freely use Propositions~\ref{nothing-in-between-is-in-between} and~\ref{prop:quantum-k-bruhat-covers} which
 characterize the covers $\lessdot_k$ and $\lessdot_k^q$.
 Note that if $q^\alpha=1$, then $u\lessdot_k u(i,j)=w$ is a classical cover, otherwise
 $u\lessdot_k^q  q_{i,j}u(i,j)=q^\alpha w$, so that $q^\alpha=q_{i,j}$.
 
 Suppose first that  $q^\alpha=1$, so that $u\lessdot_k u(i,j)=w$.
 Thus, if $a=u(i)$ and $b=u(j)$, then $a<b$ and if $i<l<j$ then either $u(l)<a$ or $b<u(l)$.
 As $w_0(c)=n{+}1{-}c$, we have that $w_0u w_0 \lessdot_{n-k} w_0 w w_0$ and  $w w_0 \lessdot_{n-k} u w_0$.
 Thus (1) implies both (3) and (4) when $q^\alpha=1$.

 Still supposing that $q^\alpha=1$, if $u(j)=b<n$, then $\ncycle u(j)=b+1\leq n$ and we see that 
 $\ncycle u \lessdot_k \ncycle w$.
 As $q^{\ncycle(u,w)}=1$, we have $\ncycle u \lessdot^q_k q^{\ncycle(u,w)} \ncycle w$.
 If $n=b=u(j)$, then for every $i<l<j$, $u(l)<a=u(i)$.
 We have that $\ncycle u(i)=a+1>\ncycle u(j)=\ncycle(n)=1$ and if  $i<l<j$, then
 $\ncycle u(i)=a{+}1>\ncycle u(l)>1=\ncycle u(j)$, which implies that $\ncycle u \lessdot^q_k q_{i,j} \ncycle w$.
 By~\eqref{Eq:qOcycle}, $q^{\ncycle(u,w)}=q_{i,j}$ and we again have that $\ncycle u \lessdot^q_k q^{\ncycle(u,w)} \ncycle w$.
 Thus when $q^\alpha=1$, (1) implies (2).
 Similar arguments show that (1) implies (2$'$), which is the statement (2), but where we use $\ncycle^{-1}$ in place of
 $\ncycle$ and replace $q^{\ncycle(u,w)}$ by $q^{\ncycle(\ncycle^{-1}(w),\ncycle^{-1}(u))}$. 

 Suppose now that $q^\alpha\neq 1$, so that $q^\alpha=q_{i,j}$ and $u\lessdot^q_k q_{i,j} w$.
 Then $a=u(i)>b=u(j)$ and if $i<l<j$, then $a>u(l)>b$.
 Multiplying by $w_0$ on the right shows that
 \[
     w w_0\ \lessdot^q_{n-k}\ q_{n+1-j,n+1-i} w w_0 (n{+}1-j,n{+}1-i)\ =\  w_0(q_{i,j})\, u w_0\,,
 \]
 and further multiplying by $w_0$ on the left shows that
 \[
    w_0uw_0\  \lessdot^q_{n-k}\  q_{n+1-j,n+1-i} w_0 u w_0  (n{+}1-j,n{+}1-i)\ =\   w_0(q_{i,j})\, w_0 w w_0\,.
 \]
 Thus (1) implies both (3) and (4) in this case.

 Let us now consider (2).
 If $a\neq n$, then as before, $\ncycle u \lessdot_k^q  q_{i,j}\ncycle w=q^{\ncycle(u,w)} q_{i,j}\ncycle w$, as $q^{\ncycle(u,w)}=1$.
 If $a=n$, then $q^{\ncycle(u,w)}=q_{i,j}^{-1}$~\eqref{Eq:qOcycle}.
 If $i<l<j$, then $n>u(l)>b=u(j)$, and we have $n\geq u(l)+1 = \ncycle u(l)>\ncycle u(j)=b+1$.
 Thus $\ncycle u \lessdot_k \ncycle u(i,j)=\ncycle w= q^{\ncycle(u,w)} q_{i,j} \ncycle w$.
 Thus (1) implies (2).
 Similar arguments show that (1) implies (2$'$), which is the statement (2), but where we use $\ncycle^{-1}$ in place of
 $\ncycle$ and  $q^{\ncycle(\ncycle^{-1}(w),\ncycle^{-1}(u))}$ in place of $q^{\ncycle(u,w)}$. 

 These arguments show that if $u\lessdot^q_k q^\alpha u(i,j)=q^\alpha w$, then (1) implies each of (2), (2$'$), (3), and (4).
 The implications  $(3)\Rightarrow(1)$ and $(4)\Rightarrow(1)$ follow as $u\mapsto u w_0$ and  $u\mapsto w_0 u w_0$ are involutions.
 Similarly, $\ncycle^{-1}\ncycle u=u$ and $q^{\ncycle(u,w)}q^{\ncycle(w,u)}=1$, thus $(1)\Rightarrow(2')$
 for $\ncycle u \lessdot^q_k q^{\ncycle(u,w)} q^\alpha\ncycle w$ is $(2)\Rightarrow(1)$.
 This completes the proof that the four statements are equivalent.
\end{proof}

We combine the equivalences of Corollary~\ref{C:symmetries-q-ops} with the interpretation of equivalent operators
corresponding to different chains in intervals of Remark~\ref{Rem:k-action}~\ref{Rem:k-action:chains}
and~\ref{Rem:k-action:uk_equiv}.

\begin{lemma}    \label{lemma:actions}
    The following statements are equivalent:
    \begin{enumerate}[label=(\alph*)]
      
    \item \label{l:actions:a}
      $\qv$ and $\qv'$ are $(u,k)$-equivalent with $\qv \kaction u = \qv' \kaction u = q^\alpha w$.
      
    \item \label{l:actions:b}
      $\ncycle \qv$ and $\ncycle \qv'$ are $(\ncycle u, k)$-equivalent, with
                     $\qv\kaction \ncycle u=\qv'\kaction \ncycle u = q^{\ncycle(u,w)} q^\alpha \ncycle w$.
    \item \label{l:actions:c}
      $w_0 \qv$ and $w_0\qv'$ are $(w_0uw_0, n{-}k)$-equivalent with
                     $w_0\qv\kactionop_{n-k} w_0uw_0=w_0\qv'\kactionop_{n-k} w_0 uw_0 = q^{w_0(\alpha)} w_0 w w_0$.
    \item \label{l:actions:d}
      $\rho \qv$ and $\rho \qv'$ are $(ww_0, n{-}k)$-equivalent with
                     $\rho\qv\kactionop_{n-k} w w_0=\rho\qv'\kactionop_{n-k} ww_0 = q^{w_0(\alpha)} u w_0$.
    \end{enumerate}
\end{lemma}

Before proving the lemma, we observe
that a composition being nonzero is preserved by $\ncycle$, $w_0$, and $\rho$, and the same for being zero.

\begin{corollary}\label{C:symmetries_bzero}
  A composition $\qv$ is equivalent to zero if and only if one (and hence all) of $\ncycle\qv$, $w_0\qv$, and $\rho\qv$ is equivalent to
  zero. 
\end{corollary}

\begin{proof}[Proof of Lemma~\ref{lemma:actions}]
  Let $u,w\in S_n$, $\alpha\in\NN^{n-1}$, and $\qv,\qv'$ be compositions of left operators such that
  $u<_k^q \qv\kaction u=\qv'\kaction u = q^\alpha w$.
  By Remark~\ref{Rem:k-action}\ref{Rem:k-action:chains},  $\qv$ and $\qv'$ each correspond to a chain in $[u,q^\alpha w]_k^q$.
  We apply the equivalences of Corollary~\ref{C:symmetries-q-ops} to each step of those chains.
  For~\ref{l:actions:c}, we have
  $w_0 uw_0 <_{n-k} w_0\qv\kactionop_{n-k} w_0uw_0 =  w_0\qv'\kactionop_{n-k} w_0uw_0 = q^{w_0(\alpha)}w_0 ww_0$.
  For~\ref{l:actions:d}, note that the chain is reversed and this shows that 
  $w w_0 <_{n-k} \rho\qv\kactionop_{n-k} w w_0 =  \rho\qv'\kactionop_{n-k} w w_0 = q^{w_0(\alpha)}u w_0$.

  The argument for~\ref{l:actions:b} is slightly more involved as the action of $\ncycle$ on monomials depends upon the cover.
  Let $u<_k^q \qv u =  q^{\alpha} w$ correspond to the saturated chain in $[u,q^\alpha w]_k^q$,
  \[
      u\ =\ q^{\alpha_0}w_0\ \lessdot_k^q\ q^{\alpha_1}w_1\ \lessdot_k^q\ q^{\alpha_2}w_2
      \ \lessdot_k^q\ \dotsb\  \lessdot_k^q\ q^{\alpha_r}w_r=q^\alpha w\,.
  \]
  By Lemma~\ref{L:symmetry_qkBO}(2), applying $\ncycle$ to $q^{\alpha_i}w_i\lessdot_k^q q^{\alpha_{i+1}}w_{i+1}$ gives
  $q^{\alpha_i}\ncycle w_i\lessdot_k^q q^{\ncycle(w_i,w_{i+1})}q^{\alpha_{i+1}}\ncycle w_{i+1}$.
  If we set $\defcolor{q^{\beta_i}}\vcentcolon=q^{\ncycle(u,w_i)}q^{\alpha_i}$, then, as 
  $q^{\ncycle(u,w_{i+1})}=q^{\ncycle(u,w_i)}q^{\ncycle(w_i,w_{i+1})}$, we obtain the saturated chain
  \[
      \ncycle u\ =\ q^{\beta_0}\ncycle w_0\ \lessdot_k^q\ q^{\beta_1}\ncycle w_1\ \lessdot_k^q\ q^{\beta_2}\ncycle w_2
      \ \lessdot_k^q\ \dotsb\  \lessdot_k^q\ q^{\beta_r}\ncycle w_r=q^{\ncycle(u,w)}q^\alpha \ncycle w\,
  \]
  in the interval $[\ncycle u,q^{\ncycle(u,w)} q^\alpha \ncycle w]_k^q$, which completes the proof.  
\end{proof}

We have the following corollary.
Here $\textit{op}$ applied to a poset gives the opposite poset, reversing all inequalities.

\begin{corollary}\label{C:equivalences}
  Let $u,w\in S_n$ and $\alpha\in\NN^{n-1}$ be such that $u <_k^q q^\alpha w$.
  Then we have the following isomorphisms of posets,
  \[
    [u,q^\alpha w]_k^q
    \ \simeq\ 
    [\ncycle u, q^{\ncycle(u,w)} q^\alpha \ncycle w]_k^q
     \ \simeq\ 
    [w_0 uw_0,q^{w_0(\alpha)} w_0 w w_0]_{n-k}^q
    \ \simeq\ 
    \textit{op}\left([ww_0,q^{w_0(\alpha)} u w_0]_{n-k}^q\right)\,.
  \]
\end{corollary}

\begin{example}
  Figure~\ref{F:three_ops} illustrates Corollary~\ref{C:equivalences} on the rightmost  interval  of Figure~\ref{F:more_minimal}. 
  \begin{figure}[htb]
  \centering
   \begin{picture}(140,110)(-60,0)
    \put(-35,0){\includegraphics{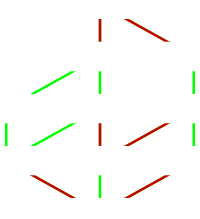}}
                                        \put(-13,100){\small$q_{1,4}q_3132\p45$}  
                                        \put(-10, 75){\small$q_{1,4}135\p42$} \put(42,75){\small$q_3532\p41$}  
    \put(-60,50){\small$q_{1,4}125\p43$} \put(-10, 50){\small$q_{1,4}134\p52$} \put(42,50){\small$q_3531\p42$}  
    \put(-60,25){\small$q_{1,4}124\p53$} \put(0, 25){\small$534\p12$}         \put(42,25){\small$q_3521\p43$}  
                                        \put(0,  0){\small$524\p13$}  
  \end{picture}
   \quad
  \begin{picture}(115,110)(-35,0)
    \put(-25,0){\includegraphics{q-int2.eps}}
                                        \put(-10,100){\small$q_{1,3}23\p541$}  
                                        \put(0, 75){\small$53\p241$}   \put(35,75){\small$q_{1,3}13\p542$}  
    \put(-35,50){\small$43\p251$}       \put(0, 50){\small$52\p341$}   \put(45,50){\small$53\p142$}  
    \put(-35,25){\small$42\p351$}       \put(0, 25){\small$51\p342$}   \put(45,25){\small$43\p152$}  
                                        \put(0,  0){\small$41\p352$}  
  \end{picture}
   \quad
  \begin{picture}(145,110)(-50,0)
    \put(-25,0){\includegraphics{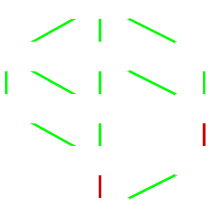}}
                                        \put( 0,100){\small$q_{1,3}25\p314$}  
    \put(-50,75){\small$q_{1,3}24\p315$} \put( 0, 75){\small$q_{1,3}15\p324$} \put(50,75){\small$q_{1,3}23\p514$}  
    \put(-50,50){\small$q_{1,3}23\p415$} \put( 0, 50){\small$q_{1,3}14\p325$} \put(50,50){\small$q_{1,3}13\p524$}  
                                        \put( 0, 25){\small$q_{1,3}13\p425$} \put(60,25){\small$53\p124$}  
                                        \put(10,  0){\small$43\p125$}  
  \end{picture}

  \caption{Actions of $\ncycle$, $\omega_0$, and $\rho$ on the rightmost interval of Figure~\ref{F:more_minimal}.}
  \label{F:three_ops}
  \end{figure}

  These operators act on each chain in the respective orders.
  Each row in Figure~\ref{F:symmetry-chains} shows 
  \begin{figure}[htb]
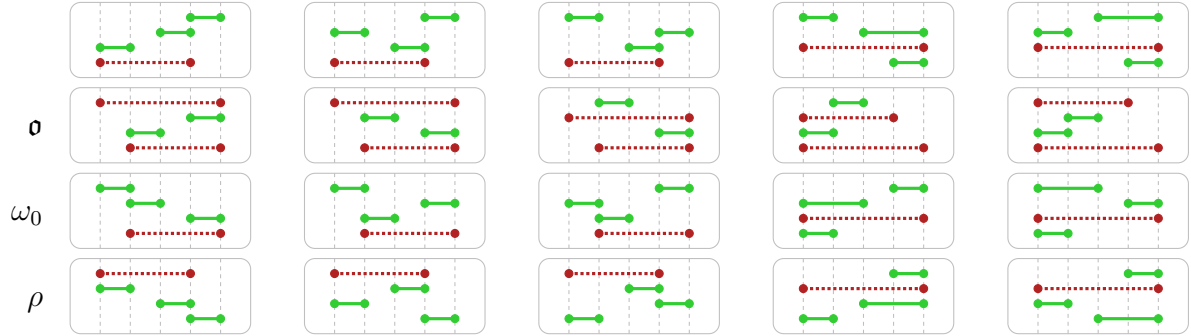

    \centering
    $\begin{array}{rccccccccc}
    &\valuediagram[padding=1]{{4,1},{1,2},{3,4},{4,5}}&\ &\valuediagram[padding=1]{{4,1},{3,4},{1,2},{4,5}}&\ &
    \valuediagram[padding=1]{{4,1},{3,4},{4,5},{1,2}}&\ &\valuediagram[padding=1]{{4,5},{5,1},{3,5},{1,2}}&\ &
    \valuediagram[padding=1]{{4,5},{5,1},{1,2},{3,5}}\\[1ex]
    \raisebox{5pt}{$\ncycle$}
    &\valuediagram[padding=1]{{5,2},{2,3},{4,5},{5,1}}&\ &\valuediagram[padding=1]{{5,2},{4,5},{2,3},{5,1}}&\ &
    \valuediagram[padding=1]{{5,2},{4,5},{5,1},{2,3}}&\ &\valuediagram[padding=1]{{5,1},{1,2},{4,1},{2,3}}&\ &
    \valuediagram[padding=1]{{5,1},{1,2},{2,3},{4,1}}\\[1ex]    
    \raisebox{5pt}{$\omega_0$}
    &\valuediagram[padding=1]{{5,2},{4,5},{2,3},{1,2}}&\ &\valuediagram[padding=1]{{5,2},{2,3},{4,5},{1,2}}&\ &
    \valuediagram[padding=1]{{5,2},{2,3},{1,2},{4,5}}&\ &\valuediagram[padding=1]{{1,2},{5,1},{1,3},{4,5}}&\ &
    \valuediagram[padding=1]{{1,2},{5,1},{4,5},{1,3}}\\[1ex]   
    \raisebox{5pt}{$\rho$}
    &\valuediagram[padding=1]{{4,5},{3,4},{1,2},{4,1}}&\ &\valuediagram[padding=1]{{4,5},{1,2},{3,4},{4,1}}&\ &
    \valuediagram[padding=1]{{1,2},{4,5},{3,4},{4,1}}&\ &\valuediagram[padding=1]{{1,2},{3,5},{5,1},{4,5}}&\ &
    \valuediagram[padding=1]{{3,5},{1,2},{5,1},{4,5}}
    \end{array}$

    \caption{Actions of $\ncycle$, $\omega_0$, and $\rho$ on chains from  rightmost interval of Figure~\ref{F:more_minimal}.}
    \label{F:symmetry-chains}
  \end{figure}
  the diagrams for all five chains in one of the four intervals, with the first row from the rightmost interval of
  Figure~\ref{F:more_minimal} and subsequent rows for the indicated interval in Figure~\ref{F:three_ops}.\hfill$\diamond$
\end{example}

As mentioned in Section~\ref{quantum-equals-classical}, we are often able to reduce to the classical (no quantum operators) case,
which is the origin of our interest in these operators $\ncycle$, $\omega_0$ and $\rho$.

\subsection{Degree two relations}
\label{S:degreeTwo}

Given a composition $\qv$ involving quantum operators, in
Section~\ref{SS:monoid} we developed the notions of equivalence to zero
($\qv\equiv\bzero$),
of $\qv$ being nonzero ($\qv\not\equiv\bzero$),
and the notion of $(u,k)$-equivalence for nonzero compositions $\qv$.
We do not however have a notion of equivalence of operators similar to the stronger one of~\cite{BS-Monoid} for classical
operators.

By Theorem~\ref{Th:shape-equivalence}, a composition $\qv$ being equivalent to $\bzero$ only depends upon its shape-equivalence class.
Furthermore, if $\qv$ is minimal and $\qv\kaction u=q^\alpha w$, then $(u,k)$-equivalence to $\qv$ only depends on
shape-equivalence (and not specifically on $u$ or $k$).
Consequently, we will simply write $\equiv$ for equivalence to $\bzero$ and $(u,k)$-equivalence of minimal compositions.

A composition $\qv_{ab}\qv_{cd}$ is minimal when $\{a,b\}\neq\{c,d\}$.
We identify all relations among minimal compositions that preserve equivalence.
These relations are  depicted in
Figures~\ref{fig:crossing-operators},
~\ref{fig:degree-2-nonzero-relations},
~\ref{fig:more-degree-2-zero-relations},
and
~\ref{fig:touching-orbit}.
This uses the notion of crossing and noncrossing sets from Section~\ref{sec:classic-mn-rule}.

\begin{lemma}\label{Lem:Rel}
 The following describe all relations among minimal compositions of two left operators; those not of the form $\qv_{ab} \qv_{ba}$ or
 $\qv_{ab} \qv_{ab}$. 
 \begin{enumerate}
   \item If $\{a,b\}\cap\{c,d\}=\emptyset$, then the operators commute,  ${\qv}_{ab}{\qv}_{cd} \equiv {\qv}_{cd}{\qv}_{ab}$.
      \begin{enumerate}[itemsep=1ex, topsep=1ex, label=(\roman*)]
   
      \item \label{Rel:1:i}
        If $\{a,b\}$ and $\{c,d\}$  are crossing, then $\qv_{ab}\qv_{cd}  \equiv \bzero$.

      \item  \label{Rel:1:ii}
        If  $a<d<c<b$ or $a>b>c>d$, then  $\qv_{cd}\qv_{ab}$ and $\qv_{ab}\qv_{cd}$ are both zero,
          \[
             \valuediagram[scale=1.1]{{1,4},{3,2}}\equiv
             \valuediagram[scale=1.1]{{2,1},{4,3}}\equiv
             \valuediagram[scale=1.1]{{3,2},{1,4}}\equiv
             \valuediagram[scale=1.1]{{4,3},{2,1}}\equiv \bzero\,.
          \]        
  
        \item  \label{Rel:1:iii}
          In all other cases when $\{a,b\}$ and $\{c,d\}$ are noncrossing, $\qv_{ab}\qv_{cd}\equiv\qv_{cd}\qv_{ab}$, 
         and both are  nonzero.
     \end{enumerate}
        
   \item If $\#(\{a,b\}\cap\{c,d\})=1$, then the operators $\qv_{ab}$ and $\qv_{cd}$ do not commute.
     \begin{enumerate}[itemsep=1ex, topsep=1ex, label=(\roman*)]

     \item  \label{Rel:2:i}
       For $a<b<c$, each of the six compositions,
       $\qv_{ba}\qv_{ac}$,   $\qv_{cb}\qv_{ba}$,   $\qv_{ac}\qv_{cb}$,   $\qv_{ac}\qv_{ba}$,   $\qv_{ba}\qv_{cb}$, and  $\qv_{cb}\qv_{ac}$
       is equivalent to $\bzero$.

     \item   \label{Rel:2:ii}
       For any distinct $a,b,c$, both $\qv_{ab}\qv_{ac}$ and $\qv_{ba}\qv_{ca}$ are  equivalent to $\bzero$.

     \item  \label{Rel:2:iii}
       The six remaining compositions with  $\#(\{a,b\}\cap\{c,d\})=1$ are nonzero and pairwise inequivalent, and are
         depicted in Figure~\ref{fig:nonzero-touching-operators}.
     \end{enumerate}
  \end{enumerate}
\end{lemma}
\begin{remark}
  The composition $\qv_{ab}\qv_{ab}$ is zero.
  When $a<b$, this is by~\eqref{Eq:monoid-relations}({\it ii}).  The case $a>b$ follows from this by cyclic shift (applying $\ncycle$) and 
  Corollary~\ref{C:symmetries_bzero}.

  The composition $\qv_{ab}\qv_{ba}$ is nonzero and $\qv_{ab}\qv_{ba}\not\equiv \qv_{ba}\qv_{ab}$.
  If $a<b$, then $\qv_{ab}\qv_{ba}\kaction u$ equals $q_k u$ if $u(k)=b$ and $u(k{+}1)=a$ and is $\bzero$ otherwise.
  If $a>b$, then $\qv_{ab}\qv_{ba}\kaction u$ equals $q_k u$ if $u(k)=a$ and $u(k{+}1)=b$ and is $\bzero$ otherwise.\hfill$\diamond$
\end{remark}

\begin{proof}[Proof of Lemma~\ref{Lem:Rel}]
  Consider compositions $\qv_{ab}\qv_{cd}$ with $\{a,b\}\cap\{c,d\}=\emptyset$.
  Then $\{a,b\}$ and $\{c,d\}$ form an ordered partition of the four-element set $\{a,b,c,d\}$.
  There are $6=\binom{4}{2}$ such.
  Each operator may be either quantum or classical, there are 24 choices in all.
  
  For (1)\ref{Rel:1:i},
  the eight compositions  $\qv_{ab}\qv_{cd}$ with $\{a,b\}$ and $\{c,d\}$ crossing are shown in Figure~\ref{fig:crossing-operators}.
\begin{figure}[htpb]
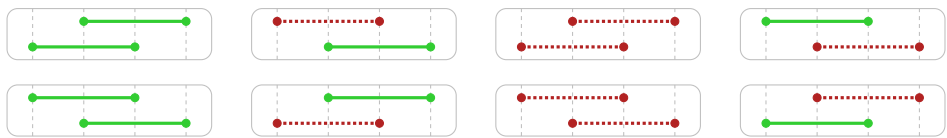

 \centering
  \begin{equation*}
    \begin{array}{l@{\hskip1.25em}l@{\hskip1.25em}l@{\hskip1.25em}l}
        \valuediagram[scale=1.7]{{1,3},{2,4}} & \valuediagram[scale=1.7]{{2,4},{3,1}}
      & \valuediagram[scale=1.7]{{3,1},{4,2}} & \valuediagram[scale=1.7]{{4,2},{1,3}}
            \\[10pt]
        \valuediagram[scale=1.7]{{2,4},{1,3}} & \valuediagram[scale=1.7]{{3,1},{2,4}}
      & \valuediagram[scale=1.7]{{4,2},{3,1}} & \valuediagram[scale=1.7]{{1,3},{4,2}}
    \end{array}
 \end{equation*}
  \caption{Crossing compositions of (1)(i).
   They all are equivalent to $\bzero$.}
 \label{fig:crossing-operators}
\end{figure}
 They form a single orbit under the joint action of $\ncycle$ and $\rho$.
 Both classical compositions (in the first column) are zero by~\eqref{Eq:monoid-relations}({\it ii}).
 Lemma~\ref{lemma:actions} implies the rest are zero.

 The four compositions in (1)\ref{Rel:1:ii} form a single orbit under $\ncycle$ by Lemma~\ref{lemma:actions}, and it suffices to show that one is equivalent to $\bzero$.
 As in Section~\ref{cohomology-ring-of-the-flag-manifold}, 
 write permutations in one-line notation, placing ``\,$\p$\,'' just after the $k$th position.
 Thus $52\pb 6134\in S_6$ has $k=2$,
 and $(\dotsc a\dotsc \pb b \dotsc)$ is a permutation $w$ with $w^{-1}(a)\leq k$ and $w(k{+}1)=b$.

 Let $a<b<c<d$ be the values involved in the first composition in (1)(ii) of operators,
 $\valuediagram{{1,4},{3,2}}=\qv_{cb}\qv_{ad}$, and suppose that $u\in S_n$ is such that $\qv_{cb} \qv_{ad}\kaction u\neq\bzero$.
 Then $u\lessdot_k (a,d)u$, and writing $v \vcentcolon= (a,d)u$, we have that $v$ has the form  $(\dotsc d\dotsc \pb\dotsc a \dotsc)$,
 with neither $b$ nor $c$ appearing between $d$ and $a$.
 Since $\qv_{cb}\kaction v\neq\bzero$, not only must $c$ be left of $\pb$ and $b$ right of $\pb$, but every value $e$ of $v$
 at a position between those of $c$ and $b$ satisfies $b<e<c$.
 But this forces $c$ and $b$ to lie between $a$ and $d$ in $v$, a contradiction.
 Thus $\qv_{cb} \qv_{ad}\equiv\bzero$.

 Consider the remaining twelve $=24-8-4$ noncrossing compositions $\qv_{ab}\qv_{cd}$ in (1).
 If $a<b$ and $c<d$, so that the operators are classical, then by~\eqref{Eq:monoid-relations}({\it i}), they are nonzero
 and $\qv_{ab}\qv_{cd}\equiv\qv_{cd}\qv_{ab}$.
 These are on the left of the two top rows in Figure~\ref{fig:degree-2-nonzero-relations}.
\begin{figure}[htpb]
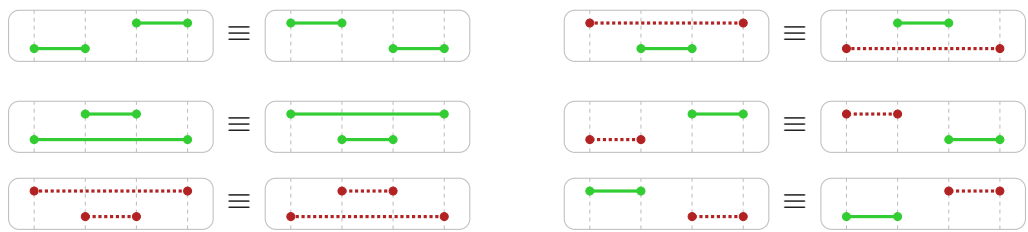

 \[
   \begin{array}{r@{\hskip5pt}c@{\hskip5pt}l@{\hskip3em}r@{\hskip5pt}c@{\hskip5pt}l}
      \valuediagram[scale=1.7]{{1,2},{3,4}} & \raisebox{3pt}{$\equiv$} & \valuediagram[scale=1.7]{{3,4},{1,2}} &
      \valuediagram[scale=1.7]{{2,3},{4,1}} & \raisebox{3pt}{$\equiv$} & \valuediagram[scale=1.7]{{4,1},{2,3}}
      \\[3ex]
      \valuediagram[scale=1.7]{{1,4},{2,3}} & \raisebox{3pt}{$\equiv$} & \valuediagram[scale=1.7]{{2,3},{1,4}} &
      \valuediagram[scale=1.7]{{2,1},{3,4}} & \raisebox{3pt}{$\equiv$} & \valuediagram[scale=1.7]{{3,4},{2,1}}
      \\[2ex]
      \valuediagram[scale=1.7]{{3,2},{4,1}} & \raisebox{3pt}{$\equiv$} & \valuediagram[scale=1.7]{{4,1},{3,2}} &
      \valuediagram[scale=1.7]{{4,3},{1,2}} & \raisebox{3pt}{$\equiv$} & \valuediagram[scale=1.7]{{1,2},{4,3}}
   \end{array}
 \]
 \caption{Relations of Lemma~\ref{Lem:Rel}(1)(iii).}
 \label{fig:degree-2-nonzero-relations}
\end{figure}
 The two relations in the first row are mapped to each other by $\ncycle$,  the four in the remaining rows form an
 orbit under $\ncycle$.
 As equivalence is preserved by $\ncycle$, this establishes (1)(iii).

 For (2),  there are 24 possibilities for $\qv_{ab}\qv_{cd}$:
 As $\#\{a,b,c,d\}=3$, the indices $\{a,b\}$ and $\{c,d\}$ are an ordered choice of two of the three subsets of $\{a,b,c,d\}$ of
 cardinality two. 
 In each of these six, either operator may be quantum or classical, giving 24 choices in all.
 Note thay they do not commute as $(a,b)(c,d)\neq(c,d)(a,b)$ when $\#(\{a,b\}\cap\{c,d\})=1$.

 Figure~\ref{fig:more-degree-2-zero-relations} shows the six compositions in (2)(i).
\begin{figure}[htpb]
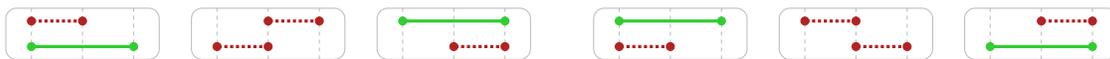

 \[
     \valuediagram[scale=1.7]{{1,3},{2,1}}\quad
     \valuediagram[scale=1.7]{{2,1},{3,2}}\quad
     \valuediagram[scale=1.7]{{3,2},{1,3}}\qquad
     \valuediagram[scale=1.7]{{2,1},{1,3}}\quad
     \valuediagram[scale=1.7]{{3,2},{2,1}}\quad
     \valuediagram[scale=1.7]{{1,3},{3,2}}
 \]
 \caption{The six compositions in (2)(i) are  all equivalent to $\bzero$.}
    \label{fig:more-degree-2-zero-relations}
\end{figure}
 They form an orbit under the joint action of $\ncycle$ and $\rho$.
 Consider the leftmost operator in Figure~\ref{fig:more-degree-2-zero-relations}.
 Writing $a<b<c$ for its support, let $u\in S_n$ be such that $\qv_{ac}\kaction u\neq\bzero$.
 Then $u\lessdot_k \qv_{ac}\kaction u=(a,c)u=\vcentcolon v$, so that $v=(\dotsc c\dotsc \pb\dotsc a \dotsc)$, and $b$ is not between $a$ and
 $c$ in $v$.
 Since  $\qv_{ba}\kaction v\neq\bzero$, $b$ must be left of $\pb$ and hence left of $c$.
 But then $c$ lies between $b$ and $a$, so that $\qv_{ba}\kaction v =\bzero$.
 Hence $\qv_{ba}\qv_{ac}=\bzero$.

 The twelve compositions of (2)(ii) are shown in Figure~\ref{fig:touching-orbit}.
\begin{figure}[htb]
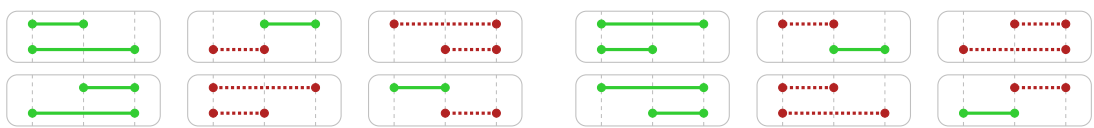

 \centering
  \[
   \begin{array}{ccccccc}
    \valuediagram[scale=1.7]{{1,3},{1,2}}&
    \valuediagram[scale=1.7]{{2,1},{2,3}}&
    \valuediagram[scale=1.7]{{3,2},{3,1}}&&
    \valuediagram[scale=1.7]{{1,2},{1,3}}&
    \valuediagram[scale=1.7]{{2,3},{2,1}}&
    \valuediagram[scale=1.7]{{3,1},{3,2}}\\[1ex]
    \valuediagram[scale=1.7]{{1,3},{2,3}}&
    \valuediagram[scale=1.7]{{2,1},{3,1}}&
    \valuediagram[scale=1.7]{{3,2},{1,2}}&&
    \valuediagram[scale=1.7]{{2,3},{1,3}}&
    \valuediagram[scale=1.7]{{3,1},{2,1}}&
    \valuediagram[scale=1.7]{{1,2},{3,2}}
  \end{array}
 \]
 \caption{The twelve compositions of (2)(ii) are all equivalent to $\bzero$.}
 \label{fig:touching-orbit}
\end{figure}
 They form a single orbit under the joint action of $\ncycle$, $w_0$, and $\rho$.
 As the classical compositions (first and fourth in each row) among them are zero by~\eqref{Eq:monoid-relations}({\it ii}), all twelve are
 zero by Lemma~\ref{lemma:actions}.

 The six remaining compositions in (2)(iii) are shown in Figure~\ref{fig:nonzero-touching-operators}.
 \begin{figure}[htb]
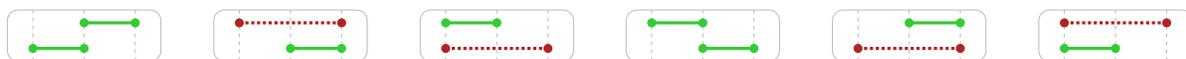

 \centering
     \valuediagram[scale=1.7]{{1,2},{2,3}}\quad
     \valuediagram[scale=1.7]{{2,3},{3,1}}\quad
     \valuediagram[scale=1.7]{{3,1},{1,2}}\quad
     \valuediagram[scale=1.7]{{2,3},{1,2}}\quad
     \valuediagram[scale=1.7]{{3,1},{2,3}}\quad
     \valuediagram[scale=1.7]{{1,2},{3,1}}
  \caption{Remaining compositions of two operators.}
  \label{fig:nonzero-touching-operators}
 \end{figure}
 They form a single orbit under the joint action of $\ncycle$ and $\rho$.
 It suffices to show that one is nonzero and inequivalent to any of the others.
 To that end, note that  $ \qv_{23} \qv_{12}\kactionop_1 1\pb23 = 3\pb12$.
 This represents the first composition in Figure~\ref{fig:nonzero-touching-operators}.
 As $\qv_{23}\kactionop_2 1\pb23 = \bzero$, $\qv_{23} \qv_{12}\not\equiv \qv_{12}\qv_{23}$, so the first and fourth are inequivalent.
 A composition of classical operators cannot be equivalent to a composition involving quantum operators, which completes the proof.
\end{proof}

\subsection{Crossing products}
\label{crossing-products}
For a composition $\qv=\qv_{a_rb_r}\cdots \qv_{a_1b_1}$, let  $\min(\qv)$ and $\max(\qv)$ be the minimum and maximum elements
of its support, $\{a_r,b_r,\dots, a_1,b_1\}$.

A composition $\qv$ defines a \defcolor{(multi)graph} whose vertices are its support.
It has one edge between $a$ and $b$ for each operator $\qv_{ab}$ appearing in $\qv$.
We say that $\qv$ is \demph{connected} if its graph is connected, that $\qv$ is a \demph{path} if its graph is a path with $a_1$ or $b_1$ as
an endpoint and $\#(\{a_i,b_i\}\cap\{a_{i+1},b_{i+1}\})=1$ for all $i$.
It is a \demph{tree} if its graph is tree.
A composition of two connected compositions is \demph{crossing} if their supports are disjoint and crossing.
Note that noncrossing products of trees are minimal.  Indeed, each tree is minimal, so a noncrossing product of trees is a noncrossing product of minimal compositions of operators, and hence is minimal.
Figure~\ref{fig:nonzero-touching-operators} shows all nonzero paths of length two and Figure~\ref{F:symmetry-chains} shows some
compositions that are trees but not paths.  
Figure~\ref{fig:crossing-operators} shows all crossing compositions $\qv_{ab} \qv_{cd}$,
and all are equivalent to zero.

\begin{lemma}
    \label{L:crossing-supports}
    A crossing composition of connected minimal operators is equivalent to zero.
\end{lemma}
\begin{proof}
  Suppose that $\qv$ and $\qv'$ are connected compositions whose supports are disjoint and crossing.
  Then there exists $l<r$ in $\supp(\qv)$ and $l'<r'$ in $\supp(\qv')$ with
  $\{l,r\}$ and $\{l',r'\}$ crossing.
  We may assume that $l < l'$, which implies that  $l < l' < r < r'$.
  As the graph of $\qv$ is connected, it contains a sequence of edges connecting $l$ to $r$.
  Similarly, the graph of $\qv'$  contains a sequence of edges connecting $l'$ to $r'$.

  In the sequence of edges connecting $l$ to $r$, let $(a,b)$ be the first edge with exactly one endpoint between $l'$ and $r'$---this
  corresponds to an operator $\qv_{ab}$ in $\qv$.
  Thus either $l'$ is between $a$ and $b$, or $r'$  is between $a$ and $b$.
  Suppose that $r'$ is between $a$ and $b$.
  Let $(c,d)$ be the first edge in the path from $l'$ to $r'$ with exactly one endpoint between $a$ and $b$---it corresponds to an operator
  $\qv_{cd}$ in $\qv'$.
  Thus $\{a,b\}$ and $\{c,d\}$ are crossing.
  (If  $l'$ is between $a$ and $b$, a similar argument gives an edge $(c,d)$ such that  $\{a,b\}$ and $\{c,d\}$ are crossing.)
  
   By Lemma~\ref{Lem:Rel}(1)(ii), $\qv_{ab} \qv_{cd} \equiv\bzero$.
   As $\supp(\qv)\cap \supp(\qv')=\emptyset$, each operator in $\qv$ commutes with each operator in $\qv'$ by
   Lemma~\ref{Lem:Rel}(1).
   By Remark~\ref{Rem:k-action}~\ref{Rem:k-action:subword}, these commutations give equivalent operators.
   Writing $\qv=\mathbf{A} \qv_{ab} \mathbf{B}$ and $\qv'=\mathbf{C}  \qv_{cd}  \mathbf{D}$, we have
    \begin{equation*}
        \qv \qv'
        \ =\  \big(\mathbf{A} \qv_{ab} \mathbf{B}
          \big)
          \big(\mathbf{C}  \qv_{cd}  \mathbf{D}
          \big)
       \  \equiv\  \mathbf{AC} \big(\qv_{ab} \qv_{cd}\big)\mathbf{BD}\  \equiv\  \bzero\,.
    \end{equation*}
   The first equivalence is by commutation of operators and the second is by Remark~\ref{Rem:k-action}\ref{Rem:k-action:zero_compose}.
\end{proof}

We single out a special case of Lemma~\ref{L:crossing-supports}.

\begin{corollary}
    \label{corollary-classical-product-times-quantum}
    Suppose that $\qu$ is a connected composition of classical operators that is minimal,
    $a, b \notin \supp(\qu)$, $a<b$, and $\qu \qv_{ba} \not\equiv \bzero$.
    If $l = \min(\qu)$ and $r = \max(\qu)$, then
    \begin{gather*}
        a < b < l < r
        \quad\text{or}\quad
        l < r < a < b
        \quad\text{or}\quad
        l < a < b < r
        \quad\text{or}\quad
        a < l < r < b\ .
    \end{gather*}
\end{corollary}
\begin{proof}
    The graphs of $\qu$ and $\qv_{ba}$ are connected, as that of $\qv_{ba}$ is a single edge.
    Since they have disjoint support and $\qu \qv_{ba} \not\equiv \bzero$, Lemma~\ref{L:crossing-supports}
    implies their supports are noncrossing.
\end{proof}

We deduce that a crossing product of a classical tree (no cycle in its graph and no repeated edges) and
a quantum operator is zero.

\begin{corollary}
    \label{T*vq=0}
    Let $\T$ be a composition of classical operators whose graph is a tree.
    If $a<b$ with $\min(\T) < a < \max(\T)$
    and $a \notin \supp(\T)$, then $\T \qv_{ba} \equiv \bzero$.
\end{corollary}
\begin{proof}
  Note that $\T$ is minimal.
  Suppose that  $\T \qv_{ba} \not\equiv \bzero$.
  Then there exist $k<n$ and a permutation $u\in S_n$ such that $\T\qv_{ba}\kaction u\neq\bzero$.
  Suppose that $\qv_{ba}\kaction u = q_{i, j} v$.
  Then $i=u^{-1}(b)\leq k <j=u^{-1}(a)$ and $a=v(i)$.
  By Proposition~\ref{prop:quantum-k-bruhat-covers}, if $i<l<j$, then $a<u(l)=v(l)<b$, and 
  by  Remark~\ref{Rem:k-action}\ref{Rem:k-action:chains}, $\T\kaction v\neq\bzero$.

  Since  $\min(\T)<a<\max(\T)$, $a\not\in\supp(\T)$, and $\T$ is a classical tree, there is an operator $\qv_{cd}$ in $\T$ with $c<a<d$.
  Write $\T=\T'' \qv_{cd} \T'$ and set $w \vcentcolon= \T' \kaction v\neq \bzero$.
  As $\T'$ is classical, $v\leq_k w$.
  By Proposition~\ref{nothing-in-between-is-in-between},
  if $l\leq k$, then $w(l)\geq v(l)$ and we have $a=w(i)$ as $a\not\in\supp(\T')$.

  We also have $\qv_{cd}\kaction w\neq\bzero$, so that $w\lessdot_k (c,d)w$.
  Consider the position $l(\leq k)$ of $c$ in $w$, $c=w(l)$.
  As $c<a<d$ and $d$ is to the right of $k$ (and $i$), we cannot have $l<i$, by Proposition~\ref{nothing-in-between-is-in-between}.
  If $i<l$, then $a<v(l)\leq w(l)=c$, contradicting $c<a$.
  Thus $\T \qv_{ba}\equiv \bzero$.
\end{proof}

\subsection{Paths}
\label{ss:Paths}
We turn our attention to nonzero paths and simply call them paths from now on.
Thus a composition $\qv=\qv_{a_r b_r}\dotsb\qv_{a_1 b_1}$ is a \demph{path} if it is nonzero, its graph is a path with $a_1$ or $b_1$ an
endpoint, and $\#(\{a_i,b_i\}\cap\{a_{i+1},b_{i+1}\})=1$ for all $i$.

\begin{figure}[htb]
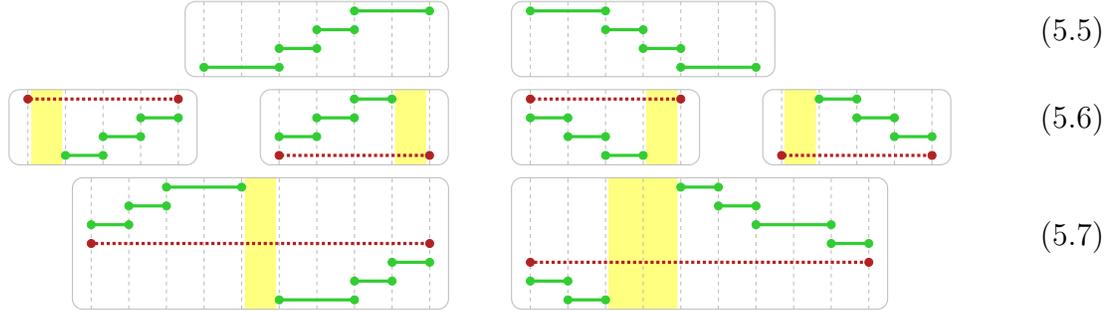

    \begin{gather}
        \valuediagram[scale=1.25, baseline=2ex]{{1,3},{3,4},{4,5},{5,7}}
        \qquad
        \valuediagram[scale=1.25, baseline=2ex]{{5,7},{4,5},{3,4},{1,3}}
        \label{Eq:classicalRowColumn}
        \\
        \valuediagram[scale=1.25, critical windows={(1,2)}, baseline=2ex]{{2,3}, {3,4}, {4,5}, {5,1}}
        \qquad
        \valuediagram[scale=1.25, critical windows={(4,5)}, baseline=2ex]{{5,1}, {1,2}, {2,3}, {3,4}}
        \qquad
        \valuediagram[scale=1.25, critical windows={(4,5)}, baseline=2ex]{{3,4}, {2,3}, {1,2}, {5,1}}
        \qquad
        \valuediagram[scale=1.25, critical windows={(1,2)}, baseline=2ex]{{5,1}, {4,5}, {3,4}, {2,3}}
        \label{Eq:InitialLemma}
        \\
        \valuediagram[scale=1.25, critical windows={(5,6)}, baseline=4ex]{{6,8},{8,9},{9,10},{10,1},{1,2},{2,3},{3,5}}
        \qquad
        \valuediagram[scale=1.25, critical windows={(3,5)}, baseline=4ex]{{2,3},{1,2},{10,1},{9,10},{7,9},{6,7},{5,6}}
        \label{Eq:quantumRowColumn}
    \end{gather}
    \caption{Some nonzero paths.
      Rows are on the left and columns are on the right.}
    \label{fig:paths}
\end{figure}

\begin{figure}[htpb]
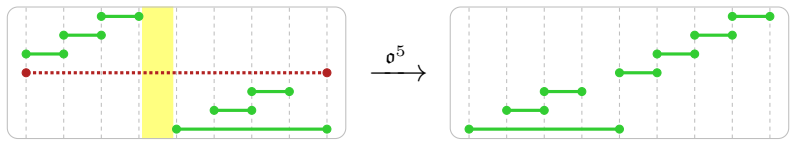

    \centering
    \begin{equation*}
        \begin{array}{c}
            \valuediagram[scale=1.25, critical windows={(4,5)}]{{5,9},{6,7},{7,8},{9,1},{1,2},{2,3},{3,4}}
        \end{array}
        \xrightarrow{~~\ncycle^5~}
        \begin{array}{c}
            \valuediagram[scale=1.25]{{1,5},{2,3},{3,4},{5,6},{6,7},{7,8},{8,9}}
        \end{array}
    \end{equation*}
    \caption{The operator on the right is a \emph{classical row} as it is the product
        of two noncrossing \emph{connected classical rows}.
        The operator on the left is a \emph{row} as it can be cyclically
    shifted to the classical row.}
    \label{fig:column-vs-classical-row}
\end{figure}

\begin{figure}[htpb]
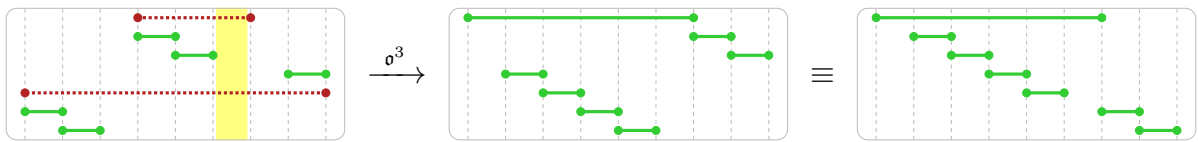

    \centering
    \begin{equation*}
        \begin{array}{c}
            \valuediagram[scale=1.25, critical windows={(6,7)}]{{2,3},{1,2},{9,1},{8,9},{5,6},{4,5},{7,4}}
        \end{array}
        \xrightarrow{~~\ncycle^3~}
        \begin{array}{c}
            \valuediagram[scale=1.25]{{5,6},{4,5},{3,4},{2,3},{8,9},{7,8},{1,7}}
        \end{array}
        \equiv
        \begin{array}{c}
            \valuediagram[scale=1.25]{{8,9},{7,8},{5,6},{4,5},{3,4},{2,3},{1,7}}
        \end{array}
    \end{equation*}
    \caption{The operator on the right is a \emph{classical column} as it is the product
        of two noncrossing \emph{connected classical columns}.
        The operator on the left is a \emph{column} as it can be cyclically
    shifted to the classical column.}
    \label{fig:column-vs-classical-column}
\end{figure}

\begin{remark}\label{rem:yellow} As in Figure~\ref{fig:paths}, a yellow (shaded) area in the diagrammatic rendering of operators highlights
  an open segment that is contained in all quantum operator(s) and no classical operators.
  Such an operator can be cyclically  shifted by powers of $\ncycle$ to a classical operator (see
  Figures~\ref{fig:column-vs-classical-row} and~\ref{fig:column-vs-classical-column}).
  Theorem~\ref{thm:quantum_path} shows this is always the case for paths.\hfill$\diamond$
\end{remark}

A composition $\qv = \qv_{a_r b_r}\dotsb\qv_{a_1 b_1}$  of classical operators
is a \demph{connected classical row} if $b_{i} = a_{i+1}$ for all $1 \leq i < r$.
That is, if $a_1 < a_2 < a_3 < \cdots < a_{r-1} < a_{r} < b_{r}$
so that $\qv=\qv_{a_rb_r}\qv_{a_{r-1} a_r}\dotsb\qv_{a_1a_2}$.
If $u$ is a permutation whose first $r+1$ values are $a_1a_2\dotsc a_rb_r$, then $\qv\kactionop_1 u\neq\bzero$, so $\qv\not\equiv\bzero$.

A \demph{classical row}  $\qv = \qv_{a_r b_r}\dotsb\qv_{a_1 b_1}$ is a noncrossing composition of connected classical rows.
Operators from different connected rows of $\qv$ commute, by~\eqref{Eq:monoid-relations}(i).
By this we mean that commuting them gives an equivalent operator.
Thus, replacing $\qv$ by an equivalent operator, we may assume that $a_1<\dotsb<a_r$, when $\qv$ is a classical row.

A \demph{row} is a composition of operators that can be cyclically shifted to a classical row.
See Figure~\ref{fig:column-vs-classical-row}.
A connected row is a particular instance of a path.

Similarly, a composition $\qv = \qv_{a_r b_r}\dotsb\qv_{a_1 b_1}$ of classical operators is a \demph{connected classical column} if
$b_{i} = a_{i-1}$ for all $1 < i \leq l$.
That is, if $a_r < a_{r-1} < \cdots < a_2 < a_1 < b_1$ so that
$\qv=\qv_{a_ra_{r-1}}\dotsb\qv_{a_2a_1}\qv_{a_1b_1}$.
If $u$ is a permutation whose first $r+1$ values are $a_ra_{r-1}\dotsc a_1b_1$, then $\qv\kactionop_r u\neq\bzero$, so $\qv\not\equiv\bzero$.

A \demph{classical column} $\qv$ is a noncrossing composition of connected classical columns.
As with rows, replacing $\qv$ by an equivalent operator, we may assume that $a_1>\dotsb>a_r$.
A \demph{column} is a composition that can be cyclically shifted to a classical column as
in Figure~\ref{fig:column-vs-classical-column}.
A connected column is an instance of a path.
Note that if $\qv$ is a column, then $\rho(\qv)$ is a row.

\begin{theorem}
   \label{thm:quantum_path}
    Every path is either a row or a column and a path contains at most one quantum operator.
    Paths can be cyclically shifted to a classical row or to a classical column.
    The only path that can be shifted to both consists of a single operator.
\end{theorem}

\begin{proof}
  We split the proof into four parts.
  Let $\qv$ be a path, which is nonzero.\smallskip

\noindent\textit{(1) Paths consisting only of classical operators.}
Suppose the path $\qv = \qv_{a_r b_r}\dotsb\qv_{a_1 b_1}$ is a composition of classical operators (for all $i$, $a_i<b_i$).
Then every pair of adjacent compositions $\qv_{a_{i+1},b_{i+1}}\qv_{a_i,b_i}$ is nonzero.
Suppose that $a_1$ is an endpoint.
As $\qv\not\equiv\bzero$, the relations~\eqref{Eq:monoid-relations}({\it ii}) imply that $a_{i+1}=b_i$ for all $i$.
Similarly, if $b_1$ is an endpoint, then we must have  $a_i=b_{i+1}$ for all $i$.
Thus, $\qv$ is either a connected classical row or a connected classical column.
If  $\qv = \qv_{a_r b_r}\dotsb\qv_{a_1 b_1}$ with for example, $a_i<b_i$ and  $a_{i+1}=b_i$ for all $i$, then
we showed that $\qv \not\equiv \bzero$, and so $\qv$ is a connected classical row.
If instead $a_i=b_{i+1}$, then applying $\rho$ or $w_0$ shows that $\qv$ is a connected classical column.
Such paths are shown in \eqref{Eq:classicalRowColumn} in Figure~\ref{fig:paths}.\smallskip

\noindent\textit{(2) Composition of a connected classical column/row with a quantum operator.}
Suppose that the path $\qv=\qv_{a_ra_{r-1}} \dotsb \qv_{a_2 a_1}\qv_{a_1b} \cdot {\color{BrickRed}\qv_{ba}}$ is a composition
of a connected classical column with single quantum operator ${\color{BrickRed}\qv_{ba}}$.
Then $a_r<a_{r-1}<\dotsb<a_2<a_1<b$, $a<b$, and $a\not\in\{a_r,\dotsc,a_2,a_1\}$.
By Corollary~\ref{T*vq=0}, we must have $a < a_r < \dotsb < a_1 < b$.
Such a path is illustrated in \eqref{Eq:InitialLemma} on the right in Figure~\ref{fig:paths}.

Applying $\ncycle^{n{+}1{-}a_r}$ to $\qv$ cyclically shifts it to a classical column.
We illustrate this.
\begin{equation*}
   \valuediagram[scale=1., critical windows={(1,2)}]{{5,1}, {4,5}, {3,4}, {2,3}}
    \ \raisebox{5pt}{$\xrightarrow{\;\ncycle\;}$}\
     \valuediagram[scale=1., critical windows={(2,3)}]{{1,2}, {5,1}, {4,5}, {3,4}}
    \ \raisebox{5pt}{$\xrightarrow{\;\ncycle\;}$}\
     \valuediagram[scale=1., critical windows={(3,4)}]{{2,3}, {1,2}, {5,1}, {4,5}}
    \ \raisebox{5pt}{$\xrightarrow{\;\ncycle\;}$}\
     \valuediagram[scale=1., critical windows={(4,5)}]{{3,4}, {2,3}, {1,2}, {5,1}}
    \ \raisebox{5pt}{$\xrightarrow{\;\ncycle\;}$}\
     \valuediagram[scale=1.]{{4,5}, {3,4}, {2,3}, {1,2}}
\end{equation*}
If instead $\qv=\qv_{a_ra_{r-1}} \dotsb \qv_{a_2 a_1}\qv_{a_1b} \cdot {\color{BrickRed}\qv_{ba}}$ is the composition of a classical
connected row and a quantum operator ${\color{BrickRed}\qv_{ba}}$, then we may apply $w_0$ to convert it to the composition of a connected
classical column and a quantum operator.
As $w_0\ncycle=\ncycle^{n-1}w_0$, arguments for a column imply this for rows.\smallskip

\noindent\textit{(3) Paths may contain at most one quantum operator.}
Suppose that $\qv$ contains at least two quantum operators.
By Figures~~\ref{fig:more-degree-2-zero-relations}~and~\ref{fig:touching-orbit},
no two quantum operators are adjacent.
Restrict $\qv$ to two of its quantum operators and the classical operators between them.
This segment remains a nonzero path.
Applying $w_0$ if necessary, we may assume that the classical segment is a column,
so that $\qv$ has the form 
 \[
    \qv\ =\ {\color{BrickRed}\qv_{b a_r}}\,\qv_{a_r a_{r-1}}\dotsb \qv_{a_2 a_1}\, {\color{BrickRed}\qv_{a_1 c}} \,,\quad\mbox{where }\
    a_r<b\,, \ a_r<\dotsb<a_2<a_1\,,\ \mbox{ and }c<a_1 \,.
 \]
Applying Part \textit{(2)} to the segment omitting $\qv_{b a_r}$ shows that $c<a_r$.
Applying $\rho$ and then $\omega_0$ to the segment omitting $\qv_{a_1 c}$, Part \textit{(2)}
shows that $a_1<b$.
Applying $\ncycle^{n+1-b}$ to $\qv$ gives:
\[
    \ncycle^{n+1-b}\qv\ =\ \qv_{\alpha_r 1}\,\qv_{\alpha_r \alpha_{l-1}}\dotsb
   \qv_{\alpha_2 \alpha_1}\, {\color{BrickRed}\qv_{\alpha_1,\gamma}}\,,\quad\mbox{where }\
    1<\gamma<\alpha_r<\dotsb<\alpha_2<\alpha_1 \,.
\]
\[
   \valuediagram[scale=1.25, xmin=1, xmax=9]
   {{6,2}, {5,6},{4,5},{3,4}, {8,3}}
   \quad \raisebox{10pt}{$\xrightarrow{\ \ncycle^{n+1-b}\ }$} \quad
   \valuediagram[scale=1.25, xmin=1, xmax=9]{{8,4}, {7,8},{6,7},{5,6}, {1,5}}
\]
As $1<\gamma<\alpha_r$, this is zero by Corollary~\ref{T*vq=0} and Part \textit{(2)}, thus  $\qv$ cannot contain two quantum
operators.

\noindent\textit{(4) Paths containing a single quantum operator.}
Let $\qv$ be a path containing a single quantum operator $\qv_{ba}$ ($a<b$).
Arguing as above, we may suppose $\qv$ that takes the form
\begin{equation*}
    \qv =
    \qv_{b_m b_{m-1}} \cdots \qv_{b_2 b_1} \qv_{b_1 b}
    \, {\color{BrickRed}\qv_{ba}}\, 
    \qv_{a a_{l}} \qv_{a_l a_{l-1}} \dotsb\qv_{a_2 a_1}
\end{equation*}
with $a < a_l < \dotsb < a_{2} < a_1$ and $b_{m} < \cdots < b_{1} < b$.
If $a_1<b_m$, then applying $\ncycle^{n{+}1{-}b_m}$ to $\qv$ cyclically shifts it to a classical column, which completes
the proof in this case.
This is illustrated on the right in~\eqref{Eq:quantumRowColumn} in Figure~\ref{fig:paths}.

Suppose that $b_{m}\leq a_1$.
Since the graph of $\qv$ is a path, $b_{m}\not\in\{a_1,\dotsc,a_r\}$.
We have $a<b_{m}$ by Part \textit{(2)}, so there is some $i$ with $a_{i+1}<b_{m}<a_i$ (if $i=l$, then $a=a_{l+1}$).
Then the cyclic shift $\ncycle^{n-b_{m}}\qv$ is a path with two quantum operators
coming from the shifts of $\qv_{a_{i+1} a_i}$ and $\qv_{b_m b_{m-1}}$,
as illustrated by
\[
  \valuediagram[baseline=2ex]  {{3,5},{2,3},{1,2}, {7,1}, {6,7},{4,6}}
   \xrightarrow{\;\ncycle\;}
  \valuediagram[baseline=2ex]  {{4,6},{3,4},{2,3},{1,2}, {7,1}, {5,7}}
   \xrightarrow{\;\ncycle\;}
  \valuediagram[baseline=2ex]  {{5,7},{4,5},{3,4},{2,3},{1,2},  {6,1}}
   \xrightarrow{\;\ncycle\;}
  \valuediagram[baseline=2ex]  {{6,1}, {5,6},{4,5},{3,4},{2,3}, {7,2}}\ .
\]
But this is a contradiction, as a nonzero path has at most two quantum operators.
\end{proof}

\subsection{Products of trees}
\label{SS:trees}
A nonzero composition $\qv=\T_k \cdots \T_2 \T_1$ is a \demph{forest} if each $\T_i$ is a tree and each pair
of trees has noncrossing supports.
A forest is minimal, so that each pair of trees commute, and commuting them gives an equivalent operator.

\begin{theorem}
    \label{main-theorem}
    Suppose that $\qv$ is a forest.
    Let $u\in S_n$ and $k<n$ be such that $\qv\kaction u\neq\bzero$.
    Then $\qv$ is $(u,k)$-equivalent to a composition $\R\C$ of a row $\R$ and a column $\C$.
    Moreover, there exists $r \in \NN$ such that $\ncycle^r(\R\C)=\ncycle^r(\R)\ncycle^r(\C)$ is classical,
    and $\ncycle^r(\R)\ncycle^r(\C)\kaction \ncycle^r(u)\neq\bzero$.
    
\end{theorem}

The proof is split into two parts.

\subsubsection{The classical case}

The case of Theorem~\ref{main-theorem} in which all the operators are
classical follows from the study of minimal permutations in~\cite[Sect.~6.2]{BS-iso}, and the discussion in
Section~\ref{minimal-permutations}.
Recall that $\calL$ is the rank function in the Grassmannian--Bruhat order on $S_n$.

\begin{proposition}
    \label{BS2002}
    If $\qv$ is a composition of classical operators that is a forest, then it is equivalent to a composition $\R\C$ of a 
    classical row and a classical column. 
\end{proposition}
\begin{proof}
  Suppose that $\zeta=(a_r,b_r)\dotsb(a_1,b_1)$ is the product of transpositions whose corresponding graph is a forest.
  Using that this graph is a forest, induction on $r$ shows that $\supp(\zeta)=\{a_r,b_r,\dotsc,a_1,b_1\}$.
  ({\it A priori}, it is only a subset.)
  The same induction shows that $\zeta$ is a product of disjoint cycles, with each cycle corresponding to a tree in the graph of
  $\zeta$. 

  Now suppose that $\qv=\qv_{a_r,b_r}\dotsb\qv_{a_1,b_1}$ is a nonzero tree with corresponding permutation
  $\zeta\vcentcolon=(a_r,b_r)\dotsb(a_1,b_1)$, which is a cycle--so $s(\zeta)=1$--and which satisfies $\calL(\zeta)=r$.
  Since $\qv$ is a tree, $\supp(\qv)=\{a_r,b_r,\dotsc,a_1,b_1\}$, which equals $\supp(\zeta)$,  has cardinality $r{+}1$.
  Thus $\calL(\zeta)=\#\supp(\zeta)-s(\zeta)$, and we conclude that  $\zeta$ is a minimal cycle.

  When $\qv$ is a forest, its corresponding permutation $\zeta$ is a noncrossing product of minimal cycles, and is therefore a
  minimal permutation.
  By Proposition~\ref{P:peaklessChains}, the interval $[e,\zeta]_\preceq$ has a peakless chain.
  Under the identification of chains in  $[e,\zeta]_\preceq$ with compositions equivalent to $\qv$ from
  Proposition~\ref{Prop:monoid}({\it 2}), this peakless chain corresponds to a composition $\R\C$ of a classical row and a
  classical column.
\end{proof}

\subsection*{Example}
We illustrate Proposition~\ref{BS2002} with a  sequence of equivalences on a classical tree.
\[
    \begin{array}{c}
        \valuediagram[scale=1, padding=1]{{3, 4}, {1, 3}, {5, 6}, {3, 5}, {5, 7}, {2, 5}}
    \end{array}
    \equiv
    \begin{array}{c}
        \valuediagram[scale=1, xmax=7, padding=1]{{5, 6},{3, 4},{2, 3},{1, 2},{2, 7},{3, 5}}
    \end{array}
    \equiv
    \begin{array}{c}
        \valuediagram[scale=1, xmax=7, padding=1]{{2,7},{3,5}} \\
        \valuediagram[scale=1, xmax=7, padding=1]{{5,6},{3,4},{2,3},{1,2}}
    \end{array}
    \equiv
    \begin{array}{c}
        \valuediagram[scale=1, xmax=7, padding=1]{{3,5}} \\[-0.75ex]
        \valuediagram[scale=1, xmax=7, padding=1]{{2,7}} \\
        \valuediagram[scale=1, xmax=7, padding=1]{{3,4},{2,3},{1,2}} \\[-0.75ex]
        \valuediagram[scale=1, xmax=7, padding=1]{{5,6}}
    \end{array}
\]
%

\subsubsection{Proof of Theorem~\ref{main-theorem}}

We use induction on the number of operators in $\qv$.
Let $\qv_{ef}$ be its rightmost operator so that $\qv=\F \qv_{ef}$, where $\F$ is a forest.
As $\qv\kaction u$ is nonzero, $\qv_{ef}\kaction u$ is nonzero.
Then there exist $\alpha\in\NN^{n-1}$ and $v\in S_n$ such that $\qv_{ef}\kaction u=q^\alpha v$, and
$\F\kaction v$ is nonzero.
By induction, $\F$ is $(v,k)$-equivalent to $\R\C$, where $\R$ is a row and $\C$ is a column, and
there exists $r \in \NN$ such that $\ncycle^r(\F)$ is $(\ncycle^r v, k)$-equivalent to $\ncycle^r(\R)\ncycle^r(\C)$, where the later is
classical. 
Thus
\begin{equation*}
    \ncycle^r(\qv) \text{ is } (\ncycle^r u,k)\text{- equivalent to } \ncycle^r(\R)\ncycle^r(\C)\qv_{cd}\,,
\end{equation*}
where $\ncycle^r(\R)$ is a classical row, $\ncycle^r(\C)$ is a classical column, and $\qv_{cd}=\ncycle^r(\qv_{ef})$.

If  $c<d$, then $\ncycle^r(\R)\ncycle^r(\C)\qv_{cd}$ is a classical forest and by
Proposition~\ref{BS2002}, it is $(\ncycle^r u, k)$-equivalent to the composition of a 
classical row and a classical column, $\R' \C'$.
But then $\qv$ is $(u, k)$-equivalent to $\ncycle^{-r}(\R')\ncycle^{-r}(\C')$, which is a composition of a row and a column that admits a cyclic
shift that is classical. 

Suppose now that $c>d$.
Set $a=d$ and $b=c$ so that $\qv_{cd}=\qv_{ba}$ with $a<b$, to conform to our running notation for quantum operators.
Since $\F$ is a forest,  $\ncycle^r(\R)\ncycle^r(\C)$ is a classical forest.
Grouping components of $\ncycle^r(\R)\ncycle^r(\C)$ together gives  noncrossing (classical) trees $\T_1,\dotsc,\T_m$ such that
$\ncycle^r(\qv)$ is $(\ncycle^r u, k)$-equivalent to $(\T_m\dotsb\T_1)\qv_{ba}$.
The proof of Theorem~\ref{main-theorem} follows from Lemma~\ref{cyclic-shift-classical-trees-times-quantum} below,
and since $\qv$ is a forest, an application of Proposition~\ref{BS2002}.

\begin{lemma}
    \label{cyclic-shift-classical-trees-times-quantum}
    Let $a<b$  and suppose that
    \begin{enumerate}[label={{\rm(\roman*)}}]
        \item
            $\T_1, \ldots, \T_m$ are disjoint, noncrossing, classical trees;

        \item
            $\T_m \cdots \T_1 \qv_{ba} \not\equiv \bzero$; and 

        \item
            the graph of $\T_m \cdots \T_1 \qv_{ba}$ is a forest.
    \end{enumerate}
    Then $\T_m \cdots \T_1 \qv_{ba}$
    admits a cyclic shift that is a product of classical operators.
\end{lemma}

\begin{proof}
  We have three cases according to  $\#(\{a,b\}\cap \supp(\T_m \cdots \T_1))$.

  Suppose that   $\emptyset=\{a,b\}\cap \supp(\T_m \cdots \T_1)$.
  Since for each $i \in [m]$, we have $\T_i \qv_{ba}\not\equiv \bzero$, Corollary~\ref{T*vq=0} implies that 
  either $\max(\T_i)<a$ or $a<\min(\T_i)$, for each $i\in[m]$.
  Hence, $\ncycle^{-a}(\T_m \cdots \T_1 \qv_{ba})$ is a composition of classical operators.
  We present an example.
  \begin{equation*}
    \begin{array}{c}
      \valuediagram[scale=1, critical windows={(5,6)}]{
                    {12,5},
                    {10,11}, {9,10},
                    {1,3}, {3,4}, {2,3},
                    {14,15}, {15,16}, {13,15},
                    {6,7}, {7,8}}
    \end{array}
        \xrightarrow{~~\ncycle^{11}~}
    \begin{array}{c}
      \valuediagram[scale=1]{
                    {7,16},
                    {5,6}, {4,5},
                    {12,14}, {14,15}, {13,14},
                    {9,10}, {10,11}, {8,10},
                    {1,2}, {2,3}}
     \end{array}
  \end{equation*}

 Suppose that  $1=\#(\{a,b\}\cap \supp(\T_m \cdots \T_1))$.
 Applying $\omega_0$ if necessary, we may assume that $b\in\supp(\T_m \cdots \T_1)$. 
 As the trees have disjoint support, $b$ lies in the support of exactly one tree.
 As the trees commute, we may replace this operator by an equivalent operator (obtained by commuting the trees) and assume that
 $b \in \supp(\T_m)$.
 By Corollary~\ref{T*vq=0}, $a < \min(\T_m)$ for otherwise $\T_m \qv_{ba} \equiv \bzero$, contradicting that 
 $\T_m \cdots \T_1\qv_{ba}$ is nonzero.

 Since $\emptyset=\{a,b\}\cap\supp(\T_{m-1}\dotsb\T_1)$, the previous argument implies that 
 for $i\in[m{-}1]$ either $\max(\T_i)<a$ or $a<\min(\T_i)$.
 Thus $\ncycle^{-a}(\T_m \cdots \T_1 \qv_{ba})$ is a composition of classical operators.
 We present an example.
 \begin{equation*}
     \begin{array}{c}
         \valuediagram[scale=1, critical windows={(5,6)}]{
             {11,5},
             {7,8}, {6,7},
             {1,3}, {3,4}, {2,3},
             {14,15}, {15,16}, {13,15},
             {9,11}, {11,12}, {10,11}}
     \end{array}
       \xrightarrow{~~\ncycle^{11}~}
     \begin{array}{c}
         \valuediagram[scale=1]{
             {6,16},
             {2,3},{1,2},
             {12,14},{14,15},{13,14},
             {9,10},{10,11},{8,10},
             {4,6},{6,7},{5,6}}
     \end{array}
 \end{equation*}

  Finally, suppose that  $\{a,b\}\subset \supp(\T_m \cdots \T_1)$.
  As the trees are disjoint,  each of $a$ and $b$ lies in the support of a unique tree.
  These cannot be the same tree, for if $a,b\in\supp(\T_i)$, then the graph of $\T_i \qv_{ba}$ is not a tree, 
  contradicting that $\T_m \cdots \T_1 \qv_{ba}$ is a forest.

  As the $\T_i$ commute, we may replace this operator by an equivalent operator and assume that $a \in \supp(\T_{m})$ and
  $b \in \supp(\T_{m-1})$.
  Let $i\in[m{-}1]$.
  As $\T_i v_{ba}\not\equiv\bzero$, Corollary~\ref{T*vq=0} implies that either $\max(\T_i)<a$ or $a<\min(\T_i)$.
  If  $\max(\T_i)<a$, then as $a\in\supp(\T_m)$, we have $\max(\T_i)<\max(\T_m)$.
  If $a<\min(\T_i)$, then as $a\in\supp(\T_m)$ and the trees $\T_i$ and $\T_m$ are noncrossing, we have
  $\max(\T_m)<\min(\T_i)$.

  Since $a\leq \max(\T_m) < b$, this implies that 
  $\ncycle^{-\max(\T_{m})}(\T_m \cdots \T_1 \qv_{ba})$ is a composition of classical operators, which completes the
  proof.
  We present an example.
 \begin{equation*}
     \begin{array}{c}
         \valuediagram[scale=1, critical windows={(7,8)}]{
             {11,6},
             {2,3},{1,2},{2,4},
             {14,16},{13,14},{14,15},
             {10,11},{9,10},{10,12},{8,9},
             {5,6},{6,7}}
     \end{array}
       \xrightarrow{~~\ncycle^{-7}~}
     \begin{array}{c}
         \valuediagram[scale=1]{
             {4,14},
             {10,11},{9,10},{10,12},
             {7,9},{6,7},{7,8},
             {3,4},{2,3},{3,5},{1,2},
             {13,14},{14,15}}
     \end{array}
     \qedhere
 \end{equation*}
\end{proof}

\subsection{Proof of Theorem~\ref{Th:nonZero}} 
\label{sec:final}

Let $[u,q^\alpha w]^q_k$ be a minimal interval, and consider a composition $\qv$ corresponding to some saturated chain in
$[u,q^\alpha w]^q_k$.
Then $\qv\kaction u= q^\alpha w$.
As explained in Section~\ref{sec:Qminim}, as $[u,q^\alpha w]^q_k$ is a minimal interval,
$\ell(q^\alpha w)-\ell(u)=\#\supp(wu^{-1}) - s(wu^{-1})$.
This implies that the graph of $\qv$ contains no cycles, and is therefore a forest.
As $\qv$ is nonzero, it is a forest.

By Theorem~\ref{main-theorem}, $\qv$ is $(u,k)$-equivalent to a composition of operators of the form
$\R\C$ where $\R$ is a row and $\C$ is a column, and there exists a
suitable cycle shift $\ncycle^r$ such that applying it to  $[u,q^\alpha w]^q_k$ gives a classical interval with
a classical chain $\ncycle^r(\R)\ncycle^r(\C)$ that is peakless.
Let $a$ be its height and $a+b-1=\ell(q^\alpha w)-\ell(u)$ its length.
 
By Postnikov's cyclic symmetry~\cite[Thm.~4]{Postnikov_symmetry} (and using the notation of~\cite{Postnikov_symmetry}),
\[
   C_{v(\lambda,k),u,w_0w}\ =\ q^{\beta(r,u,w)}C_{v(\lambda,k),\ncycle^{r} u,\ncycle^{-r}w_0w}\ =\ 
   q^{\beta(r,u,w)}C_{\ncycle^{r} u,v(\lambda,k),w_0\ncycle^{r}w}\,,
\]
where $q^{\beta(r,u,w)}$ is some monomial in the $q_i$'s determined by $u$, $w$ and the power $r$ in $\ncycle^r$.
Translating the notation of~\cite{Postnikov_symmetry} to our notation gives
\[
C_{u,v,w_0w}\ =\ \sum_\alpha q^\alpha N_{u,v}^{w,\alpha}=\ \sum_\alpha C_{u,v}^{q^\alpha w} \,.
\]

Thus if $\lambda=(b,1^{a-1})$ is the hook partition, then the peakless chain $\ncycle^r(\R)\ncycle^r(\C)$ in the classical interval 
$[\ncycle^{r} u,\ncycle^{r}w]_k$ implies that
\[
   C_{u,v(\lambda,k)}^{q^\alpha w}\ =\ C_{\ncycle^{r} u,v(\lambda,k)}^{ q^{\beta(r,u,w)}q^\alpha \ncycle^{r}w}\ =\
   C_{\ncycle^{r} u,v(\lambda,k)}^{ \ncycle^{r}w} \ne 0\,,
\]
by Proposition~\ref{P:peaklessChains} and Lemma~\ref{P:classical_hook}.
This completes the proof.\hfill$\Box$


\subsection*{Acknowledgments}
The authors want to thank the Algebraic Combinatorics Working Seminar at the Fields Institute where this project started in
2016.
Benedetti thanks Grant FAPA of the Faculty of Science of Universidad de Los Andes.
Bergeron is partially supported by NSERC and the York Research Chair in Applied Algebra.
Colmenarejo was partially supported by the AMS-Simons Travel grant.
Saliola is partially supported by NSERC.
Sottile is partially supported by the Simons Foundation and the National Science Foundation through grant DMS-2201005.

\bibliographystyle{amsplain}
\bibliography{references}
\end{document}